\documentclass[11pt]{sat}

\newif\ifsattoc\sattoctrue
\newread\testfl\immediate\openin\testfl=\jobname.toc
    \ifeof\testfl\sattocfalse\fi\closein\testfl

\newcommand{\sect}[1]{\section{#1}\setcounter{equation}{0}}

\title{
Mixed Moduli of Smoothness in $L_p$, $1<p<\infty$: \\A Survey
\footnote{This research was partially supported by the MTM 2011-27637,  RFFI 12-01-00170, NSH 979-2012-1, 2009 SGR 1303.}
}
\def\shorttitle{
Mixed Moduli of Smoothness in $L_p$, $1<p<\infty$: A Survey}

\author{M. K. Potapov, B. V. Simonov, and S. Yu. Tikhonov
}

\def\shortauthor{M. Potapov, B. Simonov, and S. Tikhonov}

\def\versiondate{11 April 2013}

\def\abstracttext{
 In this paper we survey recent developments
over the last 25 years on the
mixed fractional moduli of smoothness of periodic functions from $L_p$, $1<p<\infty$.
In particular, the paper includes  monotonicity properties, equivalence and realization results,
sharp Jackson, Marchaud, and Ul'yanov inequalities,
interrelations between the moduli of smoothness,
the Fourier coefficients, and
``angular" approximation.
The sharpness of the
results presented is discussed.
}

\def\MSCnumbers{ 26A15,  26A33, 42A10, 41A25, 41A30, 42B05, 26B05} % see http://www.ams.org/msc/

\def\keywords{Mixed fractional moduli of smoothness, fractional K-functionals,
``angular" approximation, Fourier sums, Fourier coefficients, sharp Jackson, Marchaud, and Ul'yanov inequalities}

\usepackage{amssymb}
\usepackage{amsfonts}
\usepackage{amsmath}
\usepackage{graphicx}

\RequirePackage[cp1251]{inputenc}

\newcommand\R{\mathbb{R}}

\newcommand{\T}{\mathbb{{T}}}
\newcommand{\N}{\mathbb{{N}}}

\renewcommand\ll{\lesssim}
\renewcommand\gg{\gtrsim}

\newtheorem{lemma}{Lemma}[section]
\newtheorem{theorem}{Theorem}[section]
\newtheorem{remark}{Remark}[section]

\def\startpagenumber{1}
\def\volumenumber{8} % = current year - 2013
\def\year{2013}

\def\dd{\,{\rm d}}  % for integration (making them mathop places them less well)
\def\ee{{\rm e}}  % for the base of the natural log
\def\ii{{\rm i}}  % for the imaginary unit
\def\floor#1{\lfloor#1\rfloor}

\newcommand{\beginddoc}{
\maketitle
\begin{abstract}
\abstracttext
\vskip1pt MSC: \MSCnumbers
\ifx\keywords\empty\else\vskip1pt Keywords: \keywords\fi
\end{abstract}
\insert\footins{\scriptsize
\medskip
\baselineskip 8pt
\leftline{Surveys in Approximation Theory}
\leftline{Volume \volumenumber, \year.
pp.~\thepage--\pageref{endpage}.}
\leftline{\copyright\ \year\ Surveys in Approximation Theory.}
\leftline{ISSN 1555-578X}
\leftline{All rights of reproduction in any form reserved.}
\smallskip
\par\allowbreak}
\ifsattoc\else\tableofcontents\fi}
\renewcommand\rightmark{\ifodd\thepage{\it \hfill\shorttitle\hfill}\else {\it \hfill\shortauthor\hfill}\fi}
\def\endddoc{\label{endpage}\end{document}}
\date{{\small \versiondate}}
\RequirePackage{hyperref}
\RequirePackage{srcltx}

\begin{document}
\beginddoc
\ifsattoc
\bigskip

\vskip 10mm
{\bf Contents} %\newpage

%%%%%%%%%%%%%%%%%%%%%%%%%%%  toc material
\def\toczer{0}\def\tochalf{.5}\def\tocone{1}
\def\tocindent{0}
\def\ection{section}\def\ubsection{subsection}
\def\numberline#1{\hskip\tocindent truecm{} #1\hskip1em}
\newread\testfl
\def\inputifthere#1{\immediate\openin\testfl=#1
    \ifeof\testfl\message{(#1 does not yet exist)}
    \else\input#1\fi\closein\testfl}
\countdef\counter=255
\def\diamondleaders{\global\advance\counter by 1
  \ifodd\counter \kern-10pt \fi
  \leaders\hbox to 15pt{\ifodd\counter \kern13pt \else\kern3pt \fi
  \hss.\hss}\hfill}
\newdimen\lextent
\newtoks\writestuff
\medskip
\begingroup
\small
\def\contentsline#1#2#3#4{
\def\argu{#1}
\ifx\argu\ection\let\tocindent\toczer\else
\ifx\argu\ubsection\let\tocindent\tochalf\else\let\tocindent\tocone\fi\fi
\setbox1=\hbox{#2}\ifnum\wd1>\lextent\lextent\wd1\fi}
\lextent0pt\inputifthere{\jobname.toc}\advance\lextent by 2em\relax
\def\contentsline#1#2#3#4{
\def\argu{#1}
\ifx\argu\ection\let\tocindent\toczer\else
\ifx\argu\ubsection\let\tocindent\tochalf\else\let\tocindent\tocone\fi\fi
\writestuff={#2}
\centerline{\hbox to \lextent{\rm\the\writestuff%
\ifx\empty#3\else\diamondleaders{}
\hfil\hbox to 2 em\fi{\hss#3}}}}
\inputifthere{\jobname.toc}\endgroup
\immediate\openout\testfl=\jobname.toc % to empty \jobname.toc in order to
\immediate\closeout\testfl             % get \tableofcontents to initiate
\renewcommand{\contentsname}{}         % regeneration the toc file without
\tableofcontents%\newpage               % printing a ToC.
\fi

\vskip 10mm

\sect{Introduction}
\label{sec-intr}

To open the discussion on mixed moduli of smoothness,
we  start with
 function spaces of dominating mixed smoothness.
 The Sobolev spaces of dominating mixed smoothness  were first introduced (on $\mathbb{R}^2$) by Nikol'skii \cite{1, 1111}.
He defined the space
\begin{eqnarray*}
S^{r_1,r_2}_pW({\bf R}^2)=\Bigg\{
f\in L_p({\bf R}^2): \|f\|_{S^{r_1,r_2}_pW({\bf R}^2)}&=&\|f\|_{L_p({\bf R}^2)}\\
+
\Big\|\frac{\partial^{r_1 f}}{\partial x_1^{r_1} } \Big\|_{L_p({\bf R}^2)}
&+&
\Big\|\frac{\partial^{r_2} f}{\partial x_2^{r_2} } \Big\|_{L_p({\bf R}^2)}+
\Big\|\frac{\partial^{r_1+r_2} f}{\partial x_1^{r_1}\partial x_2^{r_2} } \Big\|_{L_p({\bf R}^2)}
\Bigg\},
\end{eqnarray*}
where $1<p<\infty$, $r_1, r_2=0,1,2$.
Here, the mixed derivative $\frac{\partial^{r_1+r_2} f}{\partial x_1^{r_1}\partial x_2^{r_2} }$
plays a dominant role and it gave the name to these scales of function spaces.

Later, the
fractional Sobolev spaces with dominating mixed smoothness
(see  \cite{li2} by Lizorkin and Nikol'skii), the H\"{o}lder-Zygmund-type spaces (see  Nikol'skii  \cite{1,1111} and  Bakhvalov \cite{2}),
and
 the Besov spaces of dominating mixed smoothness were introduced (see  Amanov  \cite{Am}).
We would also like to mention the paper \cite{ba} by Babenko which
considered Sobolev spaces with dominating mixed smoothness in the
context of multivariate approximation. It transpires that spaces
with dominating mixed smoothness have several unique properties
which can be used in different settings, for example, in
multivariate approximation theory of periodic functions (see
\cite[1.3]{sch1} and \cite{teml3}) or in high-dimensional
approximation and computational mathematics (see, e.g.,
\cite{stro}).

 To define H\"{o}lder-Besov spaces (Nikol'skii-Besov)
  of dominating mixed smoothness, the notion  of the mixed modulus of smoothness is used, i.e.,
 $$\omega_{\bf{k}}(f, {\bf{t}})_p=
 \omega_{k_1,  \ldots,k_d}(f, t_1,  \ldots,t_d)_p=
 \sup_{|h_i|\le t_i, i=1, \ldots,d}\|\Delta^{\bf{k}}_{\bf{h}} f\|_p,
 $$
 where the $\bf{k}$-mixed difference is given by
 $$
 \Delta^{\bf{k}}_{\bf{h}}=
 \Delta^{{k_1}}_{{h_1}}\circ \cdots \circ
 \Delta^{{k_d}}_{{h_d}},
 $$
$$ {\bf{k}} = (k_1, \ldots, k_d), \quad {\bf{h}} = (h_1, \ldots, h_d),$$
 and $\Delta^{{k_i}}_{{h_i}}$ is the difference of order $k_i$ with step $h_i$ with respect to $x_i$; for example,
 $$\Delta^{{1}}_{{h_i}} f (x_1, \ldots,x_d) = f (x_1, \ldots, x_i+h_i, \ldots,x_d)-f (x_1, \ldots, x_i, \ldots,x_d),$$
$$ \Delta^{{k_i}}_{{h_i}} f (x_1, \ldots,x_d) = \sum_{j=0}^{k_i} (-1)^{j} \binom{k_i}{j}
f (x_1, \ldots, x_i+(k_i-j)h_i, \ldots,x_d), \qquad k_i\in\N.
$$

In their turn, moduli of smoothness of an integer order  can be naturally extended to
fractional order moduli.
The one-dimensional fractional  modulus of smoothness was introduced in the 1970's  (see \cite{butzer, 10, 11} and the monograph \cite{samko}).
 Moreover,  moduli of smoothness of { positive
orders} play an important role in Fourier analysis, approximation theory, theory of embedding theorems and some other problems (see, e.g., \cite{samko, sim-jat, real1, real3, jfaa, trebels, 11}).
 Note that one of the key results in this area of research ---
an equivalence between the modulus of smoothness and the $K$-functional ---
was proved  for the one-dimensional fractional modulus in \cite{butzer} and
for the multivariate  (non-mixed) fractional modulus in \cite{wilmes, w1}; see also \cite{dah, johnen, dah1}.
Clearly, mixed moduli of smoothness are closely related to mixed directional derivatives. Inequalities between the mixed and directional derivatives are given in, e.g., \cite{chen}.

In general, the modulus of smoothness is an important concept in
modern analysis and there are many sources providing information on
the one dimensional and multivariate (non-mixed) moduli from
different perspectives. Concerning mixed moduli of smoothness, two
old monographs \cite{Am} and \cite{timan-book} can be mentioned
where several basic properties are listed. Also, there is vast
literature on the theory of function spaces with dominating mixed
smoothness (see Section \ref{1.3} below). The main goal of this
paper is to collect the main properties of the mixed moduli of
smoothness of periodic functions from $L_p(\T^d)$, $1<p<\infty$,
from the point of view of approximation theory and Fourier analysis.

\vskip 0.5cm

This paper attempts to give a self-contained development of the
theory. Since many of the sources where the reader can find these
results are difficult to obtain and many of the results are stated
without proofs, we will provide complete proofs of all the main
results in this survey. Moreover, there are several results in
Sections 4-11 that are new, to the best of our knowledge.

For the sake of clarity, in  this survey we deal with periodic
functions on $\T^2$. We limit ourselves to this case to help the
reader follow the discussion and to put the notation and results in
a more compact form. All the results of this survey can be extended
to the case of $\T^d$, $d>2$.

Let us also mention that since any $L_p$-function $f$ on $\T^2$ can be written as
$$
f(x,y)= F(x,y)+\phi(x)+\psi(y)+c,
$$
where
$F\in L_p^0(\T^2)$, i.e.,
 $\int_\T F \dd x  = \int_\T F \dd y=0$ and since
 $\omega_{\alpha_1,\alpha_2}(f; \delta_1,\delta_2)_p=\omega_{\alpha_1,\alpha_2}(F; \delta_1,\delta_2)_p$,
 it suffices to deal with functions from $L_p^0(\T^2)$.

\vskip 0.5cm

\subsection{How this survey is organized}
After auxiliary results and notation given in Sections 2 and 3,
in Section 4 we collect the main properties of the mixed moduli, mainly, various monotonicity properties
and  direct and inverse type approximation theorems.
In Section 5 we prove a constructive characterization of the mixed moduli of smoothness which is a  realization result (see, e.g., \cite{ditzian}).
This result provides us with a useful tool to obtain the results of the later sections.
In particular, this allows us to show the equivalence between the mixed modulus of
smoothness  and the corresponding $K$-functional in Section 6.

In Section 7 two-sided estimates of the mixed moduli of smoothness in terms of the Fourier coefficients are given.
In Section 8 we deal with sharp inequalities between the mixed moduli of smoothness of functions and their derivatives, i.e.,
$\omega_{\bf{k}}(f, {\bf{t}})_p$ and $\omega_{\bf{l}}(f^{\bf{r}}, {\bf{t}})_p$.
Section 9  gives sharp order two-sided estimates of the mixed moduli of smoothness of $L_p$ functions in terms of their ``angular" approximations.

In Section 10 we study sharp relationships between
$\omega_{\bf{k}}(f, {\bf{t}})_p$ and $\omega_{\bf{l}}(f,
{\bf{t}})_p$. One part of this relation is usually called the sharp
Marchaud inequality (see, e.g., \cite{dai, march, ima}), another is
equivalent to the sharp Jackson inequality (\cite{jat-dai, dai}). It
is well known that these results are closely connected to the
results of Section 8 because of Jackson and Bernstein-Stechkin type
inequalities. Finally, in Section 11, we discuss sharp Ul'yanov's
inequality, i.e., sharp relationships between $\omega_{\bf{k}}(f,
{\bf{t}})_p$ and $\omega_{\bf{l}}(f, {\bf{t}})_q$ for $p<q$ (see
\cite{sim-jat}).

In Sections 7-10, we deal with two-sided estimates for the mixed
moduli of smoothness. In order to show sharpness of these estimates
we will introduce special function classes so that for functions
from these classes the two-sided estimates become equivalences.

\vskip 0.5cm

\subsection{What is not included in this survey}\label{1.3}
In this paper, we restrict ourselves to questions which were not covered by previous expository papers and which were  actively developed
over the last 25 years. For example, we do not discuss questions which are quite naturally linked to the mixed moduli of smoothness such as
\begin{enumerate}
\item[$\cdot$] Different types of convergence of multiple Fourier series
 (see Chapter I in the surveys \cite{4, zhizhi-s} and the papers \cite{dav, dyach1996});
\item[$\cdot$] Absolute convergence of multiple Fourier series (see Chapter X in the surveys
\cite{4, zhizhi-s} and the papers \cite{mus, veres});
\item[$\cdot$] Summability theory of multiple Fourier series (see the paper \cite{timan1} and the monograph \cite{zhizhi});
\item[$\cdot$]
Interrelations between the total, partial and mixed moduli of smoothness; derivatives
(see \cite{brud, dit, kolyada, deriv-st, timan2});
\item[$\cdot$]Fourier coefficients of functions from certain smooth spaces (see, e.g., \cite{ant, berisha, 4});
\item[$\cdot$]
Conjugate multiple Fourier series (see, e.g., the book \cite{zhizhi-b} and Chapter VIII in the survey \cite{4});
\item[$\cdot$]
Representation and approximation of multivariate functions
(see, e.g.,  \cite{bazar, dav, de1, de2, ram, sun} and Chapter 11 of the recent book \cite{trigub}); in particular,
for  Whitney type results see \cite{dung};
\item[$\cdot$]
Approximate characteristics of functions, entropy, and widths
 (see \cite{han, pust, teml0, teml1,teml2, wang}).
\end{enumerate}
Also, we do not deal with the questions of
\begin{enumerate}
\item[$\cdot$]
The theory of function spaces with dominating mixed smoothness,
\end{enumerate}
in particular, with characterization, representation, embeddings theorems,  characterization of approximation spaces, $m$-term approximation,
which are fast growing topics nowadays.
Let us only mention a few basic older papers \cite{li1, li2, li3}, the monograph \cite{besov} by Besov, Il'in and Nikol'skii,
the monograph by Schmeisser and Triebel \cite{book1}, the recent book by Triebel \cite{book2}, and the 2006's survey
\cite{sch1} on this topic. The reader might also be interested in the recent work of researchers from the Jena school  \cite{ha, krbec, sch, sch3, tri1, vyb, vyb1}; see also  \cite{ho2, ho1}.
 The coincidence of the Fourier-analytic  definition of the spaces of dominating mixed smoothness and the definition in terms of differences is given in the paper \cite{ullrich}.

\vskip 1.5cm
\sect{Definitions and notation}

Let $L_p = L_p(\T^2)$, $1 < p < \infty$, be the space of measurable functions $f$ of two variables that
are $2\pi$-periodic in each variable and  such that
$$
\|f\|_{L_p(\T^2)}=\left(\int\limits_{0}^{2\pi}\int\limits^{2\pi}_{0}|f(x,y)|^p\dd x \dd y\right)^{1/p}<\infty.
$$
Let also $L_p^0(\T^2)$ be the collection of $f\in L_p(\T^2)$ such that
 $$\int\limits^{2\pi}_{0}f(x,y)\dd y=0  \qquad \qquad   \textnormal{for \,\, a.e.}~x$$
 and
 $$\int\limits^{2\pi}_{0}f(x,y)\dd x =0 \qquad \qquad   \textnormal{for \,\, a.e.}~y.$$

If $F(f,\delta_1,\delta_2)>0$ and  $G(f,\delta_1,\delta_2)>0$ for all $\delta_1, \delta_2>0$,
then writing
$F(f,\delta_1,\delta_2)\ll G(f,\delta_1,\delta_2)$ means
 that there exists a constant $C$, independent of $f,\delta_1,\delta_2$ such that $F(f,\delta_1,\delta_2)\leq C G(f,\delta_1,\delta_2).$
 Note that $C$ may depend on unessential parameters (clear from context), and may change form line
to line.
If  $F(f,\delta_1,\delta_2)\ll G(f,\delta_1,\delta_2)$ and $G(f,\delta_1,\delta_2)\ll F(f,\delta_1,\delta_2)$
  simultaneously, then we will write  $F(f,\delta_1,\delta_2) \asymp G(f,\delta_1,\delta_2) $.

\vskip 0.5cm

\subsection{The best  angular approximation}
By  $s_{m_1,\infty}(f),$ $s_{\infty, m_2}(f)$, and  $s_{m_1, m_2}(f)$ we denote the partial sums
of the Fourier series of a function $f\in L^p(\T^2)$, i.e.,
 $$
 s_{m_1, \infty}(f)=\frac1{\pi}\int\limits^{2\pi}_{0}f(x+t_1,y)D_{m_1}(t_1)\dd t_1,
 $$
 $$
 s_{\infty, m_2}(f)=\frac1{\pi}\int\limits^{2\pi}_{0}f(x,y+t_2)D_{m_2}(t_2)\dd t_2,
 $$
 $$
 s_{m_1, m_2}(f)=\frac1{\pi^2}\int\limits^{2\pi}_{0}\int\limits^{2\pi}_{0}f(x+t_1,y+t_2)D_{m_1}(t_1)D_{m_2}(t_2)\dd t_1\dd t_2,
 $$
where $D_m$ is the Dirichlet kernel, i.e.,
 $$D_m(t)=\frac{\sin{(m+\frac12)t}}{2\sin{\frac{t}{2}}}, \qquad  m=0,1,2, \ldots $$

As a means of approximating a function $f\in L^p(\T^2)$, we will use
the so called  best (two-dimensional) angular approximation
$Y_{m_1,m_2}(f)_{L_p(\T^2)}$ which is also sometimes called
``approximation by an angle'' (\cite{6}). By definition,
 $$
 Y_{m_1,m_2}(f)_{L_p(\T^2)}=\inf\limits_{T_{m_1,\infty},T_{\infty,m_2}}\|f-T_{m_1,\infty}-T_{\infty,m_2}\|_{L_p(\T^2)},
 $$
where the function $T_{m_1,\infty} \in L_p(\T^2)$ is a trigonometric polynomial of degree at most $m_1$ in $x$, and the
function $T_{\infty,m_2}\in L_p(\T^2)$ is a trigonometric polynomial of degree at most $m_2$ in $y$.

\vskip 0.5cm

\subsection{The mixed moduli of smoothness}
For a function $f\in L_p(\T^2)$, the difference of order $\alpha_1>0$ with respect to the variable $x$
and the difference of order $\alpha_2>0$ with respect to the variable $y$  are defined  as follows:
$$
\Delta^{\alpha_1}_{h_1}(f)=\sum\limits^{\infty}_{\nu_1=0}(-1)^{\nu_1}\left(^{\alpha_1}_{\nu_1}\right)f(x+(\alpha_1-\nu_1)h_1,y)
$$
and, respectively,
$$
\Delta^{\alpha_2}_{h_2}(f)=\sum\limits^{\infty}_{\nu_2=0}(-1)^{\nu_2}\left(^{\alpha_2}_{\nu_2}\right)f(x,y+(\alpha_2-\nu_2)h_2),
$$
where
$\left(^{\alpha}_{\nu}\right)=1$  for  $\nu=0, \; \left(^{\alpha}_{\nu}\right)=\alpha$  for  $\nu=1 $, $\left(^{\alpha}_{\nu}\right)=\frac{\alpha(\alpha-1)...(\alpha-\nu+1)}{\nu!}$   for  $\nu\geq 2$.

Denote by
$\omega_{\alpha_1,\alpha_2}(f,\delta_1,\delta_2)_{L_p(\T^2)} $
the mixed modulus of smoothness of a function $f \in L_p(\T^2)$ of orders
$\alpha_1>0$ and  $\alpha_2>0$
with respect to the variables $x$ and $y$, respectively, i.e.,
$$
\omega_{\alpha_1,\alpha_2}(f,\delta_1,\delta_2)_{L_p(\T^2)}=\sup\limits_{|h_i|\leq \delta_i, i=1,2}\|\Delta_{h_1}^{\alpha_1}(\Delta_{h_2}^{\alpha_2}(f))\|_{L_p(\T^2)}.
$$
We remark that
$
\|\Delta_{h_1}^{\alpha_1}(\Delta_{h_2}^{\alpha_2}(f))\|_{L_p(\T^2)}\le C(\alpha_1,\alpha_2) \|f\|_{L_p(\T^2)},
$
where $ C(\alpha_1,\alpha_2)\le 2^{\floor{\alpha_1}+\floor{\alpha_2}+2}.$

\vskip 0.5cm

\subsection{$K$-functional}
First, let us recall the definition of the fractional integral and fractional derivative in the sense of Weyl of a function $f$ defined on $\mathbb{T}$.
If the Fourier series of a function  $f\in L^1(\mathbb{T})$
  is given by
$$\sum_{n\in \mathbb{Z}} c_n \ee^{\ii nx}, \qquad c_0=0,$$
then the fractional integral or order $\rho > 0$ of $f$ is defined by (see, e.g., \cite[Ch.~XII]{18})
$$
I^\alpha f(x):= \big(f*\psi_\rho\big)(x)=\frac{1}{2\pi}\int_0^{2\pi} f(t) \psi_\rho(x-t) \dd t,
$$
where
$$
\psi_\alpha(x)=\sum_{{}^{n\in \mathbb{Z}}_{n\ne 0}} \frac{\ee^{\ii nx}}{(\ii n)^\rho}.
$$
To define the fractional derivative or order $\rho > 0$ of $f$, we put
$n:=\floor{\rho}+1$ and
$$
f^{(\rho)}(x):= \frac{\dd^n}{\dd x ^n} I^{n-\rho} f(x).
$$

By $f^{(\rho_1,\rho_2)}$  we will denote the Weyl derivative  of order $\rho_1\ge 0$ with
respect to $x$ and of order $\rho_2\ge 0$ with respect to $y$ of the function
$f \in L_1^0(\T^2)$.

Denote  by
$W_p^{(\alpha_1,0)}$ the Weyl class, i.e., the set of functions  $f \in L_p^0(\T^2)$ such that
 $f^{(\alpha_1,0)} \in L_p^0(\T^2)$. Similarly,
$ W_p^{(0, \alpha_2)}$ is the set of functions  $f \in L_p^0(\T^2)$ such that $f^{(0, \alpha_2)} \in
L_p^0(\T^2).$ Moreover,
$ W_p^{(\alpha_1,\alpha_2)}$ is the set of functions  $f \in L_p^0(\T^2)$ such that
$f^{(\alpha_1,\alpha_2)} \in L_p^0(\T^2).$

The mixed
$K$-functional of a function  $f \in L_p^0(\T^2)$ is given by
\begin{eqnarray*} %K(f, t_1, t_2) &=&
 K(f, t_1, t_2, \alpha_1, \alpha_2, p) &=&
\inf\limits_{g_1 \in W_p^{(\alpha_1,0)},g_2 \in W_p^{(0,
\alpha_2)}, g \in W_p^{(\alpha_1,\alpha_2)}} \Big[ \|f - g_1
- g_2- g \|_{L_p(\T^2)} \\&+& t_1^{\alpha_1} \|g_1^{(\alpha_1,0)}
\|_{L_p(\T^2)} + t_2^{\alpha_2} \|g_2^{(0, \alpha_2)} \|_{L_p(\T^2)}
+ t_1^{\alpha_1} t_2^{\alpha_2} \|g^{(\alpha_1,\alpha_2)}
\|_{L_p(\T^2)} \Big].
\end{eqnarray*}

\vskip 0.5cm

\subsection{Special classes of functions}
We  define the function class  $M_p, 1<p<\infty$, as the set of functions $f\in L_p^0(\T^2)$ such that the Fourier series of $f$ is given by
$\sum\limits^{\infty}_{\nu_1=1}\sum\limits^{\infty}_{\nu_2=1}a_{\nu_1,\nu_2}\cos{\nu_1
x}\cos{\nu_2y}$, where
\begin{equation}\label{fgg}
a_{\nu_1,\nu_2} -a_{\nu_1+1,\nu_2}-
a_{\nu_1,\nu_2+1}+a_{\nu_1+1,\nu_2+1}\geq 0
\end{equation}
for any integers $\nu_1$ and $\nu_2$. Note that  (\ref{fgg}) implies
\begin{equation}\label{mon}
a_{n, m_1}\ge a_{n, m_2} \qquad m_1 \le m_2\quad\mbox{and}\quad
a_{n_1, m}\ge a_{n_2, m} \qquad n_1 \le n_2.
\end{equation}

We also define the function class  $\Lambda_p, 1<p<\infty$, as the set of functions $f\in L_p^0(\T^2)$ such that the Fourier series of $f$ is given by
 $\sum\limits^{\infty}_{\mu_1=0}\sum\limits^{\infty}_{\mu_2=0}\lambda_{\mu_1,\mu_2}\cos{2^{\mu_1}
x}\cos{2^{\mu_2}y}$, where $\lambda_{\mu_1,\mu_2}\in\mathbb{R}.$

\vskip 1.5cm
\sect{Auxiliary results}
\subsection{Jensen and Hardy inequalities}

\begin{lemma} \label{l2.9}
 \cite[Ch.~1]{7}
Let $ a_k \geq 0, 0 < \alpha \leq \beta <
\infty.$ Then

$$ \left( \sum\limits_{k=1}^{\infty} a_k^{\beta}
\right)^{{1}/{\beta}} \leq \left( \sum\limits_{k=1}^{\infty}
a_k^{\alpha} \right)^{{1}/{\alpha}}.  $$
\end{lemma}

\begin{lemma} \label{l2.8}
 \cite{8}  Let $ a_k \geq 0, b_k \geq 0.$

(A). Suppose
$\sum\limits_{k=1}^n a_k = a_n \gamma_n.$
If \, $ 1 \leq p < \infty, $ then

$$ \sum\limits_{k=1}^{\infty} a_k \Big( \sum\limits_{n=k}^{\infty}
b_n \Big)^p \ll \sum\limits_{k=1}^{\infty} a_k (b_k \gamma_k)^p.
$$

If \, $0 < p \leq 1$, then

$$ \sum\limits_{k=1}^{\infty} a_k \Big( \sum\limits_{n=k}^{\infty}
b_n \Big)^p \gg \sum\limits_{k=1}^{\infty} a_k (b_k \gamma_k)^p.
$$

(B). Suppose $ \ \sum\limits_{k=n}^{\infty}
a_k = a_n \beta_n.$
If \, $ 1 \leq p < \infty, $ then

$$ \sum\limits_{k=1}^{\infty} a_k \Big( \sum\limits_{n=1}^{k}
b_n \Big)^p \ll \sum\limits_{k=1}^{\infty} a_k (b_k \beta_k)^p.
$$

If \, $0 < p \leq 1$, then

$$ \sum\limits_{k=1}^{\infty} a_k \Big( \sum\limits_{n=1}^{k}
b_n \Big)^p \gg \sum\limits_{k=1}^{\infty} a_k (b_k \beta_k)^p.
$$
\end{lemma}

\vskip 0.5cm

\subsection{Results on angular approximation }

\begin{lemma} \label{l2.1}
\cite{6} Let $f\in L_p^0(\T^2), 1<p<\infty, n_i=0,1,2,\ldots, i=1,2. $ Then
$$
\|f-s_{n_1, \infty}(f)-s_{\infty, n_2}(f)+s_{n_1,n_2}(f)\|_{L_p(\T^2)}\ll
Y_{n_1,n_2}(f)_{L_p(\T^2)}.
$$
\end{lemma}

\begin{lemma} \label{l2.10}
 \cite{9}
Let $f\in L_p^0(\T^2), 1<p<q<\infty,
\theta=\frac1p-\frac1q, N_i=0,1,2, \ldots, i=1,2. $ Then
$$
Y_{2^{N_1}-1,2^{N_2}-1}(f)_{L_q(\T^2)}\ll  \left\{\sum\limits^{\infty}_{\nu_1=N_1}\sum\limits^{\infty}_{\nu_2=N_2}2^{(\nu_1+\nu_2)\theta q}Y^q_{2^{\nu_1}-1,2^{\nu_2}-1}(f)_{L_p(\T^2)}\right\}^{\frac1q}.
$$
\end{lemma}
Note that similar results for functions on $\mathbb{R}^d$ can be found in \cite{tom1}.

\vskip 0.5cm

\subsection{
Fourier coefficients of $L_p(\T^2)$-functions, Multipliers, and Littlewood-Paley theorem}

\begin{lemma}\label{l2.2} (The Marcinkiewicz multiplier theorem, \cite[Ch.~1]{7})
 Let the Fourier series of a function $f\in L_p^0(\T^2), 1<p<\infty, $ be
\begin{eqnarray}\label{e2.1}
\sum\limits^{\infty}_{n_1=1}\sum\limits^{\infty}_{n_2=1}\big(a_{n_1,n_2}\cos{n_1x}
\cos{n_2y}&+&b_{n_1,n_2}\sin{n_1x}\cos{n_2y}
\nonumber
\\ +\quad c_{n_1,n_2}\cos{n_1x}\sin{n_2y}&+& d_{n_1,n_2}\sin{n_1x}\sin{n_2y}\big)
=:
\sum\limits^{\infty}_{n_1=1}\sum\limits^{\infty}_{n_2=1}A_{n_1,n_2}(x,y).
\end{eqnarray}
%Let the number sequence  $\{\vartheta_{n_1,n_2}\}_{n_1,n_2=1}^\infty$ satisfy
Let the number sequence  $(\vartheta_{n_1,n_2})_{n_1,n_2=1}^\infty$ satisfy
$$
|\vartheta_{n_1,n_2}|\leq M, \;
\sum\limits^{2^{n_1}}_{m_1=2^{n_1-1}+1}|\vartheta_{m_1,n_2}-\vartheta_{m_1+1,n_2}|\leq
M, \;
\sum\limits^{2^{n_2}}_{m_2=2^{n_2-1}+1}|\vartheta_{n_1,m_2}-\vartheta_{n_1,m_2+1}|\leq
M
$$
and
$$
\sum\limits^{2^{n_1}}_{m_1=2^{n_1-1}+1}\sum\limits^{2^{n_2}}_{m_2=2^{n_2-1}+1}|
\vartheta_{m_1,m_2}-\vartheta_{m_1+1,m_2}-\vartheta_{m_1,m_2+1}+\vartheta_{m_1+1,m_2+1}|\leq
M
$$
for
some finite $M$ and  any  $n_i\in \N, i=1,2. $ Then the trigonometric series
$\sum\limits^{\infty}_{n_1=1}\sum\limits^{\infty}_{n_2=1}\vartheta_{n_1,n_2}A_{n_1,n_2}(x,y)$
is the Fourier series of a function $\phi\in L_p^0(\T^2)$ and
$$\|\phi\|_{L_p(\T^2)} \ll  \|f\|_{L_p(\T^2)}.$$
\end{lemma}

\begin{lemma} \label{l2.3} (The Littlewood-Paley theorem,
\cite[Ch.~1]{7})
Let the Fourier series of a function $f\in L_p^0(\T^2), 1<p<\infty, $ be given by
 (\ref{e2.1}).
Let
$\Delta_{0,0}:=A_{1,1}(x,y)$, \;
$$\Delta_{m_1,0}:=\sum\limits^{2^{m_1}}_{\nu_1=2^{m_1-1}+1}A_{\nu_1,1}(x,y)
\quad \textnormal{\it for} \;\;    m_1\in \N,
\quad
\Delta_{0,m_2}:=\sum\limits^{2^{m_2}}_{\nu_1=2^{m_2-1}+1}A_{1,\nu_2}(x,y)
\quad \textnormal{\it for}\;\;    m_2\in \N,$$
and
$$\Delta_{m_1,m_2}:=\sum\limits^{2^{m_1}}_{\nu_1=2^{m_1-1}+1}\sum\limits^{2^{m_2}}_{\nu_1=2^{m_2-1}+1}
A_{\nu_1,\nu_2}(x,y) \quad \textnormal{\it for} \;\;  m_1\in \N \quad \textnormal{\it and} \;\;   m_2\in \N.$$
Then
$$
\|f\|_{L_p(\T^2)}\asymp \left(\int\limits^{2\pi}_{0}
\int\limits^{2\pi}_{0}\left(\sum\limits^{\infty}_{\nu_1=0}
\sum\limits^{\infty}_{\nu_2=0}\Delta^2_{\nu_1,\nu_2}
\right)^{{p}/{2}} \dd x \dd y \right)^{1/p}.
$$
\end{lemma}

\begin{lemma} \label{l2.4}
(The Hardy-Littlewood-Paley theorem, \cite{4}) Let the Fourier series of a function $f\in L_1^0(\T^2)$ be given by
 (\ref{e2.1}).
\\
(A). Let $2\leq p<\infty$ and
$$
I:=\left(\sum\limits^{\infty}_{n_1=1}\sum\limits^{\infty}_{n_2=1}(|a_{n_1,n_2}|+|b_{n_1,n_2}|+|c_{n_1,n_2}|+|d_{n_1,n_2}|)^p(n_1n_2)^{p-2} \right)^{1/p}<\infty.
$$
Then  $f\in L_p^0(\T^2)$ and $ \|f\|_{L_p(\T^2)}\ll  I $.
\\
(B).  Let $f\in L_p^0(\T^2), 1<p\leq 2$. Then $I\ll  \|f\|_{L_p(\T^2)}.$
\end{lemma}

\begin{lemma} \label{l2.6}
 Let $f\in M_p, 1<p< \infty,$  $r_i \ge 0,
i=1,2$. Then
\begin{eqnarray}\label{hl1}
\|f\|_{L_p(\T^2)}\asymp \left(\sum\limits^{\infty}_{\nu_1=1}
\sum\limits^{\infty}_{\nu_2=1}a_{\nu_1,\nu_2}^{p}(\nu_1\nu_2)^{p-2} \right)^{{1}/{p}}
\end{eqnarray}
and
\begin{eqnarray}\label{hl2}
\|f^{(r_1, r_2)}\|_{L_p(\T^2)}\asymp
\left(\sum\limits^{\infty}_{\nu_1=1}\sum\limits^{\infty}_{\nu_2=1}a_{\nu_1,\nu_2}^{p}\nu_1^{r_1p+p-2}
\nu_2^{r_2p+p-2} \right)^{{1}/{p}}.
\end{eqnarray}
\end{lemma}
{\bf Proof.}
The proof of (\ref{hl1}) is given in \cite{15} (see also \cite{5}).

Let us verify (\ref{hl2}).
If  $ 2 \leq p < \infty,$ then the estimate from above in (\ref{hl2}) follows from Lemma \ref{l2.4} (A).
If  $1 < p <2,$ then Lemmas \ref{l2.2} and \ref{l2.3} imply
$$
I^p = \|f^{(r_1, r_2)}\|_{L_p(\T^2)}^p \asymp \int\limits_0^{2 \pi}
\int\limits_0^{2 \pi} \Big( \sum\limits_{\nu_1=0}^{\infty}
\sum\limits_{\nu_2=0}^{\infty} 2^{2 ( \nu_1 r_1 + \nu_2 r_2 )}
\Delta_{\nu_1, \nu_2}^2 \Big)^{{p}/{2}} \dd x  \dd y.
$$
Since $\frac{p}{2} < 1,$ using Lemma \ref{l2.9}, we get
$$
I^p \ll \sum\limits_{\nu_1=0}^{\infty}
\sum\limits_{\nu_2=0}^{\infty} 2^{p ( \nu_1 r_1 + \nu_2 r_2 )} \|
 \Delta_{\nu_1, \nu_2} \|_p^p.
 $$
In the paper  \cite{5} it was shown that  $ \|
 \Delta_{\nu_1, \nu_2} \|_p^p \ll 2^{(\nu_1 + \nu_2)p-1}
 a_{\floor{2^{\nu_1 -1}}+1,\floor{2^{\nu_2 -1}}+1 }^p.$
 Hence
$$ I^p \ll \sum\limits_{\nu_1=0}^{\infty}
\sum\limits_{\nu_2=0}^{\infty} 2^{p ( \nu_1 r_1 + \nu_2 r_2 ) +
(\nu_1 + \nu_2) p -1 }
 a_{\floor{2^{\nu_1 -1}}+1,\floor{2^{\nu_2 -1}}+1 }^p
$$
and by (\ref{mon})
$$ I^p \ll \sum\limits_{\nu_1=1}^{\infty}
\sum\limits_{\nu_2=1}^{\infty} a_{\nu_1, \nu_2}^p \nu_1^{(r_1+1)p-2} \nu_2^{(r_2+1)p-2}.
$$
Thus, we have proved the part ``$\ll$'' in  (\ref{hl2}).

To show the estimate from below, if $ 1 < p \leq 2$  then we simply use Lemma  \ref{l2.4} (B). If $ 2 < p < \infty,$
 we use the inequality
$$ \|f^{(r_1, r_2)} \|_p^p \gg \sum\limits_{\nu_1=1}^{\infty}
\sum\limits_{\nu_2=1}^{\infty} (\nu_1 \nu_2)^{-2} \Big(
\sum\limits_{\mu_1=\nu_1}^{\infty}
\sum\limits_{\mu_2=\nu_2}^{\infty} a_{\mu_1, \mu_2} \mu_1^{r_1}
\mu_2^{r_2} \Big)^p
$$
from the paper \cite{17}. Therefore, by (\ref{mon}),
$$ \|f^{(r_1, r_2)} \|_p^p \gg \sum\limits_{\nu_1=1}^{\infty}
\sum\limits_{\nu_2=1}^{\infty} a_{\nu_1, \nu_2}^p
\nu_1^{(r_1+1)p-2} \nu_2^{(r_2+1)p-2}.
$$
\hfill $\Box$

\begin{lemma} \label{l2.7}
 Let $f\in \Lambda_p, 1<p< \infty $. Then
\begin{eqnarray}\label{e2.4}
\|f\|_{L_p(\T^2)}\asymp
\left(\sum\limits^{\infty}_{\mu_1=0}\sum\limits^{\infty}_{\mu_2=0}\lambda_{\mu_1,\mu_2}^{2}\right)^{{1}/{2}}.
\end{eqnarray}
\end{lemma}
This lemma is well known in one dimension (\cite[Ch.~V, \S 8]{18}) but we failed to find its multivariate version. For the sake of completeness  we give a simple proof of this result.
\\[0.5cm]
{\bf Proof.}  Lemma  \ref{l2.3} yields
$$ I := \|f\|_{L_p(\T^2)}
\asymp \Bigg( \int\limits_0^{2 \pi} \int\limits_0^{2 \pi}
\Big(\sum\limits_{\nu_1=0}^{\infty} \sum\limits_{\nu_2=0}^{\infty}
 \Delta_{\nu_1, \nu_2}^2 \Big)^{{p}/{2}} \dd x  \dd y
\Bigg)^{{1}/{p}}.
$$
Since $\Delta_{\nu_1,\nu_2}=\lambda_{\nu_1,\nu_2} \cos{2^{\nu_1} x}
\cos{2^{\nu_2}y} $ for $ f \in \Lambda_p, $ we get
\begin{eqnarray}\label{e2.5}
I &\asymp& \Bigg(
\int\limits^{2\pi}_{0}\int\limits^{2\pi}_{0}\left(\sum\limits^{\infty}_{\nu_1=0}
\sum\limits^{\infty}_{\nu_2=0} \lambda_{\nu_1,\nu_2}^2 \big(
\cos{2^{\nu_1} x}\cos{2^{\nu_2}y}\big)^2
 \right)^{{p}/{2}} \dd x  \dd y
\Bigg)^{{1}/{p}}\\
\nonumber&\ll&
\left(\sum\limits^{\infty}_{\nu_1=0}\sum\limits^{\infty}_{\nu_2=0}
\lambda_{\nu_1,\nu_2}^2 \right)^{{1}/{2}}.
\end{eqnarray}
Let us now verify the estimate from below.
If $ 1 < p < 2,$ using Minkowski's inequality in (\ref{e2.5}), we have
$$ I \gg \Bigg( \sum\limits^{\infty}_{\nu_1=0}\sum\limits^{\infty}_{\nu_2=0}
\lambda_{\nu_1,\nu_2}^2 \Bigg( \int\limits_0^{2 \pi} \int\limits_0^{2
\pi} \big| \cos{2^{\nu_1} x}\cos{2^{\nu_2}y}\big|^p \dd x  \dd y
\Bigg)^{{2}/{p}}  \Bigg)^{{1}/{2}} \gg \Big(
\sum\limits^{\infty}_{\nu_1=0}\sum\limits^{\infty}_{\nu_2=0}
\lambda_{\nu_1,\nu_2}^2 \Big)^{{1}/{2}}.
$$
If  $2 \leq p < \infty, $ then $ I = \|f\|_{L_p(\T^2)} \gg
\|f\|_{L_2(\T^2)} \gg \Big(
\sum\limits^{\infty}_{\nu_1=0}\sum\limits^{\infty}_{\nu_2=0}
\lambda_{\nu_1,\nu_2}^2 \Big)^{{1}/{2}}.$
\hfill $\Box$

\vskip 0.5cm
\subsection{Auxiliary results for functions on $\T$}

Below we collect several useful results for functions of one variable.
As usual, $L_p (\T)$ is the collection of  $2\pi$-periodic measurable functions
 $f$ such that
$
\|f\|_{L_p (\T)} =\left(\int\limits_{0}^{2\pi}|f(x)|^p \, \dd x \,\right)^{1/p}<\infty
$
and
$L_p^0 (\T)$ is the collection of $f\in L_p(\T)$ such that $\int\limits_{0}^{2\pi}f(x)\dd x =0$.

Let $s_n(f)$ be the $n$-th partial sum of the Fourier series $f\in L_p(\T)$, i.e.,
$$
s_n(f)=s_n(f,x)=\frac1{\pi}\int\limits_{0}^{2\pi}f(x+t)\frac{\sin{(n+\frac12)t}}{2\sin{\frac{t}{2}}}\dd t.
$$
Let also $f^{(\rho)}$ be the Weyl derivative of
order $\rho>0$ of the function $f$.

For  $f\in L_p$ we define the difference of positive order $\alpha$ as follows
$$
\Delta_h^{\alpha}(f) = \sum\limits_{\nu=0}^{\infty} (-1)^{\nu}
\binom{\alpha}{\nu} f(x + (\alpha - \nu)h).
$$
We let  $\omega_{\alpha}(f,\delta)_{L_p(\T)}$ denote
the modulus of smoothness of $f$ of positive order $\alpha$ (\cite{butzer, 10, 11}), i.e.,
$$
\omega_{\alpha} (f, \delta)_{L_p(\T)} := \sup\limits_{|h|\leq
\delta} \| \Delta_h^{\alpha}(f) \|_{L_p(\T)}.
$$

\begin{lemma} \label{l3.3}
 \cite{butzer, 10}
 Let $f, g\in L_p^0(\T)$, $1<p<\infty$, and $\alpha>0, \beta>0$. Then

\begin{itemize}
\item[{\rm (a)}]
 $\triangle_h^\alpha(f+g) = \triangle_h^{\alpha} f + \triangle_h^{\alpha} g$;

\item[{\rm (b)}]
 $\triangle_h^\alpha(\triangle_h^\beta f) = \triangle_h^{\alpha+\beta} f$;

\item[{\rm (c)}]  $\|\triangle_h^{\alpha} f\|_{L_p(\T)} \ll \|f\|_{L_p(\T)}$.

\end{itemize}
\end{lemma}

\begin{lemma} \label{l3.4}
 \cite{butzer, 10}
Let $1<p<\infty$,  $\alpha>0,$ and $T_n$ be a trigonometric polynomial of degree at most  $n$,  $n\in \N $. Then

\begin{itemize}
\item[{\rm (a)}]
  we have for any  $ 0<|h|\leq \frac{\pi}{n}$
 $$
 \|\triangle_h^{\alpha} T_n\|_{L_p(\T)} \ll n^{-\alpha}\|T^{(\alpha)}_n\|_{L_p(\T)};
 $$

\item[{\rm (b)}]
we have
  $$
  \|T^{(\alpha)}_n\|_{L_p(\T)}\ll n^{\alpha}\|\triangle_{\frac{\pi}{n}}^{\alpha} T_n\|_{L_p(\T)}.
 $$
\end{itemize}
\end{lemma}

\begin{lemma} \label{l3.1}
 \cite{7}
Let $f \in L_p^0(\T)$, $1 < p < \infty$. Then
$$
  \|s_n(f)\|_{L_p(\T)}\ll \|f\|_{L_p(\T)}, \qquad n \in \mathbb{N}.
$$
\end{lemma}

\begin{lemma} \label{l3.5}
 \cite{14}
 Let $f\in L_p^0(\T)$, $1<p<q<\infty$,  $\theta:=\frac1p-\frac1q$, $n=0,1,2, \ldots$ Then

 $$
 \|f-s_{2^n}(f)\|_{L_q(\T)} \ll \left\{\sum\limits^{\infty}_{\nu=n}2^{\theta \nu q}\|f-s_{2^{\nu}}(f)\|_{L_p(\T)}^q\right\}^{1/q}.
 $$
\end{lemma}

\begin{lemma} \label{l3.6} (The Hardy-Littlewood inequality for fractional integrals,  \cite{18})

Let $f\in L_p^0(\T)$, $1<p<q<\infty$,  $\theta:=\frac1p-\frac1q$, $\alpha>0.$ Then
 $$
 \|s_{n}^{(\alpha)}(f)\|_{L_q(\T)} \ll \|s_{n}^{(\alpha+\theta)}(f)\|_{L_p(\T)}, \qquad n \in \mathbb{N}.
 $$
\end{lemma}

\vskip 1.5cm

\sect{Basic properties of the mixed moduli of smoothness}

We collect the main properties of the mixed moduli of smoothness of $L_p(\T^2)$-functions, $1<p<\infty$, in the following result.

\begin{theorem}\label{th6.1}
Let $f, g\in L_p(\T^2)$,  $1 < p < \infty,$  $\alpha_i> 0,$ $i=1,2.$ Then
\begin{itemize}
\item[{}] $$
\omega_{\alpha_1,\alpha_2}(f,\delta_1,0)_{L_p(\T^2)}=\omega_{\alpha_1,\alpha_2}(f,0,\delta_2)_{L_p(\T^2)}
=\omega_{\alpha_1,\alpha_2}(f,0,0)_{L_p(\T^2)}=0;
\leqno{\rm (1)}
$$

\item[{}] $$
\omega_{\alpha_1,\alpha_2}(f+g,\delta_1,\delta_2)_{L_p(\T^2)}\ll
     \omega_{\alpha_1,\alpha_2}(f,\delta_1,\delta_2)_{L_p(\T^2)}+
\omega_{\alpha_1,\alpha_2}(g,\delta_1,\delta_2)_{L_p(\T^2)}; \leqno{\rm (2)}
$$

\item[{}]
$$ \omega_{\alpha_1,\alpha_2}(f,\delta_1,\delta_2)_{L_p(\T^2)}\ll
\omega_{\alpha_1,\alpha_2}(f,t_1,t_2)_{L_p(\T^2)} \leqno{\rm (3)}
$$
for  $ 0<\delta_i\le t_i, i=1,2$;

\item[{}] $$
\frac{\omega_{\alpha_1,\alpha_2}(f,\delta_1,\delta_2)_{L_p(\T^2)}}{\delta_1^{\alpha_1}\delta_2^{\alpha_2}}\ll
\frac{\omega_{\alpha_1,\alpha_2}(f,t_1,t_2)_{L_p(\T^2)}}{t_1^{\alpha_1}t_2^{\alpha_2}}\leqno{\rm (4)}
$$ for  $0<t_i\le\delta_i\le 1, i=1,2$;

\item[{}]
$$\omega_{\alpha_1,\alpha_2}(f,\lambda_1\delta_1,\lambda_2\delta_2)_{L_p(\T^2)}\ll
\lambda_1^{\alpha_1}
\lambda_2^{\alpha_2}
\omega_{\alpha_1,\alpha_2}(f,\delta_1,\delta_2)_{L_p(\T^2)}
\leqno{\rm (5)}
$$
 for  $ \lambda_i>1, \,i=1,2;$

\item[{}]
$$
{\omega_{\beta_1,\beta_2}(f,\delta_1,\delta_2)_{L_p(\T^2)}}\ll
{\omega_{\alpha_1,\alpha_2}(f,\delta_1,\delta_2)_{L_p(\T^2)}}
\leqno{\rm (6)}
$$ for  $0<\alpha_i<\beta_i, \,i=1,2$;

\item[{}]
$$
{\omega_{\alpha_1,\alpha_2}(f,\delta_1,\delta_2)_{L_p(\T^2)}} \ll
\delta_1^{\alpha_1} \delta_2^{\alpha_2}\int_{\delta_1}^1\int_{\delta_2}^1
\frac{\omega_{\beta_1,\beta_2}(f,t_1,t_2)_{L_p(\T^2)}}{t_1^{\alpha_1}t_2^{\alpha_2}}
\frac{\dd t_1}{t_1}\frac{\dd t_2}{t_2}
\leqno{\rm (7)}
$$ for  $0<\alpha_i<\beta_i,$
$0<\delta_i\le\frac12$,
$\,i=1,2$ (Marchaud's inequality);

\item[{}]
$$
 \frac{\omega_{\alpha_1, \alpha_2}(f,\delta_1,\delta_2)_{L_p(\T^2)}}{ \delta_1^{\alpha_1}\delta_2^{\alpha_2}} \ll
\frac{\omega_{\beta_1,\beta_2}(f, \delta_1,\delta_2)_{L_p(\T^2)} }{ \delta_1^{\beta_1}\delta_2^{\beta_2}}\leqno{\rm (8)}
$$ for  $0<\alpha_i<\beta_i, \,i=1,2$;

\item[{}]
$$
\omega_{\beta_1+r_1, \beta_2+r_2}
(f, \delta_1,\delta_2)_{L_p(\T^2)}\ll
\delta_1^{r_1}\delta_2^{r_2}
\omega_{\beta_1, \beta_2}
(f^{(r_1, r_2)}, \delta_1, \delta_2)_{L_p(\T^2)}
\leqno{\rm (9)}
$$
for $\beta_i, r_i>0$,$\, i=1,2$;

\item[{}]
$$
\omega_{\beta_1, \beta_2} (f^{(r_1, r_2)}, \delta_1,\delta_2)_{L_p(\T^2)}
\ll
 \int\limits_0^{\delta_1} \int\limits_0^{\delta_2}
t_1^{-r_1}t_2^{-r_2} \omega_{\beta_1+r_1,\beta_2+r_2}(f, t_1,t_2)_{L_p(\T^2)}
\frac{\dd t_1}{t_1} \frac{\dd t_2}{t_2}
\leqno{\rm (10)}
$$
for $\beta_i, r_i>0$, $\, i=1,2$.

    \end{itemize}

\end{theorem}

\begin{remark}
Note that sharp versions of inequalities given in  (6)--(10) can be found in
Sections 8 and 10 below.
\end{remark}

{\bf Proof of Theorem \ref{th6.1}.}
Properties  (1), (2), and (3) follow from Lemma \ref{l3.3}  (a),
$$\sum\limits^{\infty}_{\nu=0}(-1)^{\nu}\left(^{\alpha}_{\nu}\right)=0,
$$
 and the definition of modulus of smoothness.

To prove  (4), we write

\begin{eqnarray*}
\frac{\omega_{\alpha_1, \alpha_2}(f, \delta_1,
 \delta_2)_{L_p(\T^2)}}{\delta_1^{\alpha_1} \delta_2^{\alpha_2}} &\asymp&
\frac{K(f, \delta_1, \delta_2, \alpha_1, \alpha_2,
p)}{\delta_1^{\alpha_1} \delta_2^{\alpha_2}} \\
 &\ll& \frac{K(f, t_1, t_2, \alpha_1, \alpha_2, p)}{t_1^{\alpha_1}
t_2^{\alpha_2}} \asymp \frac{\omega_{\alpha_1, \alpha_2}(f, t_1,
t_2)_{L_p(\T^2)}}{t_1^{\alpha_1} t_2^{\alpha_2}},
\end{eqnarray*}
where the equivalence between the modulus of smoothness and the $K$-functional is given by  Theorem \ref{th5.1} below.

Property (4) yields (5).
Note that Lemma  \ref{l3.3}  (b), (c) implies
 $\|\triangle_h^{\beta} f\|_{L_p(\T)} \ll \|\triangle_h^{\alpha} f\|_{L_p(\T)}$ for $0<\alpha<\beta$. Then
property (6) follows.

 The Marchaud inequality (7) can be easily shown using the direct and inverse approximation inequalities given below in this section by Theorem \ref{ll2.1}.
Property (8) is a simple consequence of (4) and (7).
Finally, properties (9) and (10) follow directly from Theorem \ref{T10} below.
\hfill $\Box$

Theorem \ref{th6.1} was partially proved in the papers  \cite{19, 23}.
For the one-dimensional case see \cite{butzer, sim-jat, 10}.

\vskip 0.5cm

\subsection{Jackson and Bernstein-Stechkin type inequalities}
The direct and inverse type results for periodic functions on $\T^2$ using the mixed modulus of smoothness are given by the following result.

\begin{theorem} \label{ll2.1}
Let $f\in L_p^0(\T^2)$, $1<p<\infty,$
$n_1, n_2=0,1,2,\ldots$, $\alpha_1, \alpha_2>0$. Then
\begin{eqnarray}\label{jackson-rr}
Y_{n_1,n_2}(f)_{L_p(\T^2)}&\ll& \omega_{\alpha_1,\alpha_2}(f,\frac{\pi}{n_1+1},\frac{\pi}{n_2+1})_{L_p(\T^2)}
\\
\nonumber
&\ll& \frac{1}{(n_1+1)^{\alpha_1}(n_2+1)^{\alpha_2}}
\sum\limits^{n_1+1}_{\nu_1=1}
\sum\limits^{n_2+1}_{\nu_2=1}\nu_1^{\alpha_1-1}\nu_2^{\alpha_2-1}Y_{\nu_1-1,\nu_2-1}(f)_{L_p(\T^2)}.
\end{eqnarray}
\end{theorem}
Theorem \ref{ll2.1} was proved for integers $\alpha_1, \alpha_2$ in the paper \cite{6}. In the general case,  Theorem \ref{ll2.1} follows from
Theorem \ref{th8.2} which is sharp versions of Jackson and Bernstein-Stechkin inequalities.
For non-mixed moduli of smoothness, see \cite[Ch.~5]{7} and \cite[Chs. V-VI]{timan-book}.

\begin{remark}
Note that the results of Theorem \ref{th6.1} and Theorem \ref{ll2.1} also hold in $L_p(\T^2)$, $p=1,\infty$; see \cite{isaac}
for  Theorem \ref{th6.1} (1)-(8) and Theorem \ref{ll2.1}.
Moreover, the Jackson inequality (\ref{jackson-rr}) is true in $L_p(\T^2)$, $0<p<1$; see \cite{runov}.
\end{remark}

\vskip 1.5cm
\sect{Constructive characteristic of the mixed moduli of smoothness}

\begin{theorem}\label{th4.1}
 Let $f \in L_p^0(\T^2)$, $1 < p < \infty,$ $\alpha_i > 0,$ $n_i \in \N,$ $i=1,2.$ Then
\begin{eqnarray*}
\omega_{\alpha_1,\alpha_2}\Big(f,\frac{\pi}{n_1},\frac{\pi}{n_2}\Big)_{L_p(\T^2)}  &\asymp&
n^{- \alpha_1}_1n^{- \alpha_2}_2 ||s_{n_1,n_2}^{(\alpha_1,\alpha_2)}(f)||_{L_p(\T^2)} + n^{- \alpha_1}_1 ||s_{n_1,\infty}^{(\alpha_1,0)}(f-s_{\infty,n_2}(f))||_{L_p(\T^2)}
\\
&+&n^{- \alpha_2}_2
||s_{\infty,n_2}^{(0,\alpha_2)}(f-s_{n_1,\infty}(f))||_{L_p(\T^2)}\\&+&||f-s_{n_1,\infty}(f)-s_{\infty,n_2}(f)+s_{n_1,n_2}(f)||_{L_p(\T^2)}.
\end{eqnarray*}
\end{theorem}
{\bf Proof.}
Using properties of the norm   we get, for any  $h_i$ and $n_i\in \N, i=1,2 $,
\begin{eqnarray*}
 \|\triangle_{h_1}^{\alpha_1}(\triangle_{h_2}^{\alpha_2}(f)) \|_{L_p(\T^2)}&\leq& \|\triangle_{h_1}^{\alpha_1}(\triangle_{h_2}^{\alpha_2}(f-s_{n_1,\infty}(f)-s_{\infty,n_2}(f)+s_{n_1,n_2}(f))) \|_{L_p(\T^2)}
\\&+&
\|\triangle_{h_1}^{\alpha_1}(\triangle_{h_2}^{\alpha_2}(s_{n_1,\infty}(f-s_{\infty,n_2}))) \|_{L_p(\T^2)}
\\&+&
\|\triangle_{h_1}^{\alpha_1}(\triangle_{h_2}^{\alpha_2}(s_{\infty,n_2}(f-s_{n_1,\infty}))) \|_{L_p(\T^2)}
\\
&+&
\|\triangle_{h_1}^{\alpha_1}(\triangle_{h_2}^{\alpha_2}(s_{n_1,n_2}(f))) \|_{L_p(\T^2)}
=:
I_1+I_2+I_3+I_4.
\end{eqnarray*}
First we estimate $I_1$ from above. Denote
$\varphi(x,y):=f-s_{n_1,\infty}(f)-s_{\infty,n_2}(f) +
s_{n_1,n_2}(f).$ By Lemma  \ref{l3.3}  (c), we have for a.e. $y$
$$
\left(\int\limits^{2\pi}_{0}|\triangle_{h_1}^{\alpha_1}(\triangle_{h_2}^{\alpha_2}(\varphi))|^p\dd x  \right)^{1/p}\ll  \left(\int\limits^{2\pi}_{0}|\triangle_{h_2}^{\alpha_2}(\varphi)|^p\dd x  \right)^{1/p}.
$$
Then
$$
\int\limits^{2\pi}_{0}\int\limits^{2\pi}_{0}|\triangle_{h_1}^{\alpha_1}(\triangle_{h_2}^{\alpha_2}(\varphi))|^p\dd
x \dd y \ll
\int\limits^{2\pi}_{0}\int\limits^{2\pi}_{0}|\triangle_{h_2}^{\alpha_2}(\varphi)|^p\dd
x \dd y .
$$
Therefore, $I_1\ll
\|\triangle_{h_2}^{\alpha_2}(\varphi)\|_{L_p(\T^2)}=: I_5$.
Using Lemma \ref{l3.3}  (c), we have for a.e.  $x$
$$
\left(\int\limits^{2\pi}_{0}|\triangle_{h_2}^{\alpha_2}(\varphi)|^p\dd y \right)^{1/p}\ll  \left(\int\limits^{2\pi}_{0}|\varphi|^p\dd y \right)^{1/p}.
$$
Then
$$
\int\limits^{2\pi}_{0}\int\limits^{2\pi}_{0}|\triangle_{h_2}^{\alpha_2}(\varphi)|^p\dd y
\dd x \ll \int\limits^{2\pi}_{0} \int\limits^{2\pi}_{0}|\varphi|^p\dd y
\dd x
$$
and  $I_5\ll \|\varphi\|_{L_p(\T^2)}$.
Thus,  $I_1\ll
\|f-s_{n_1,\infty}(f)-s_{\infty,n_2}(f)+s_{n_1,n_2}(f)\|_{L_p(\T^2)}$.

Similarly, to estimate  $I_2$ from above, we denote $\psi:=f-s_{\infty,n_2}(f)$.
By Lemma \ref{l3.3}  (c),  for a.e. $x$,
$$
\left(\int\limits^{2\pi}_{0}|\triangle_{h_1}^{\alpha_1}(\triangle_{h_2}^{\alpha_2}(s_{n_1,\infty}(\psi)))|^p\dd y
\right)^{1/p}\ll
\left(\int\limits^{2\pi}_{0}|\triangle_{h_1}^{\alpha_1}(s_{n_1,\infty}(\psi))|^p\dd y
\right)^{1/p}.
$$
Hence,
$$
\int\limits^{2\pi}_{0}\int\limits^{2\pi}_{0}|\triangle_{h_1}^{\alpha_1}(\triangle_{h_2}^{\alpha_2}(s_{n_1,\infty}(\psi)))|^p\dd y
\dd x \ll
\int\limits^{2\pi}_{0}\int\limits^{2\pi}_{0}|\triangle_{h_1}^{\alpha_1}(s_{n_1,\infty}(\psi))|^p\dd y
\dd x
$$
and $I_2\ll
\|\triangle_{h_1}^{\alpha_1}(s_{n_1,\infty}(\psi))\|_{L_p(\T^2)}=: I_6$.
Now by using  Lemma \ref{l3.4}  (a), we get for a.e. $y$  and any
$h_1\in(0,\frac{\pi}{n_1})$
$$
\left(\int\limits^{2\pi}_{0}|\triangle_{h_1}^{\alpha_1}(s_{n_1,\infty}(\psi))|^p\dd x
\right)^{1/p}\ll  n_1^{-\alpha_1}
\left(\int\limits^{2\pi}_{0}|s^{(\alpha_1,0)}_{n_1,\infty}(\psi)|^p\dd x
\right)^{1/p}.
$$
Then we get
$$
\int\limits^{2\pi}_{0}\int\limits^{2\pi}_{0}|\triangle_{h_1}^{\alpha_1}(s_{n_1,\infty}(\psi))|^p\dd x
\dd y\ll  n_1^{-\alpha_1}
\int\limits^{2\pi}_{0}\int\limits^{2\pi}_{0}|s^{(\alpha_1,0)}_{n_1,\infty}(\psi)|^p\dd x
\dd y.
$$
Hence,  $I_6\ll  n_1^{-\alpha_1}\|s_{n_1,\infty}^{(\alpha_1,0)}(\psi)\|_{L_p(\T^2)}$.
We have shown that, for
$0<|h_1|<\frac{\pi}{n_1}$,
$$
I_2\ll  n_1^{-\alpha_1}\|s_{n_1,\infty}^{(\alpha_1,0)}(f-s_{\infty,n_2}(f))\|_{L_p(\T^2)}.
$$
Similarly, we obtain
\begin{eqnarray*}
I_3&\ll&  n_2^{-\alpha_2}\|s_{\infty,n_2}^{(0,\alpha_2)}(f-s_{n_1,\infty}(f))\|_{L_p(\T^2)},
\\
I_4&\ll&  n_1^{-\alpha_1}n_2^{-\alpha_2}\|s_{n_1,n_2}^{(\alpha_1,\alpha_2)}(f)\|_{L_p(\T^2)},
\end{eqnarray*}
for  $0<|h_1|<\frac{\pi}{n_1},$ $0<|h_2|<\frac{\pi}{n_2}$.

Finally,
\begin{eqnarray*}
\omega_{\alpha_1,\alpha_2}(f,\frac{\pi}{n_1},\frac{\pi}{n_2})_{L_p(\T^2)}
&\ll&
||f-s_{n_1,\infty}(f)-s_{\infty,n_2}(f)+s_{n_1,n_2}(f)||_{L_p(\T^2)}
\\
&+&
n^{- \alpha_1}_1 ||s_{n_1,\infty}^{(\alpha_1,0)}(f-s_{\infty,n_2}(f))||_{L_p(\T^2)}+n^{- \alpha_2}_2 ||s_{\infty,n_2}^{(0,\alpha_2)}(f-s_{n_1,\infty}(f))||_{L_p(\T^2)}
\\
&+&n^{- \alpha_1}_1n^{- \alpha_2}_2 ||s_{n_1,n_2}^{(\alpha_1,\alpha_2)}(f)||_{L_p(\T^2)},
\end{eqnarray*}
and the upper estimate follows.

To prove the estimate from below, we use Lemma \ref{l2.1}, Theorem \ref{ll2.1}, and properties of the mixed modulus of integer order
\begin{eqnarray*}
A_1&:=&
||f-s_{n_1,\infty}(f)-s_{\infty,n_2}(f)+s_{n_1,n_2}(f)||_{L_p(\T^2)}\ll
Y_{n_1,n_2}(f)_{L_p(\T^2)}
\\
&\ll&
\omega_{\floor{\alpha_1}+1,\floor{\alpha_2}+1}(f,\frac{\pi}{n_1+1},\frac{\pi}{n_2+1})_{L_p(\T^2)}\ll \omega_{\floor{\alpha_1}+1,\floor{\alpha_2}+1}(f,\frac{\pi}{n_1},\frac{\pi}{n_2})_{L_p(\T^2)}.
\end{eqnarray*}
By Lemma \ref{l3.3}  (b), we get
$$
A_1\leq \sup\limits_{|h_i|\leq
\frac{\pi}{n_i},i=1,2}\|\triangle_{h_1}^{\floor{\alpha_1}+1-\alpha_1}(\triangle_{h_2}^{\floor{\alpha_2}+1-\alpha_2}
(\triangle_{h_1}^{\alpha_1}(\triangle_{h_2}^{\alpha_2}
(f))))\|_{L_p(\T^2)}.
$$
Using Lemma  \ref{l3.3}  (c),

$$
A_1\leq \sup\limits_{|h_i|\leq
\frac{\pi}{n_i},i=1,2}\|\triangle_{h_1}^{\alpha_1}(\triangle_{h_2}^{\alpha_2}
(f))\|_{L_p(\T^2)}=\omega_{\alpha_1,\alpha_2}(f,\frac{\pi}{n_1},\frac{\pi}{n_2})_{L_p(\T^2)}.
$$
Now let us estimate
$$
A_2
:=\|s_{n_1,\infty}^{(\alpha_1,0)}(f-s_{\infty,n_2}(f))\|_{L_p(\T^2)}.
$$
Defining  $\gamma(x,y):= f(x,y)-s_{\infty,n_2}(f)$ and using Lemma
 \ref{l3.4}  (b), we have for a.e. $y$
$$
\left(\int\limits^{2\pi}_{0}|s_{n_1,\infty}(\gamma)|^p\dd x \right)^{1/p}\ll
n_1^{\alpha_1}
\left(\int\limits^{2\pi}_{0}|\triangle^{\alpha_1}_{\frac{\pi}{n_1}}(s_{n_1,\infty}(\gamma))|^p\dd x \right)^{1/p}.
$$
Hence,
$$
\int\limits^{2\pi}_{0}\int\limits^{2\pi}_{0}|s_{n_1,\infty}(\gamma)|^p\dd x \dd y\ll
n_1^{\alpha_1 p}
\int\limits^{2\pi}_{0}\int\limits^{2\pi}_{0}|\triangle^{\alpha_1}_{\frac{\pi}{n_1}}(s_{n_1,\infty}(\gamma))|^p\dd x \dd y,
$$
and
$$
A_2\ll n_1^{\alpha_1
}\|s_{n_1,\infty}(\triangle^{\alpha_1}_{\frac{\pi}{n_1}}(\gamma))\|_{L_p(\T^2)}.
$$
Lemma \ref{l3.1} implies
$$
A_2\ll n_1^{\alpha_1
}\|\triangle^{\alpha_1}_{\frac{\pi}{n_1}}(f-s_{\infty,n_2}(f))\|_{L_p(\T^2)}.
$$
Defining  $\triangle^{\alpha_1}_{\frac{\pi}{n_1}}(f) =: F$, we get
 $A_2\ll n_1^{\alpha_1 }\|F-s_{\infty,n_2}(F)\|_{L_p(\T^2)}$.
Since $s_{0,\infty}(F)=s_{0,n_2}(F)=0$, then
$$
A_2\ll n_1^{\alpha_1
}\|F-s_{0,\infty}(F)-s_{\infty,n_2}(F)+s_{0,n_2}(F)\|_{L_p(\T^2)}.
$$
Therefore by Lemma \ref{l2.1}, Theorem \ref{ll2.1}, and the properties of the mixed moduli of smoothness
$$
A_2\ll n_1^{\alpha_1
}\omega_{
\floor{\alpha_1}+1,\floor{\alpha_2}+1
}(F,\pi,\frac{\pi}{n_2+1})_{L_p(\T^2)}\ll
n_1^{\alpha_1
}\omega_{
\floor{\alpha_1}+1,\floor{\alpha_2}+1
}(F,\pi,\frac{\pi}{n_2})_{L_p(\T^2)}.
$$
Using Lemma \ref{l3.3} (b), we have
$$
A_2\ll  n_1^{\alpha_1 } \sup\limits_{|h_1|\leq\pi, |h_2|\leq
\frac{\pi}{n_2}}\|\triangle_{h_1}^{
\floor{\alpha_1}+1
}(\triangle_{h_2}^{
\floor{\alpha_2}+1-\alpha_2}
(\triangle_{h_2}^{\alpha_2}(F)))\|_{L_p(\T^2)}.
$$
Now Lemma \ref{l3.3} (c) yields
$$
A_2\ll  n_1^{\alpha_1 } \sup\limits_{ |h_2|\leq \frac{\pi}{n_2}}\|
\triangle_{h_2}^{\alpha_2}(F)\|_{L_p(\T^2)}=n_1^{\alpha_1 }
\sup\limits_{ |h_2|\leq \frac{\pi}{n_2}}\|
\triangle_{h_2}^{\alpha_2}(\triangle_{\frac{\pi}{n_1}}^{\alpha_1}(f))\|_{L_p(\T^2)}
\ll n_1^{\alpha_1 }
\omega_{\alpha_1,\alpha_2}(f,\frac{\pi}{n_1},\frac{\pi}{n_2})_{L_p(\T^2)}.
$$
Similarly, one can show that
$$
A_3:=\|s_{\infty,n_2}^{(0,\alpha_2)}(f-s_{n_1,\infty}(f))\|_{L_p(\T^2)}\ll
n_2^{\alpha_2}
\omega_{\alpha_1,\alpha_2}(f,\frac{\pi}{n_1},\frac{\pi}{n_2})_{L_p(\T^2)}
$$
and
$$
A_4:=\|s_{n_1,n_2}^{(\alpha_1,\alpha_2)}(f)\|_{L_p(\T^2)}\ll
n_1^{\alpha_1 }n_2^{\alpha_2}
\omega_{\alpha_1,\alpha_2}(f,\frac{\pi}{n_1},\frac{\pi}{n_2})_{L_p(\T^2)}.
$$
Finally,
$$
||f-s_{n_1,\infty}(f)-s_{\infty,n_2}(f)+s_{n_1,n_2}(f)||_{L_p(\T^2)}+
n^{- \alpha_1}_1
||s_{n_1,\infty}^{(\alpha_1,0)}(f-s_{\infty,n_2}(f))||_{L_p(\T^2)}+
$$
$$
 n^{- \alpha_2}_2||s_{\infty,n_2}^{(0,\alpha_2)}(f-s_{n_1,\infty}(f))||_{L_p(\T^2)}+
n^{- \alpha_1}_1n^{- \alpha_2}_2 ||s_{n_1,n_2}^{(\alpha_1,\alpha_2)}(f)||_{L_p(\T^2)}
\ll  \omega_{\alpha_1,\alpha_2}(f,\frac{\pi}{n_1},\frac{\pi}{n_2})_{L_p(\T^2)},
$$
i.e., the required  estimate from below.
\hfill $\Box$

%Theorem \ref{th4.1} was stated in the papers  \cite{21, 19} without the proof.
Theorem \ref{th4.1} was stated in the papers  \cite{21, 19} without proof.
This statement  is called the realization result;
in  dimension one, see \cite{ditzian} for the moduli of smoothness of integer order and \cite{sbornik} for the fractional case.
For  the non-mixed moduli of smoothness of functions on $\R^d$ see, e.g., \cite[(5.3)]{dit-ul}.

\vskip 1.5cm

\sect{
The mixed moduli of smoothness and the $K$-functionals
}

\begin{theorem}\label{th5.1}
Let $ f \in L_p^0(\T^2)$, $1 < p <\infty,$ $\alpha_i
> 0, 0 < \delta_i \leq \pi, i =1,2. $ Then
\begin{eqnarray}\label{e5.1}
\omega_{\alpha_1, \alpha_2}(f, \delta_1, \delta_2 )_{L_p(\T^2)}
\asymp K(f, \delta_1, \delta_2, \alpha_1, \alpha_2, p).
\end{eqnarray}
\end{theorem}
{\bf Proof.}
For any  $ \delta_i \in (0, \pi ]$ take integers
$n_i$  such that $ \frac{\pi}{n_i+1} < \delta_i \leq
\frac{\pi}{n_i}, i = 1,2. $  If $ f \in L_p^0(\T^2),$ then
$$
(s_{n_1+1,\infty}(f)-s_{n_1+1,n_2+1}(f))\in
W_p^{(\alpha_1,0)};\quad(s_{\infty,n_2+1}(f)-s_{n_1+1,n_2+1}(f))\in
W_p^{(0,\alpha_2)};
$$
$$
s_{n_1+1,n_2+1}(f)\in W_p^{(\alpha_1,\alpha_2)}.
$$
Then it is clear that
\begin{align*}
 &K(f, \delta_1, \delta_2, \alpha_1, \alpha_2, p)
 \\&\le\|f-(s_{n_1+1,\infty}(f)-s_{n_1+1,n_2+1}(f))-(s_{\infty,n_2+1}(f)-s_{n_1+1,n_2+1}(f))-
s_{n_1+1,n_2+1}(f)\|_{L_p(\mathbb{T}^2)}
\\
&+\delta_1^{\alpha_1}\|s_{n_1+1,\infty}^{(\alpha_1,0)}(f)-s_{n_1+1,n_2+1}^{(\alpha_1,0)}(f)\|_{L_p(\mathbb{T}^2)}
+\,
\delta_2^{\alpha_2}\|s_{\infty,n_2+1}^{(0,\alpha_2)}(f)-s_{n_1+1,n_2+1}^{(0,\alpha_2)}(f)\|_{L_p(\mathbb{T}^2)}
\\
&\qquad\qquad\qquad\qquad\qquad\qquad\qquad\qquad\qquad\qquad\qquad\qquad+\delta_1^{\alpha_1}\delta_2^{\alpha_2}\|s_{n_1+1,n_2+1}^{(\alpha_1,\alpha_2)}(f)\|_{L_p(\mathbb{T}^2)}
\\&\ll
\|f-s_{n_1+1,\infty}(f)-s_{\infty,n_2+1}(f)+s_{n_1,n_2}(f)\|_{L_p(\mathbb{T}^2)}
+\,
n_1^{-\alpha_1}\|s_{n_1+1,\infty}^{(\alpha_1,0)}(f-s_{\infty,n_2+1}(f))\|_{L_p(\mathbb{T}^2)}
\\
&\qquad\qquad\qquad\qquad+\,
n_2^{-\alpha_2}\|s_{\infty,n_2+1}^{(0,\alpha_2)}(f-s_{n_1+1,\infty}(f))\|_{L_p(\mathbb{T}^2)}+
n_1^{-\alpha_1}n_2^{-\alpha_2}\|s_{n_1+1,n_2+1}^{(\alpha_1,\alpha_2)}(f)\|_{L_p(\mathbb{T}^2)}.
\end{align*}
By Theorem  \ref{th4.1}, the last expression is bounded  by $\omega_{\alpha_1, \alpha_2} \big(f,
\frac{\pi}{n_1+1}, \frac{\pi}{n_2+1} \big)_{L_p(\mathbb{T}^2)}$ and therefore

\begin{eqnarray}\label{e5.2}
 K(f, \delta_1, \delta_2, \alpha_1, \alpha_2, p) \ll \ll
\omega_{\alpha_1, \alpha_2} \big(f, \delta_1, \delta_2
\big)_{L_p(\mathbb{T}^2)}, \end{eqnarray}
and  the part
``$\gg$''
 in  estimate  (\ref{e5.1})  follows.

Let us prove
the part ``$\ll$''
 in  estimate  (\ref{e5.1}).
Take any functions $g_1 \in
W_p^{(\alpha_1,0)}$, $g_2 \in W_p^{(0, \alpha_2)}$, and $g \in
W_p^{(\alpha_1, \alpha_2) }$. Then using Theorem \ref{th6.1} (2), we get

\begin{eqnarray*} \omega_{\alpha_1, \alpha_2}(f, \delta_1, \delta_2 )_{L_p(\mathbb{T}^2)}
 &\ll&
 \omega_{\alpha_1,\alpha_2}(f -  g_1 -g_2 -g , \delta_1, \delta_2 )_{L_p(\mathbb{T}^2)}
   +  \omega_{\alpha_1,\alpha_2}( g_1, \delta_1, \delta_2 )_{L_p(\mathbb{T}^2)}
   \\ &+&
  \omega_{\alpha_1,\alpha_2}( g_2, \delta_1, \delta_2 )_{L_p(\mathbb{T}^2)} +
  \omega_{\alpha_1,\alpha_2}(  g , \delta_1, \delta_2 )_{L_p(\mathbb{T}^2)} =: J_1 + J_2 +J_3 +J_4 .
\end{eqnarray*}
By Lemma \ref{l3.3}, we have  $J_1 \ll \|f - g_1 - g_2 -
g\|_{L_p(\mathbb{T}^2)}. $

To estimate $J_2$, for any  $ \delta_i \in (0, \pi] $ we take
integers $n_i$ such that $
\frac{\pi}{2^{n_{i+1}}} < \delta_i \leq \frac{\pi}{2^{n_i}}, \ i
=1,2. $
We consider $B_2 := \omega_{\alpha_1, \alpha_2} \big(g_1,
\frac{\pi}{2^{n_1}}, \frac{\pi}{2^{n_2}}
\big)_{L_p(\mathbb{T}^2)}.$
Lemma  \ref{l3.3} yields
\begin{eqnarray*}
B_2 &\ll& \omega_{\alpha_1, \alpha_2} \big(g_1 - s_{2^{n_1},
\infty}(g_1) , \frac{\pi}{2^{n_1}}, \frac{\pi}{2^{n_2}}
\big)_{L_p(\mathbb{T}^2)} + \omega_{\alpha_1, \alpha_2}
\big(s_{2^{n_1},\infty}(g_1), \frac{\pi}{2^{n_1}},
\frac{\pi}{2^{n_2}} \big)_{L_p(\mathbb{T}^2)}
\\ &\ll&
\|g_1 - s_{2^{n_1}, \infty}(g_1)\|_{L_p(\mathbb{T}^2)} +
\sup\limits_{|h_1|\leq \frac{\pi}{2^{n_1}}} \|
\Delta_{h_1}^{\alpha_1}(s_{2^{n_1},\infty}(g_1))
\|_{L_p(\mathbb{T}^2)} =: J_{21} + J_{22}.
\end{eqnarray*}
Using Lemma  \ref{l3.4}  (a) and Lemma  \ref{l3.1}, we get for a.e.  $y$ and
 $0<h_1\le\frac{\pi}{2^{n_1}}$
   $$
   \Bigg( \int\limits_0^{2 \pi} \big| \Delta_{h_1}^{\alpha_1}
   s_{2^{n_1},
    \infty}(g_1) \big|^p \dd x  \Bigg)^{1/p} \ll 2^{-n_1 \alpha_1}
    \Bigg( \int\limits_0^{2 \pi} \big|  s_{2^{n_1},
    \infty}^{(\alpha_1,0)}(g_1) \big|^p \dd x  \Bigg)^{1/p} \ll 2^{-n_1 \alpha_1}
    \Bigg( \int\limits_0^{2 \pi} \big| g_1^{(\alpha_1,0)} \big|^p \dd x
    \Bigg)^{1/p}.
    $$
Then the inequality
$$
\int\limits_0^{2\pi}\int\limits_0^{2\pi}|\Delta_{h_1}^{\alpha_1}s_{2^{n_1},\infty}(g_1)|^p\dd x \dd y\ll2^{-n_1\alpha_1p}\int\limits_0^{2\pi}\int\limits_0^{2\pi}|g_1^{(\alpha_1,0)}|^p\dd x \dd y
$$
implies
$J_{22}\ll2^{-n_1\alpha_1}\|g_1^{(\alpha_1,0)}\|_{L_p(\mathbb{T}^2)}$.
Since
 $g_1\in L_p^0(\mathbb{T}^2)$, then Lemmas
\ref{l2.2}  and \ref{l2.3} give
$$
J_{21}\ll\Bigg(\int\limits_0^{2\pi}\int\limits_0^{2\pi}\Big(\sum\limits_{\nu_1=n_1}^\infty
\sum\limits_{\nu_2=0}^\infty\Delta_{\nu_1,\nu_2}^2\Big)^{{p}/{2}}\dd x \dd y\Bigg)^{{1}/{p}}
\ll2^{-n_1\alpha_1}\Bigg(\int\limits_0^{2\pi}\int\limits_0^{2\pi}\Big(\sum\limits_{\nu_1=n_1}^\infty
\sum\limits_{\nu_2=0}^\infty2^{2\nu_1\alpha_1}\Delta_{\nu_1,\nu_2}^2\Big)^{{p}/{2}}\dd x \dd y\Bigg)^{{1}/{p}}.
$$
Using again Lemmas \ref{l2.2} and \ref{l2.3}, and further Lemma \ref{l3.1}, we get
\begin{eqnarray*}
J_{21} &\ll&
2^{-n_1\alpha_1}
\|g_1^{(\alpha_1,0)}-s_{2^{n_1}-1,\infty}(g_1^{(\alpha_1,0)})\|_{L_p(\mathbb{T}^2)}
\\
&\ll&
2^{-n_1\alpha_1}
\Big(\|g_1^{(\alpha_1,0)}\|_{L_p(\mathbb{T}^2)}+\|s_{2^{n_1}-1,\infty}(g_1^{(\alpha_1,0)})\|_{L_p(\mathbb{T}^2)}
\Big)
\ll
2^{-n_1\alpha_1}
\|g_1^{(\alpha_1,0)}\|_{L_p(\mathbb{T}^2)}.
\end{eqnarray*}
The estimates for $J_{21}$ and $J_{22}$ imply
$$ B_2 \ll 2^{- n_1 \alpha_1}  \|g_1^{(\alpha_1,0)}\|_{L_p(\mathbb{T}^2)}.$$
Using properties of moduli of smoothness, we get $
\omega_{\alpha_1, \alpha_2} \big(g_1, \delta_1, \delta_2
\big)_{L_p(\mathbb{T}^2)} \ll \omega_{\alpha_1, \alpha_2}
\big(g_1, \frac{\pi}{2^{n_1}}, \frac{\pi}{2^{n_2}}
\big)_{L_p(\mathbb{T}^2)} $ and
$$ J_2 \ll 2^{- n_1 \alpha_1}
\|g_1^{(\alpha_1,0)}\|_{L_p(\mathbb{T}^2)}. $$
Similarly,
$$ J_3 \ll 2^{- n_2 \alpha_2}
\|g_2^{(0, \alpha_2)}\|_{L_p(\mathbb{T}^2)}\qquad\mbox{and} \qquad
J_4 \ll 2^{- n_1 \alpha_1 - n_2 \alpha_2} \|g^{(\alpha_1,
\alpha_2)}\|_{L_p(\mathbb{T}^2)}. $$
Finally, combining
estimates for $J_1, J_2, J_3$, and $J_4,$ we get
 \begin{eqnarray*}
 \omega_{\alpha_1, \alpha_2} \big(g_1, \delta_1, \delta_2
\big)_{L_p(\mathbb{T}^2)} &\ll& \|f - g_1 - g_2 -
g\|_{L_p(\mathbb{T}^2)} + \delta_1^{\alpha_1}
\|g_1^{(\alpha_1,0)}\|_{L_p(\mathbb{T}^2)}  \\ &+&
\delta_2^{\alpha_2} \|g_2^{(0, \alpha_2)}\|_{L_p(\mathbb{T}^2)} +
\delta_1^{\alpha_1} \delta_2^{\alpha_2}
\|g^{(\alpha_1,\alpha_2)}\|_{L_p(\mathbb{T}^2)}.
 \end{eqnarray*}
Since the last inequality holds for any  $g_1\in W_p^{(\alpha_1,0)}$, $g_2\in W_p^{(0, \alpha_2)},$ and $g \in W_p^{(\alpha_1,\alpha_2)},$ we get
\begin{eqnarray}\label{e5.3}
 \omega_{\alpha_1, \alpha_2} \big(g_1, \delta_1, \delta_2
\big)_{L_p(\mathbb{T}^2)} \ll
K(f,\delta_1,\delta_2,\alpha_1,\alpha_2,p) \end{eqnarray}
and therefore the proof of  the part ``$\ll$'' in estimate  (\ref{e5.1})  follows.
\hfill $\Box$

In the case of integers $\alpha_1$ and $\alpha_2$ Theorem \ref{th5.1} was proved in  the paper \cite{12} for
$ 1 \leq p \leq \infty$, and in the paper  \cite{13} for  $p =\infty$ using different methods.
In the one-dimensional case and in the multivariate case for non-mixed moduli of smoothness,
the equivalence between the moduli of smoothness and the corresponding $K$-functionals was proved in \cite{johnen} (see also \cite[p. 339]{ben-sha}).

\vskip 1.5cm

\sect{The mixed moduli of smoothness of $L_p$-functions and their Fourier coefficients}

The classical Riemann-Lebesgue lemma states that  the Fourier coefficients  of an $L_p(\T^2)$ function tend to $0$ as $|n|\to \infty$. Its quantitative version is written as follows:
$$
\rho_{n_1,n_2}
\ll
\omega_{\alpha_1,\alpha_2}
\Big(f, \frac{1}{n_1},\frac1{n_2} \Big)_{L_p(\T^2)},
$$
where
 $$
 \rho_{n_1,n_2}=|a_{n_1,n_2}|+|b_{n_1,n_2}|+|c_{n_1,n_2}|+|d_{n_1,n_2}|$$
and the Fourier series of $f\in L_p^0(\T^2), 1<p<\infty,$ is given by (\ref{e2.1}).
We extend this estimate by writing the following two-sided inequalities using
weighted tail-type sums of the Fourier series.

\begin{theorem}\label{th7.1}
Let $f \in L_p^0$, $ 1 < p < \infty,$ and the Fourier series of $f$ be given by  (\ref{e2.1}).
  Let  $\tau:=\max(2,p)$,
$\theta:=\min(2,p),$  $\alpha_1 > 0,\alpha_2 > 0$, and $ n_1 \in \N, n_2 \in \N$. Then
\begin{eqnarray}\label{e56.3}
 \mathbb{I}(\theta)\ll
\omega_{\alpha_1,\alpha_2}
\Big(f, \frac{1}{n_1},\frac1{n_2} \Big)_{L_p(\T^2)} \ll
\mathbb{I}(\tau),
\end{eqnarray}
where
$$
\mathbb{I}(s) := \frac{1}{n^{\alpha_1}_1}\frac{1}{n_2^{\alpha_2}}
\Bigg\{
\sum\limits_{\nu_1=1}^{n_1}\sum\limits_{\nu_2=1}^{n_2} \rho_{\nu_1,\nu_2}^s
\nu_1^{(\alpha_1+1)s-2}\nu_2^{(\alpha_2+1)s-2}
 \Bigg\}^{{1}/{s}} +
 \frac{1}{n^{\alpha_1}_1}
\Bigg\{ \sum\limits_{\nu_1=1}^{n_1}\sum\limits_{\nu_2=n_2+1}^{\infty} \rho_{\nu_1,\nu_2}^s
\nu_1^{(\alpha_1+1)s-2}\nu_2^{s-2}
 \Bigg\}^{{1}/{s}}
$$
$$
+\frac{1}{n_2^{\alpha_2}}
\Bigg\{ \sum\limits_{\nu_1=n_1+1}^{\infty}\sum\limits_{\nu_2=1}^{n_2} \rho_{\nu_1,\nu_2}^s
\nu_1^{s-2}\nu_2^{(\alpha_2+1)s-2}
 \Bigg\}^{{1}/{s}} +
 \Bigg\{ \sum\limits_{\nu_1=n_1+1}^{\infty}\sum\limits_{\nu_2=n_2+1}^{\infty} \rho_{\nu_1,\nu_2}^s
(\nu_1\nu_2)^{s-2}
 \Bigg\}^{{1}/{s}}
$$ for $1<s<\infty.$
\end{theorem}

The next two theorems provide sharper results than (\ref{e56.3}) for special classes of functions defined in Section 2.4. In particular, these results show that the parameters
 $\tau=\max(2,p)$ and $\theta=\min(2,p)$  in (\ref{e56.3}) cannot be extended.

\begin{theorem}\label{th7.2}
 Let $f \in M_p, 1 < p < \infty, \alpha_1 > 0, \alpha_2 > 0,  n_1 \in \N,  n_2 \in \N$.
 Then
$$
\omega_{\alpha_1,\alpha_2} \Big(f, \frac{1}{n_1}, \frac{1}{n_2} \Big)_{L_p(\T^2)} \asymp
\frac{1}{n_1^{\alpha_1}n_2^{\alpha_2}} \Bigg\{ \sum\limits_{\nu_1=1}^{n_1} \sum\limits_{\nu_2=1}^{n_2} a_{\nu_1,\nu_2}^p
\nu_1^{(\alpha_1+1)p -2}\nu_2^{(\alpha_2+1)p -2} \Bigg\}^{{1}/{p}}
$$
$$
+
\frac{1}{n_1^{\alpha_1}} \Bigg\{ \sum\limits_{\nu_1=1}^{n_1} \sum\limits_{\nu_2=n_2+1}^{\infty} a_{\nu_1,\nu_2}^p \nu_1^{(\alpha_1
+1)p -2}\nu_2^{p -2} \Bigg\}^{{1}/{p}}
+
\frac{1}{n_2^{\alpha_2}} \Bigg\{ \sum\limits_{\nu_1=n_1+1}^{\infty} \sum\limits_{\nu_2=1}^{n_2} a_{\nu_1,\nu_2}^p \nu_1^{p -2}\nu_2^{(\alpha_2+1)p -2} \Bigg\}^{{1}/{p}}
 $$
 $$
 +\Bigg\{ \sum\limits_{\nu_1=n_1+1}^{\infty} \sum\limits_{\nu_2=n_2+1}^{\infty} a_{\nu_1,\nu_2}^p (\nu_1\nu_2)^{p -2} \Bigg\}^{{1}/{p}}.
 $$
\end{theorem}

\begin{theorem}\label{th7.3}
 Let $f \in \Lambda_p, 1 < p < \infty, \alpha_1 > 0, \alpha_2>0, n_i=0,1,2,\ldots,   i=1,2$.
 Then
$$ \omega_{\alpha_1,\alpha_2} \Big(f, \frac{1}{2^{n_1}}, \frac{1}{2^{n_2}} \Big)_p \asymp
\frac{1}{2^{n_1\alpha_1+n_2 \alpha_2}} \Bigg\{
\sum\limits_{\mu_1=0}^{n_1}\sum\limits_{\mu_2=0}^{n_2}
\lambda_{\mu_1,\mu_2}^2 2^{2(\mu_1\alpha_1+\mu_2 \alpha_2)}
\Bigg\}^{{1}/{2}}
 $$
 $$
+\frac{1}{2^{n_1\alpha_1}} \Bigg\{
\sum\limits_{\mu_1=0}^{n_1}\sum\limits_{\mu_2=n_2+1}^{\infty}
\lambda_{\mu_1,\mu_2}^2 2^{2\mu_1\alpha_1} \Bigg\}^{\frac{1}{2}}
+\frac{1}{2^{n_2\alpha_2}} \Bigg\{
\sum\limits_{\mu_1=n_1+1}^{\infty}\sum\limits_{\mu_2=0}^{n_2}
\lambda_{\mu_1,\mu_2}^2 2^{2\mu_2\alpha_2} \Bigg\}^{{1}/{2}}
 $$
$$
+\Bigg\{ \sum\limits_{\mu_1=n_1+1}^{\infty}\sum\limits_{\mu_2=n_2+1}^{\infty} \lambda_{\mu_1,\mu_2}^2  \Bigg\}^{{1}/{2}}.
$$
\end{theorem}
{\bf Proof  of Theorem \ref{th7.1}.}
By Theorem \ref{th4.1}, we have
$$ I := \omega_{\alpha_1,\alpha_2} \Big(f, \frac{1}{n_1},\frac{1}{n_2} \Big)_{L_p(\T^2)}  \asymp
n^{-\alpha_1}_1n_2^{-\alpha_2}
 || s_{n_1,n_2}^{(\alpha_1,\alpha_2)}(f) ||_{L_p(\T^2)}
$$
$$ +
n^{-\alpha_1}_1
 || s_{n_1,\infty}^{(\alpha_1,0)}(f-s_{\infty,n_2}(f)) ||_{L_p(\T^2)} +
n_2^{-\alpha_2}
 || s_{\infty,n_2}^{(0,\alpha_2)}(f-s_{n_1,\infty}(s)) ||_{L_p(\T^2)}+
$$

$$
+\|(f-s_{n_1,\infty}(f)-s_{\infty,n_2}(f)+s_{n_1,n_2}(f))\|_{L_p(\T^2)}.
$$
Let first $ 2 \leq p < \infty.$ Then taking into account Lemmas  \ref{l2.2} and \ref{l2.4} (A), we get

\begin{eqnarray*}
I\, &\ll&\,
n_1^{-\alpha_1}n_2^{-\alpha_2} \Big(  \sum\limits_{\nu_1=1}^{n_1} \sum\limits_{\nu_2=1}^{n_2}\rho_{\nu_1,\nu_2}^{p}
\nu_1^{(\alpha_1+1)p-2} \nu_2^{(\alpha_2+1)p-2} \Big)^{{1}/{p}}
\\&+&n_1^{-\alpha_1}  \Big(  \sum\limits_{\nu_1=1}^{n_1} \sum\limits_{\nu_2=n_2+1}^{\infty}\rho_{\nu_1,\nu_2}^{p}
\nu_1^{(\alpha_1+1)p-2} \nu_2^{p-2} \Big)^{{1}/{p}}
\\&+&n_2^{-\alpha_2} \Big(  \sum\limits_{\nu_1=n_1+1}^{\infty} \sum\limits_{\nu_2=1}^{n_2}\rho_{\nu_1,\nu_2}^{p}
\nu_1^{p-2} \nu_2^{(\alpha_2+1)p-2} \Big)^{{1}/{p}}
\\&+&\Big(  \sum\limits_{\nu_1=n_1+1}^{\infty} \sum\limits_{\nu_2=n_2+1}^{\infty}\rho_{\nu_1,\nu_2}^{p}
(\nu_1\nu_2)^{p-2} \Big)^{{1}/{p}} .
\end{eqnarray*}
Therefore, for  $p\geq 2$ we show the estimate from above in Theorem \ref{th7.1}.

If  $1<p<2$, then  H\"{o}lder's inequality gives
$$\omega_{\alpha_1,\alpha_2}\left(f,\frac1{n_1},\frac1{n_2}\right)_{L_p(\T^2)}\ll
\omega_{\alpha_1,\alpha_2}\left(f,\frac1{n_1},\frac1{n_2}\right)_{L_2(\T^2)}.$$
Thus, the estimate from above is obtained.

To prove the estimate from below, if $1<p\le 2$ , we use Lemmas  \ref{l2.2} and \ref{l2.4} (B):
\begin{eqnarray*}
&I\gg&
n_1^{-\alpha_1}n_2^{-\alpha_2} \Big(  \sum\limits_{\nu_1=1}^{n_1} \sum\limits_{\nu_2=1}^{n_2}\rho_{\nu_1,\nu_2}^{p}
\nu_1^{(\alpha_1+1)p-2} \nu_2^{(\alpha_2+1)p-2} \Big)^{{1}/{p}}
\\&+&n_1^{-\alpha_1}  \Big(  \sum\limits_{\nu_1=1}^{n_1} \sum\limits_{\nu_2=n_2+1}^{\infty}\rho_{\nu_1,\nu_2}^{p}
\nu_1^{(\alpha_1+1)p-2} \nu_2^{p-2} \Big)^{{1}/{p}}
\\&+&n_2^{-\alpha_2} \Big(  \sum\limits_{\nu_1=n_1+1}^{\infty} \sum\limits_{\nu_2=1}^{n_2}\rho_{\nu_1,\nu_2}^{p}
\nu_1^{p-2} \nu_2^{(\alpha_2+1)p-2} \Big)^{{1}/{p}}
\\&+&\Big(  \sum\limits_{\nu_1=n_1+1}^{\infty} \sum\limits_{\nu_2=n_2+1}^{\infty}\rho_{\nu_1,\nu_2}^{p}
(\nu_1\nu_2)^{p-2} \Big)^{{1}/{p}} .
\end{eqnarray*}
Let now  $2< p<\infty$.  Applying H\"{o}lder's inequality once again, we get
$$\omega_{\alpha_1,\alpha_2}\left(f,\frac1{n_1},\frac1{n_2}\right)_{L_p(\T^2)}\gg \omega_{\alpha_1,\alpha_2}\left(f,\frac1{n_1},\frac1{n_2}\right)_{L_2(\T^2)}.$$
\hfill $\Box$
\\[0.5cm]
{\bf Proof of Theorem \ref{th7.2}.}
Let us denote
\begin{eqnarray*}
I_1&:=&n_1^{-\alpha_1}n_2^{-\alpha_2}\|s^{(\alpha_1,\alpha_2)}_{n_1,n_2}(f)\|_{L_p(\T^2)},
\; \;
I_2:=n_1^{-\alpha_1}\|s^{(\alpha_1,0)}_{n_1,\infty}(f-s_{\infty,n_2}(f))\|_{L_p(\T^2)},
\\
I_3&:=&n_2^{-\alpha_2}\|s^{(0,\alpha_2)}_{\infty,n_2}(f-s_{n_1,\infty})\|_{L_p(\T^2)},
\; \;
I_4:=\|f-s_{n_1,\infty}(f)-s_{\infty,n_2}(f)+s_{n_1,n_2}\|_{L_p(\T^2)}
\end{eqnarray*}
and

\begin{eqnarray*}
A_1&:=& n_1^{-\alpha_1}n_2^{-\alpha_2} \Big(
\sum\limits_{\nu_1=1}^{n_1}
\sum\limits_{\nu_2=1}^{n_2}a_{\nu_1,\nu_2}^{p}
\nu_1^{(\alpha_1+1)p-2} \nu_2^{(\alpha_2+1)p-2}
\Big)^{{1}/{p}},
\\A_2&:=&n_1^{-\alpha_1}  \Big(  \sum\limits_{\nu_1=1}^{n_1} \sum\limits_{\nu_2=n_2+1}^{\infty}a_{\nu_1,\nu_2}^{p}
\nu_1^{(\alpha_1+1)p-2} \nu_2^{p-2} \Big)^{{1}/{p}} ,
\\A_3&:=&n_2^{-\alpha_2} \Big(  \sum\limits_{\nu_1=n_1+1}^{\infty} \sum\limits_{\nu_2=1}^{n_2}a_{\nu_1,\nu_2}^{p}
\nu_1^{p-2} \nu_2^{(\alpha_2+1)p-2} \Big)^{{1}/{p}},
\\A_4&:=&\Big(  \sum\limits_{\nu_1=n_1+1}^{\infty} \sum\limits_{\nu_2=n_2+1}^{\infty}a_{\nu_1,\nu_2}^{p}
(\nu_1\nu_2)^{p-2} \Big)^{{1}/{p}} .
\end{eqnarray*}

Let us establish the interrelation between
$I_1, I_2, I_3, I_4$ and $A_1, A_2, A_3, A_4.$ By Lemma \ref{l2.6}, we  have
\begin{equation}\label{f6.1}
I_1\asymp A_1.
\end{equation}
To estimate  $I_2$  and  $A_2$, we use
\begin{equation}\label{f6.11}
\eta_1(x,y):=\sum\limits^{n_1}_{\nu_1=1}\sum\limits^{\infty}_{\nu_2=1}a^{*}_{\nu_1,\nu_2}\cos{\nu_1 x}\cos{\nu_2 y}
\end{equation}
and
\begin{equation}\label{f6.111}
\eta_2(x,y):=\sum\limits^{n_1}_{\nu_1=1}\sum\limits^{n_2}_{\nu_2=1}a^{*}_{\nu_1,\nu_2}\cos{\nu_1
x}\cos{\nu_2 y},
\end{equation}
where
\[a^{*}_{\nu_1,\nu_2}:=
\begin{cases}
a_{\nu_1,\nu_2} \; \textnormal{for} \; 1 \le \nu_1 \le n_1, \nu_2>n_2 \\
a_{\nu_1,n_2} \;  \textnormal{for} \; 1 \le \nu_1 \le n_1, 1 \le \nu_2
\le n_2.
\end{cases}
\]
Then
\begin{equation}\label{f6.2}
s^{(\alpha_1,0)}_{n_1,\infty}(f-s_{\infty,n_2}(f))=(\eta_1-\eta_2)^{(\alpha_1,0)}.
\end{equation}
Let us first estimate   $I_2$. It is clear that
$$
I_2\ll n_1^{-\alpha_1}\, \left\{\|\eta_1^{(\alpha_1,0)}\|_{L_p(\T^2)}+\|\eta_2^{(\alpha_2,0)}\|_{L_p(\T^2)}\right\}.
$$
Now Lemma \ref{l2.6} implies that
\begin{eqnarray}\nonumber
I_2&\le & n_1^{-\alpha_1} \Big(  \sum\limits_{\nu_1=1}^{n_1}
\sum\limits_{\nu_2=1}^{\infty}(a^*_{\nu_1,\nu_2})^{p}
\nu_1^{(\alpha_1+1)p-2} \nu_2^{p-2} \Big)^{{1}/{p}}
\\
\nonumber&+&n_1^{-\alpha_1}  \Big(  \sum\limits_{\nu_1=1}^{n_1} \sum\limits_{\nu_2=1}^{n_2}(a^*_{\nu_1,\nu_2})^{p}
\nu_1^{(\alpha_1+1)p-2} \nu_2^{p-2} \Big)^{{1}/{p}}
\\\nonumber&\ll &n_1^{-\alpha_1} \Big(  \sum\limits_{\nu_1=1}^{n_1} \sum\limits_{\nu_2=n_2+1}^{\infty}a_{\nu_1,\nu_2}^{p}
\nu_1^{(\alpha_1+1)p-2} \nu_2^{p-2} \Big)^{{1}/{p}}
\\\nonumber&+&n_1^{-\alpha_1} \Big(  \sum\limits_{\nu_1=1}^{n_1}a_{\nu_1,n_2}^{p}
\nu_1^{(\alpha_1+1)p-2} n_2^{p-1} \Big)^{{1}/{p}}
\\
\label{f6.3}
&\ll&  A_2+n_1^{-\alpha_1}n_2^{-\alpha_2} \Big(
\sum\limits_{\nu_1=1}^{n_1}
\sum\limits_{\nu_2=
\floor{\frac{n_2}{2}}+1}^{n_2}a_{\nu_1,\nu_2}^{p}
\nu_1^{(\alpha_1+1)p-2} \nu_2^{(\alpha_2+1)p-2}
\Big)^{{1}/{p}}\ll  A_2+A_1.
\end{eqnarray}
Estimating $A_2$ from above, we write
$$
A_2\ll n_1^{-\alpha_1} \Big(  \sum\limits_{\nu_1=1}^{n_1}
\sum\limits_{\nu_2=1}^{\infty}(a^*_{\nu_1,\nu_2})^{p}
\nu_1^{(\alpha_1+1)p-2} \nu_2^{p-2} \Big)^{{1}/{p}}.
$$
 Lemma \ref{l2.6} and formulas (\ref{f6.11})-(\ref{f6.111}) imply
$$
A_2\ll n_1^{-\alpha_1} \|\eta^{(\alpha_1,0)}_{1}\|_{L_p(\T^2)}\ll I_2+n_1^{-\alpha_1}\|\eta^{(\alpha_1,0)}_{2}\|_{L_p(\T^2)}
$$
and
\begin{eqnarray*}
A_2&\ll&
I_2+n_1^{-\alpha_1}\Big(  \sum\limits_{\nu_1=1}^{n_1}
\sum\limits_{\nu_2=1}^{n_2}(a^*_{\nu_1,\nu_2})^{p}
\nu_1^{(\alpha_1+1)p-2} \nu_2^{p-2} \Big)^{{1}/{p}}
\\
&\ll& I_2+n_1^{-\alpha_1}\Big(  \sum\limits_{\nu_1=1}^{n_1}
a_{\nu_1,n_2}^{p} \nu_1^{(\alpha_1+1)p-2} n_2^{p-1}
\Big)^{{1}/{p}}\ll I_2+A_1.
\end{eqnarray*}
Now we use (\ref{f6.1}) to get
\begin{equation}\label{f6.4}
A_2\ll I_2+I_1.
\end{equation}
Similarly, we get
\begin{equation}\label{f6.5}
I_3\ll A_1+A_3
\end{equation}
and
\begin{equation}\label{f6.6}
A_3\ll I_3+I_1.
\end{equation}
Thus, it is shown that
\begin{equation}\label{ff6.6}
A_1+A_2+A_3\asymp I_1+I_2+I_3.
\end{equation}

Now we consider
$ B:=f-s_{n_1,\infty}-s_{\infty,n_1}+s_{n_1,n_2}$
and
$$
\varphi_1(x,y):=\sum\limits^{\infty}_{\nu_1=1}\sum\limits^{\infty}_{\nu_2=1}b_{\nu_1,\nu_2}\cos{\nu_1 x}\cos{\nu_2 y},
$$
$$
\varphi_2(x,y):=\sum\limits^{n_1}_{\nu_1=1}\sum\limits^{\infty}_{\nu_2=n_2+1}b_{\nu_1,\nu_2}\cos{\nu_1
x}\cos{\nu_2 y},
$$
$$
\varphi_3(x,y):=\sum\limits^{\infty}_{\nu_1=n_1+1}\sum\limits^{n_2}_{\nu_2=1}b_{\nu_1,\nu_2}\cos{\nu_1 x}\cos{\nu_2 y},
$$
$$
\varphi_4(x,y):=\sum\limits^{n_1}_{\nu_1=1}\sum\limits^{n_2}_{\nu_2=1}b_{\nu_1,\nu_2}\cos{\nu_1 x}\cos{\nu_2 y},
$$
where
\[b_{\nu_1,\nu_2}:=
\begin{cases}
a_{\nu_1, \nu_2} \qquad \textnormal{for} \qquad  \nu_1>n_1, \nu_2>n_2,\\
a_{\nu_1, n_2} \qquad  \textnormal{for}  \qquad \nu_1>n_1,  1 \le \nu_2 \le n_2,\\
 a_{n_1,\nu_2} \qquad \textnormal{for} \qquad 1 \le \nu_1 \le n_1, \nu_2>n_2, \\
 a_{n_1, n_2} \qquad  \textnormal{for}  \qquad  1 \le \nu_1 \le n_1,1 \le \nu_2 \le n_2.
\end{cases}
\]
Then
\begin{equation}\label{f6.7}
B=\varphi_1-\varphi_2-\varphi_3+\varphi_4.
\end{equation}

We first estimate
$$
I_4=\|B\|_{L_p(\T^2)}\ll \|\varphi_1\|_{L_p(\T^2)}+\|\varphi_2\|_{L_p(\T^2)}+\|\varphi_3\|_{L_p(\T^2)}+\|\varphi_4\|_{L_p(\T^2)}.
$$
Lemma \ref{l2.6} yields
$$
\|\varphi_1\|_{L_p(\T^2)} \asymp \sum\limits^{\infty}_{\nu_1=1}\sum\limits^{\infty}_{\nu_2=1}b^p_{\nu_1,\nu_2}(\nu_1\nu_2)^{p-2}
\asymp \sum\limits^{\infty}_{\nu_1=n_1+1}\sum\limits^{\infty}_{\nu_2=n_2+1}a^p_{\nu_1,\nu_2}(\nu_1\nu_2)^{p-2}
$$
$$
+\sum\limits^{\infty}_{\nu_1=n_1+1}a^p_{\nu_1,n_2}\nu_1^{p-2}n_2^{p-1}+\sum\limits^{\infty}_{\nu_2=n_2+1}a^p_{n_1,\nu_2}\nu_2^{p-2}n_1^{p-1}
+a^p_{n_1,n_2}(n_1n_2)^{p-1}=: J_1+J_2+J_3+J_4.
$$
To estimate $\|\varphi_2\|_{L_p(\T^2)}$, we consider
$$
\varphi_{21}(x,y):=\sum\limits^{n_1}_{\nu_1=1}\sum\limits^{\infty}_{\nu_2=1}b^*_{\nu_1,\nu_2}\cos{\nu_1 x}\cos{\nu_2 y},
$$
$$
\varphi_{22}(x,y):=\sum\limits^{n_1}_{\nu_1=1}\sum\limits^{n_2}_{\nu_2=1}b^*_{\nu_1,\nu_2}\cos{\nu_1 x}\cos{\nu_2 y},
$$
where
\[b_{\nu_1,\nu_2}^*:=
\begin{cases}
b_{\nu_1,\nu_2} \; \textnormal{for}  \;   1 \le \nu_1 \le n_1, \nu_2>n_2,\\
b_{\nu_1,n_2} \;  \textnormal{for}  \; 1 \le \nu_1 \le n_1,  1 \le \nu_2 \le n_2.\\
\end{cases}
\]
Then   $\|\varphi_2\|_{L_p(\T^2)}\ll  \|\varphi_{21}\|_{L_p(\T^2)}+
\|\varphi_{22}\|_{L_p(\T^2)}.$ Using Lemma \ref{l2.6}, we have
\begin{eqnarray*}
\|\varphi_{21}\|_{L_p(\T^2)} &\asymp& \sum\limits^{n_1}_{\nu_1=1}\sum\limits^{\infty}_{\nu_2=1}(b^*_{\nu_1,\nu_2})^p(\nu_1\nu_2)^{p-2}\ll
 \sum\limits^{n_1}_{\nu_1=1}\sum\limits^{\infty}_{\nu_2=n_2+1}b^p_{\nu_1,\nu_2}(\nu_1\nu_2)^{p-2}
\\
&+& \sum\limits^{n_1}_{\nu_1=1}b^p_{\nu_1,n_2}\nu_1^{p-2}n_2^{p-1}\ll
 \sum\limits^{\infty}_{\nu_2=n_2+1}a^p_{n_1,\nu_2}\nu_2^{p-2}n_1^{p-1}+a_{n_1,n_2}^p(n_1n_2)^{p-1}
 =: J_3+J_4,
\\
\|\varphi_{22}\|_{L_p(\T^2)} &\asymp& \sum\limits^{n_1}_{\nu_1=1}\sum\limits^{n_2}_{\nu_2=1}(b^*_{\nu_1,\nu_2})^p(\nu_1\nu_2)^{p-2}\ll
a_{n_1,n_2}^p(n_1n_2)^{p-1}\ll  J_4.
\end{eqnarray*}
Therefore, $\|\varphi_2\|_{L_p(\T^2)}^p\ll  J_3+J_4.$
Similarly, $\|\varphi_3\|_{L_p(\T^2)}^p\ll  J_2+J_4$ and $\|\varphi_4\|_{L_p(\T^2)}^p\ll  J_4.$
Combining these estimates, we get
\begin{eqnarray*}
I_4 &\ll& \left\{ \sum\limits^{\infty}_{\nu_1=n_1+1}\sum\limits^{\infty}_{\nu_2=n_2+1} a^p_{\nu_1,\nu_2}(\nu_1\nu_2)^{p-2} \right\}^{1/p}+
\left\{ \sum\limits^{\infty}_{\nu_1=n_1+1} a^p_{\nu_1,n_2}\nu_1^{p-2} \right\}^{1/p} n_2^{1-\frac1p}
\\&+&\left\{ \sum\limits^{\infty}_{\nu_2=n_2+1} a^p_{n_1,\nu_2}\nu_2^{p-2} \right\}^{1/p} n_1^{1-\frac1p}+a_{n_1,n_2} (n_1n_2)^{1-\frac1p}.
\end{eqnarray*}
Note that the latter  inequality holds also if $n_1=0$ and/or
 $n_2=0$.
It is easy to verify that the last estimate implies
\begin{equation}\label{f6.9}
I_4 \ll A_1+A_2+A_3+A_4.
\end{equation}

Let us estimate $A_4$ from above. It is clear that
$$
A_4\ll \left\{ \sum\limits^{\infty}_{\nu_1=1}\sum\limits^{\infty}_{\nu_2=1} b^p_{\nu_1,\nu_2}(\nu_1\nu_2)^{p-2} \right\}^{1/p}.
$$
Lemma \ref{l2.6} yields
$A_4\ll  \|\varphi_1\|_{L_p(\T^2)}$ and, by
(\ref{f6.7}),  we get
$
A_4\ll  \|B\|_{L_p(\T^2)}+\|\varphi_2\|_{L_p(\T^2)}+\|\varphi_3\|_{L_p(\T^2)}+\|\varphi_4\|_{L_p(\T^2)}.
$
Then
\begin{eqnarray*}
A_4&\le& I_4+J_2+J_3+J_4\le I_4+\left\{\sum\limits^{\infty}_{\nu_1=1}a^p_{\nu_1,n_2}\nu_1^{p-2} \right\}^{1/p} n_{2}^{1-\frac1p}
\\
&+&
\left\{\sum\limits^{\infty}_{\nu_2=n_2+1}a^p_{n_1,\nu_2}\nu_2^{p-2}
\right\}^{1/p}
n_{1}^{1-\frac1p}+a_{n_1,n_2}(n_1n_2)^{1-\frac1p}\ll
I_4+A_2+A_3+A_1.
\end{eqnarray*}
Using estimates  (\ref{ff6.6}) and (\ref{f6.9}), we have
\begin{equation}\label{f6.10}
A_4 \le I_1+I_2+I_3+I_4.
\end{equation}
This, (\ref{ff6.6}) and (\ref{f6.9}) give us
\begin{equation*}
  I_1+I_2+I_3+I_4\asymp A_1+A_2+A_3+A_4 .
\end{equation*}

Since Theorem \ref{th4.1} implies
$$
\omega_{\alpha_1,\alpha_2}\left(f,\frac1{n_1},\frac1{n_2}\right)_{L_p(\T^2)}^p \asymp I_1+I_2+I_3+I_4,
$$
the proof of Theorem \ref{th7.2} is now complete.
\hfill $\Box$
\\[0.5cm]
{\bf Proof of Theorem \ref{th7.3}.}
Theorem \ref{th4.1} implies
\begin{eqnarray*}
I&:=&\omega_{\alpha_1,\alpha_2}\left(f,\frac1{2^{n_1}},\frac1{2^{n_2}}\right)_{L_p(\T^2)}^p  \asymp 2^{-n_1\alpha_1}2^{-n_2\alpha_2}
\|s_{2^{n_1},2^{n_2}}^{(\alpha_1,\alpha_2)}\|_{L_p(\T^2)}\\
&+&2^{-n_1\alpha_1}\|s_{2^{n_1},\infty}^{(\alpha_1,0)}(f-s_{\infty,2^{n_2}}(f))\|_{L_p(\T^2)}+
2^{-n_2\alpha_2}\|s_{\infty,2^{n_2}}^{(0,\alpha_2)}(f-s_{2^{n_1},\infty}(f))\|_{L_p(\T^2)}\\
&+&\|f-s_{2^{n_1},\infty}(f)-s_{\infty,2^{n_2}}(f)+s_{2^{n_1},2^{n_2}}(f)\|_{L_p(\T^2)}.
\end{eqnarray*}
Lemmas  \ref{l2.2} and \ref{l2.7} yield
\begin{eqnarray*}
I &\asymp& 2^{-n_1\alpha_1-n_2\alpha_2}\left\{ \sum\limits^{n_1}_{\mu_1=0}
\sum\limits^{n_2}_{\mu_2=0} \lambda_{\mu_1,\mu_2}^2 2^{2(\alpha_1\mu_1+\alpha_2\mu_2)} \right\}^{1/2}+
2^{-n_1\alpha_1} \left\{ \sum\limits^{n_1}_{\mu_1=0}
\sum\limits^{\infty}_{\mu_2=n_1+1} \lambda_{\mu_1,\mu_2}^2 2^{2\alpha_1\mu_1} \right\}^{1/2}\\
&+&2^{-n_2\alpha_2}\left\{ \sum\limits^{\infty}_{\mu_1=n_1+1}
\sum\limits^{n_2}_{\mu_2=0} \lambda_{\mu_1,\mu_2}^2 2^{2\alpha_2\mu_2} \right\}^{1/2}+
\left\{ \sum\limits^{\infty}_{\mu_1=n_1+1}
\sum\limits^{\infty}_{\mu_2=n_2+1} \lambda_{\mu_1,\mu_2}^2  \right\}^{1/2},
\end{eqnarray*}
i.e., the required equivalence.
\hfill $\Box$

Theorem \ref{th7.1} is stated in the paper \cite{23}, while Theorems  \ref{th7.2}-\ref{th7.3} can be found in the papers \cite{23,21,22}.
The one-dimensional version of Theorem \ref{th7.2} for functions with general monotone coefficients is proved in \cite{gorb}.
 Theorem \ref{th7.3} is given in \cite{berisha} in the case of integer order moduli.

\vskip 1.5cm

\sect{The mixed moduli of smoothness of $L_p$-functions and their derivatives}
Let $1<p<\infty$.
We start this section with two well-known relations for the one-dimensional modulus of smoothness:
for $f,$ $f^{(k)}\in L_p(\T),$
 we have the estimate (see
\cite[p.~46]{de-lo})
$$
\omega_{r+k}(f,\delta)_{L_p(\T)}\ll \delta^k\omega_{r} (f^{(k)},\delta)_{L_p(\T)},
\qquad\text{\rm where}\qquad r,k\in\N
$$
and its weak inverse (see \cite[p.~178]{de-lo})
$$
\omega_{r}(f^{(k)},\delta)_{L_p(\T)} \ll \int^\delta_0 \;\frac{\omega
_{r+k}(f,u)_{L_p(\T)}}{u^{k+1}}\;du, \qquad\text{\rm where}\qquad r,k\in\N.
$$
The analogues of these  inequalities for the mixed moduli of smoothness are given by
\begin{eqnarray*}
\delta_1^{-r_1}\delta_2^{-r_2}
\omega_{\beta_1+r_1, \beta_2+r_2}
(f, \delta_1,\delta_2)_{L_p(\T^2)}&\ll&
\omega_{\beta_1, \beta_2}
(f^{(r_1, r_2)}, \delta_1, \delta_2)_{L_p(\T^2)}
\\
&\ll&
 \int\limits_0^{\delta_1} \int\limits_0^{\delta_2}
t_1^{-r_1}t_2^{-r_2} \omega_{\beta_1+r_1,\beta_2+r_2}(f, t_1,t_2)_{L_p(\T^2)}
\frac{\dd t_1}{t_1} \frac{\dd t_2}{t_2},
\end{eqnarray*}
where $\beta_j, r_j>0,\quad j=1,2$ (see Theorem \ref{th6.1}; properties (9) and (10)).

Now we provide a generalization of these estimates (see \cite{deriv-st}). Note that this idea goes back to the papers by Besov
\cite{be} and Marcinkiewicz \cite{ma}. We also remark that the right-hand side inequality was done for non-mixed moduli in
\cite{dit-studia, trebels}. Moreover, a similar estimate for mixed moduli of smoothness of functions defined on $\mathbb{R}^d$ was proved in \cite{tom}.

\begin{theorem}\label{T10}
Let $f \in L_p^0(\T^2), 1 < p  < \infty,$
$ \theta := \min(2,p), \tau:=\max(2,p) $,
$\beta_1, \beta_2, r_1, r_2>0$, and $\delta_1,\delta_2\in (0,\frac12).$
 If
\begin{equation*} \label{th01}
\left( \int\limits^{1}_{0}\int\limits^{1}_{0}
t_1^{-r_1\theta-1} t_2^{-r_2\theta-1}  \omega^\theta
_{r_1+\beta_1, r_2+\beta_2}(f,t_1, t_2)_{L_p(\T^2)} \,\dd t_1\,  \dd t_2
\right)^{{1}/{\theta}} \,<\infty,
\end{equation*}
then  $f$ has a mixed derivative $f^{(r_1,r_2)}\in {L_p(\T^2)}$ in the sense of Weyl and
\begin{equation} \label{th02}
\omega_{\beta_1, \beta_2} (f^{(r_1,r_2)},\delta_1, \delta_2)_{L_p(\T^2)}
\ll
\left\{
\int\limits_{0}^{\delta_1}\int\limits_{0}^{\delta_2}
t_1^{-r_1\theta-1} t_2^{-r_2\theta-1}
 \omega^\theta_{r_1+\beta_1, r_2+\beta_2}(f,t_1, t_2)_{L_p(\T^2)} \,\dd t_1\,  \dd t_2
\right\}^{1/{\theta}}.
\end{equation}
If  $f\in {L_p(\T^2)}$ has a mixed derivative $f^{(r_1,r_2)}\in {L_p(\T^2)}$ in the sense of Weyl,
then
$$
\left\{
\int\limits_{0}^{\delta_1}\int\limits_{0}^{\delta_2}
t_1^{-r_1\tau-1} t_2^{-r_2\tau-1}
 \omega^\tau_{r_1+\beta_1, r_2+\beta_2}(f,t_1, t_2)_{L_p(\T^2)} \,\dd t_1\,  \dd t_2
\right\}^{1/{\tau}}
\ll
 \omega_{\beta_1, \beta_2}
(f^{(r_1,r_2)},\delta_1, \delta_2)_{L_p(\T^2)}.
$$
\end{theorem}

Now we deal with  the function classes $M_p$ and $\Lambda_p$  defined in Section 2.4.
\begin{theorem}\label{T11}
Let $f \in M_p, 1 < p < \infty,$
$\beta_1, \beta_2, r_1, r_2>0,$ and $\delta_1,\delta_2\in (0,\frac12)$. Then
\begin{equation*}
\omega_{\beta_1, \beta_2} (f^{(r_1,r_2)},\delta_1, \delta_2)_{L_p(\T^2)}
\asymp
\left\{
\int\limits_{0}^{\delta_1}\int\limits_{0}^{\delta_2}
t_1^{-r_1p-1} t_2^{-r_2p-1}
 \omega^p_{r_1+\beta_1, r_2+\beta_2}(f,t_1, t_2)_{L_p(\T^2)} \,\dd t_1\,  \dd t_2
\right\}^{1/{p}}.
\end{equation*}
\end{theorem}

\begin{theorem}\label{T12}
Let $f \in \Lambda_p, 1 < p < \infty,$
$\beta_1, \beta_2, r_1, r_2>0,$ and $\delta_1,\delta_2\in (0,\frac12)$. Then
\begin{equation*}
\omega_{\beta_1, \beta_2} (f^{(r_1,r_2)},\delta_1, \delta_2)_{L_p(\T^2)}
\asymp
\left\{
\int\limits_{0}^{\delta_1}\int\limits_{0}^{\delta_2}
t_1^{-2r_1-1} t_2^{-2r_2-1}
 \omega^2_{r_1+\beta_1, r_2+\beta_2}(f,t_1, t_2)_{L_p(\T^2)} \,\dd t_1\,  \dd t_2
\right\}^{1/{2}}.
\end{equation*}
\end{theorem}
{\bf Proof of Theorem \ref{T10}.}
We choose integers $n_1, n_2$ such that $2^{-n_i}< \delta_i \le
2^{-n_i+1},\,i=1,2$. By Lemmas \ref{l2.2}, \ref{l2.3}, and Theorem \ref{th4.1}, we have
\begin{eqnarray*}
\omega_{\beta_1, \beta_2}
(f^{(r_1,r_2)},\delta_1,\delta_2)_{{p}} &\ll&
 \left\|
 \sum\limits_{\nu_1=2^{n_1}+1}^{\infty}\sum\limits_{\nu_2=n_2+1}^{\infty}
\nu_1^{r_1} \nu_2^{r_2} A_{\nu_1 \nu_2} (x_1, x_2)
\right\|_{{p}}
\\
&+&
 2^{-n_1\beta_1}
 \left\|
\sum\limits_{\nu_1=1}^{2^{n_1}}
 \sum\limits_{\nu_2=2^{n_2}+1}^{\infty}
\nu_1^{r_1+\beta_1} \nu_2^{r_2} A_{\nu_1 \nu_2} (x_1, x_2)
\right\|_{{p}}
\\
&+&
 2^{-n_2\beta_2}
 \left\|
  \sum\limits_{\nu_1=2^{n_1}+1}^{\infty}
 \sum\limits_{\nu_2=1}^{2^{n_2}}
\nu_1^{r_1} \nu_2^{r_2+\beta_2} A_{\nu_1 \nu_2} (x_1, x_2)
\right\|_{{p}}
\\
&+& 2^{-n_1\beta_1}   2^{-n_2\beta_2}
 \left\|
\sum\limits_{\nu_1=1}^{2^{n_1}} \sum\limits_{\nu_2=1}^{2^{n_2}}
\nu_1^{r_1+\beta_1} \nu_2^{r_2+\beta_2} A_{\nu_1 \nu_2} (x_1, x_2)
\right\|_{{p}}
\\
&=:& I_1+I_2+I_3+I_4.
\end{eqnarray*}
where $A_{\nu_1 \nu_2}$ is given by (\ref{e2.1}).
Estimating $I_2$ as follows
$$
I_2 \ll 2^{-n_1\beta_1}
 \left\{
 \int\limits_{0}^{2\pi}
 \int\limits_{0}^{2\pi} \left[ \sum\limits_{\nu_1=0}^{n_1}
\sum\limits_{\nu_2=n_2+1}^{\infty} 2^{2\nu_1 (r_1+\beta_1) + 2\nu_2
r_2} \triangle_{\nu_1 \nu_2}^2 \right]^{p/2} \dd x_1
 \dd x_2 \right\}^{1/p},
$$
we use
\begin{equation}\label{vsp1}
2^{2\nu_j r_j} \asymp \left( \sum\limits_{\xi_j=n_j+1}^{\nu_j}
2^{\xi_j r_j \theta} \right)^{2/\theta},  \quad j=1,2,
\end{equation}
and Minkowski's inequality. Hence,
$$%\begin{multline*}
I_{2}\ll 2^{-n_1\beta_1}
 \left\{
 \int\limits_{0}^{2\pi}
 \int\limits_{0}^{2\pi}
\left[ \sum\limits_{\xi_2=n_2+1}^{\infty} 2^{\xi_2 r_2\theta} \left(
\sum\limits_{\nu_2=\xi_2}^{\infty} \sum\limits_{\nu_1=0}^{n_1}
2^{2\nu_1 (r_1+\beta_1) } \triangle_{\nu_1 \nu_2}^2
\right)^{\theta/2} \right]^{p/\theta} \dd x_1
 \dd x_2 \right\}^{1/p}.
$$%\end{multline*}
Applying again Minkowski's inequality for
 sums and integrals, we have
$$
I_{2}\ll
 \left\{
2^{-n_1\beta_1\theta} \sum\limits_{\xi_2=n_2+1}^{\infty} 2^{\xi_2
r_2\theta} \left(
 \int\limits_{0}^{2\pi}
 \int\limits_{0}^{2\pi}
\left<
\sum\limits_{\nu_1=0}^{n_1}
\sum\limits_{\nu_2=\xi_2}^{\infty} 2^{2\nu_1 (r_1+\beta_1)}
\triangle_{\nu_1 \nu_2}^2 \right>^{p/2}
\dd x_1
\dd x_2 \right)^{\theta/p}
\right\}^{1/\theta}.
$$
Using now
\begin{equation}\label{vsp2}
2^{-n_j \beta_j\theta} \asymp \sum\limits_{\xi_j=n_j+1}^{\infty}
2^{-\xi_j \beta_j\theta},
 \quad j=1,2,
\end{equation}
as well as Lemmas \ref{l2.2}, \ref{l2.3},  and Theorem \ref{th4.1}, we get
\begin{eqnarray*}
I_{2}&\ll&
 \left\{
\sum\limits_{\xi_1=n_1}^{\infty} 2^{-\xi_1 \beta_1\theta}
\sum\limits_{\xi_2=n_2}^{\infty} 2^{\xi_2 r_2\theta}
 \left\|
\sum\limits_{\nu_1=1}^{ 2^{\xi_1} } \sum\limits_{ \nu_2=2^{\xi_2}+1
}^{\infty} \nu_1^{r_1+\beta_1} A_{\nu_1 \nu_2} (x_1, x_2)
\right\|_{{p}} ^\theta \right\}^{1/\theta}
\\
&\ll&
 \left\{
\sum\limits_{\xi_1=n_1}^{\infty} \sum\limits_{\xi_2=n_2}^{\infty}
2^{\xi_1 r_1\theta} 2^{\xi_2 r_2\theta} \omega^\theta_{r_1+\beta_1,
r_2+\beta_2}\left(f,\frac{1}{2^{\xi_1}},
\frac{1}{2^{\xi_2}}\right)_{{p}} \right\}^{1/\theta}.
\end{eqnarray*}
For the estimate of $I_{3}$ we use the same reasoning:
\begin{eqnarray*}
I_{3}&\ll&
 \left\{
\sum\limits_{\xi_1=n_1}^{\infty} \sum\limits_{\xi_2=n_2}^{\infty}
2^{\xi_1 r_1\theta} 2^{-\xi_2 \beta_2\theta}
 \left\|
\sum\limits_{\nu_1=2^{\xi_1}+1}^{\infty}
\sum\limits_{\nu_2=1}^{2^{\xi_2}} \nu_2^{r_2+\beta_2} A_{\nu_1
\nu_2} (x_1, x_2) \right\|_{{p}} ^\theta
\right\}^{1/\theta}
\\
&\ll&
 \left\{
\sum\limits_{\xi_1=n_1}^{\infty} \sum\limits_{\xi_2=n_2}^{\infty}
2^{\xi_1 r_1\theta} 2^{\xi_2 r_2\theta} \omega^\theta_{r_1+\beta_1,
r_2+\beta_2}\left(f,\frac{1}{2^{\xi_1}},
\frac{1}{2^{\xi_2}}\right)_{
{p}} \right\}^{1/\theta}.
\end{eqnarray*}
Now we estimate $I_{1}$:
\begin{eqnarray*}
I_1 &\ll &
 \left\{
 \int\limits_{0}^{2\pi}
 \int\limits_{0}^{2\pi} \left[
\sum\limits_{\nu_1=n_1+1}^{\infty}
\sum\limits_{\nu_2=n_2+1}^{\infty} 2^{2\nu_1 r_1 + 2\nu_2 r_2}
\triangle_{\nu_1 \nu_2}^2 \right]^{p/2} \dd x_1
 \dd x_2 \right\}^{1/p}
\\
&\ll&
 \left\{
\sum\limits_{\xi_2=n_2}^{\infty} 2^{\xi_2 r_2\theta} \left(
 \int\limits_{0}^{2\pi}
 \int\limits_{0}^{2\pi}
\left< \sum\limits_{\nu_1=n_1+1}^{\infty}
\sum\limits_{\nu_2=\xi_2+1}^{\infty} 2^{2 \nu_1 r_1 } \triangle_{\nu_1
\nu_2}^2 \right>^{p/2} \dd x_1  \dd x_2
\right)^{\theta/p} \right\}^{1/\theta}
\\
&\ll&
 \left\{
\sum\limits_{\xi_1=n_1}^{\infty} \sum\limits_{\xi_2=n_2}^{\infty}
2^{\xi_1 r_1\theta} 2^{\xi_2 r_2\theta} \left(
 \int\limits_{0}^{2\pi}
 \int\limits_{0}^{2\pi}
\left< \sum\limits_{\nu_1=\xi_1+1}^{\infty}
\sum\limits_{\nu_2=\xi_2+1}^{\infty} \triangle_{\nu_1 \nu_2}^2
\right>^{p/2} \dd x_1  \dd x_2
\right)^{\theta/p} \right\}^{1/\theta}
\\
&\ll&
 \left\{
\sum\limits_{\xi_1=n_1}^{\infty} \sum\limits_{\xi_2=n_2}^{\infty}
2^{\xi_1 r_1\theta} 2^{\xi_2 r_2\theta}
 \left\|
\sum\limits_{\nu_1=2^{\xi_1}+1}^{\infty}
\sum\limits_{\nu_2=2^{\xi_2}+1}^{\infty} A_{\nu_1 \nu_2} (x_1, x_2)
\right\|_{{p}} ^\theta \right\}^{1/\theta}
\\
&\ll&
 \left\{
\sum\limits_{\xi_1=n_1}^{\infty} \sum\limits_{\xi_2=n_2}^{\infty}
2^{\xi_1 r_1\theta} 2^{\xi_2 r_2\theta} \omega^\theta_{r_1+\beta_1,
r_2+\beta_2}\left(f,\frac{1}{2^{\xi_1}},
\frac{1}{2^{\xi_2}}\right)_{{p}} \right\}^{1/\theta}.
\end{eqnarray*}
To estimate  $I_4$, we use the expression (\ref{vsp2}) twice:
\begin{eqnarray*}
I_4& \ll & 2^{-n_1\beta_1} 2^{-n_2\beta_2}
 \left\{
 \int\limits_{0}^{2\pi}
 \int\limits_{0}^{2\pi} \left[ \sum\limits_{\nu_1=0}^{n_1}
\sum\limits_{\nu_2=0}^{n_2} 2^{2\nu_1 (r_1+\beta_1) + 2\nu_2
(r_2+\beta_2)} \triangle_{\nu_1 \nu_2}^2 \right]^{p/2} \dd x_1
 \dd x_2 \right\}^{1/p}
\\
&\ll&
 \left\{
\sum\limits_{\xi_1=n_1}^{\infty} \sum\limits_{\xi_2=n_2}^{\infty}
2^{-\xi_1 \beta_1\theta} 2^{-\xi_2 \beta_2\theta}
 \left\|
\sum\limits_{\nu_1=1}^{ 2^{\xi_1} } \sum\limits_{\nu_2=1}^{
2^{\xi_2} } \nu_1^{r_1+\beta_1} \nu_2^{r_2+\beta_2} A_{\nu_1 \nu_2}
(x_1, x_2) \right\|_{{p}} ^\theta \right\}^{1/\theta}
\\
&\ll&
 \left\{
\sum\limits_{\xi_1=n_1}^{\infty} \sum\limits_{\xi_2=n_2}^{\infty}
2^{\xi_1 r_1\theta} 2^{\xi_2 r_2\theta} \omega^\theta_{r_1+\beta_1,
r_2+\beta_2}\left(f,\frac{1}{2^{\xi_1}},
\frac{1}{2^{\xi_2}}\right)_{{p}} \right\}^{1/\theta}.
\end{eqnarray*}
Collecting estimates for  $I_j, j=1,2,3,4,$ we finally get
$$
\omega_{\beta_1, \beta_2}
(f^{(r_1,r_2)},\delta_1,\delta_2)_{{p}} \ll
 \left\{
\sum\limits_{\xi_1=n_1+1}^{\infty}
\sum\limits_{\xi_2=n_2+1}^{\infty} 2^{\xi_1 r_1\theta} 2^{\xi_2
r_2\theta} \omega^\theta_{r_1+\beta_1,
r_2+\beta_2}\left(f,\frac{1}{2^{\xi_1}},
\frac{1}{2^{\xi_2}}\right)_{{p}} \right\}^{1/\theta}
$$
and (\ref{th02}) follows.

Now let us prove the reverse inequality. We denote
\begin{eqnarray*}
K
&:=&
\int\limits_{0}^{\delta_1}\int\limits_{0}^{\delta_2}
t_1^{-r_1\tau-1} t_2^{-r_2\tau-1}
 \omega^\tau_{r_1+\beta_1, r_2+\beta_2}(f,t_1, t_2)_{L_p(\T^2)} \,\dd t_1\,  \dd t_2
\\
&\ll&
 \sum\limits_{\xi_1=n_1}^{\infty}
\sum\limits_{\xi_2=n_2}^{\infty} 2^{\xi_1 r_1\tau} 2^{\xi_2 r_2\tau}
\omega^\tau_{r_1+\beta_1,
r_2+\beta_2}\left(f_,\frac{1}{2^{\xi_1}},
\frac{1}{2^{\xi_2}}\right)_{{p}}.
\end{eqnarray*}
Using  Lemmas \ref{l2.2}, \ref{l2.3},  and Theorem \ref{th4.1},
\begin{eqnarray*}
K
&\ll&
\sum\limits_{\xi_1=n_1}^{\infty} \sum\limits_{\xi_2=n_2}^{\infty}
2^{\xi_1 r_1\tau} 2^{\xi_2 r_2\tau}
 \left\|
\sum\limits_{\nu_1=2^{\xi_1}+1}^{ \infty }
\sum\limits_{\nu_2=2^{\xi_2}+1}^{ \infty } A_{\nu_1 \nu_2} (x_1,
x_2) \right\|_{{p}} ^\tau
\\
&+& \sum\limits_{\xi_1=n_1}^{\infty}
\sum\limits_{\xi_2=n_2}^{\infty} 2^{\xi_1 r_1\tau} 2^{-\xi_2
\beta_2\tau}
 \left\|
\sum\limits_{\nu_1=2^{\xi_1}+1}^{ \infty } \sum\limits_{\nu_2=1}^{
2^{\xi_2} } \nu_2^{r_2+\beta_2} A_{\nu_1 \nu_2} (x_1, x_2)
\right\|_{{p}} ^\tau
\\
&+& \sum\limits_{\xi_1=n_1}^{\infty}
\sum\limits_{\xi_2=n_2}^{\infty} 2^{-\xi_1 \beta_1\tau} 2^{\xi_2
r_2\tau}
 \left\|
\sum\limits_{\nu_1=1}^{ 2^{\xi_1} }
\sum\limits_{\nu_2=2^{\xi_2}+1}^{ \infty } \nu_1^{r_1+\beta_1}
A_{\nu_1 \nu_2} (x_1, x_2) \right\|_{{p}} ^\tau
\\
&+& \sum\limits_{\xi_1=n_1}^{\infty}
\sum\limits_{\xi_2=n_2}^{\infty} 2^{-\xi_1 \beta_1\tau} 2^{-\xi_2
\beta_2\tau}
 \left\|
\sum\limits_{\nu_1=1}^{ 2^{\xi_1} } \sum\limits_{\nu_2=1}^{
2^{\xi_2} } \nu_1^{r_1+\beta_1} \nu_2^{r_2+\beta_2} A_{\nu_1 \nu_2}
(x_1, x_2) \right\|_{{p}}^\tau
\\
&=:&
K_1+K_2+K_3+K_4.
\end{eqnarray*}
Let us estimate $K_2$.
\begin{eqnarray*}
K_2 &\ll& \sum\limits_{\xi_1=n_1}^{\infty}
\sum\limits_{\xi_2=n_2}^{\infty} 2^{\xi_1 r_1\tau} 2^{-\xi_2
\beta_2\tau}
 \left\|
\sum\limits_{\nu_1=2^{\xi_1}+1}^{ \infty } \sum\limits_{\nu_2=1}^{
2^{n_2} } \nu_2^{r_2+\beta_2} A_{\nu_1 \nu_2} (x_1, x_2)
\right\|_{{p}} ^\tau +
\\
&+& \sum\limits_{\xi_1=n_1}^{\infty}
\sum\limits_{\xi_2=n_2}^{\infty} 2^{\xi_1 r_1\tau} 2^{-\xi_2
\beta_2\tau}
 \left\|
\sum\limits_{\nu_1=2^{\xi_1}+1}^{ \infty }
\sum\limits_{\nu_2=2^{n_2}+1}^{ 2^{\xi_2} } \nu_2^{r_2+\beta_2}
A_{\nu_1 \nu_2} (x_1, x_2) \right\|_{{p}} ^\tau
\\
&=:& K_{21}+K_{22}.
\end{eqnarray*}
Using  (\ref{vsp2}) and Lemmas  \ref{l2.2}, \ref{l2.3}, we get
$$
K_{21} \ll 2^{-n_2 \beta_2 \tau} \sum\limits_{\xi_1=n_1}^{\infty}
2^{\xi_1 r_1\tau} \left(
 \int\limits_{0}^{2\pi}
 \int\limits_{0}^{2\pi}
\left<
\sum\limits_{\nu_1=\xi_1+1}^{\infty}
\sum\limits_{\nu_2=0}^{n_2} 2^{2 \nu_2 (r_2+\beta_2)}
\triangle_{\nu_1 \nu_2}^2 \right>^{p/2} \dd x_1
 \dd x_2 \right)^{\tau/p}.
$$
Using Minkowski's inequality, (\ref{vsp1}), Lemmas \ref{l2.2}, \ref{l2.3},  and Theorem \ref{th4.1}, we have
\begin{eqnarray*}
K_{21} &\ll& 2^{-n_2 \beta_2 \tau} \left\{
 \int\limits_{0}^{2\pi}
 \int\limits_{0}^{2\pi}
\left[ \sum\limits_{\xi_1=n_1}^{\infty} 2^{\xi_1 r_1\tau} \left<
\sum\limits_{\nu_1=\xi_1+1}^{\infty} \sum\limits_{\nu_2=0}^{n_2}
2^{2 \nu_2 (r_2+\beta_2)} \triangle_{\nu_1 \nu_2}^2
\right>^{\tau/2}
\right]^{p/\tau} \dd x_1
 \dd x_2 \right\}^{\tau/p}
\\
&\ll& 2^{-n_2 \beta_2 \tau} \left\{
 \int\limits_{0}^{2\pi}
 \int\limits_{0}^{2\pi}
\left[ \sum\limits_{\nu_1=n_1+1}^{\infty} \left<
\sum\limits_{\xi_1=n_1}^{\nu_1} 2^{\xi_1 r_1\tau}
\right>^{2/\tau} \sum\limits_{\nu_2=0}^{n_2} 2^{2 \nu_2
(r_2+\beta_2)} \triangle_{\nu_1 \nu_2}^2 \right]^{p/2} \dd x_1
 \dd x_2 \right\}^{\tau/p}
\\
&\ll& 2^{-n_2 \beta_2 \tau} \left\{
 \int\limits_{0}^{2\pi}
 \int\limits_{0}^{2\pi}
\left[ \sum\limits_{\nu_1=n_1+1}^{\infty}
\sum\limits_{\nu_2=0}^{n_2} 2^{2 \nu_1 r_1} 2^{2 \nu_2
(r_2+\beta_2)} \triangle_{\nu_1 \nu_2}^2 \right]^{p/2} \dd x_1
 \dd x_2 \right\}^{\tau/p}
\\
&\ll& 2^{-n_2 \beta_2 \tau}
 \left\|
\sum\limits_{\nu_1=2^{n_1}+1}^{ \infty }
\sum\limits_{\nu_2=1}^{2^{n_2} } \nu_1^{r_1} \nu_2^{r_2+\beta_2}
A_{\nu_1 \nu_2} (x_1, x_2) \right\|_{{p}}^\tau
\\
 &\ll &
\omega^\tau_{\beta_1, \beta_2}
(f^{(r_1,r_2)},\delta_1,\delta_2)_{{p}}.
\end{eqnarray*}
Estimating  $K_{22}$ as above, we get:
\begin{eqnarray*}
K_{22} &\ll& \sum\limits_{\xi_2=n_2}^{\infty} 2^{-\xi_2 \beta_2
\tau} \left\{
 \int\limits_{0}^{2\pi}
  \int\limits_{0}^{2\pi}
\left[ \sum\limits_{\nu_1=n_1+1}^{\infty}
\sum\limits_{\nu_2=n_2+1}^{\xi_2} 2^{2\nu_1 r_1} 2^{2 \nu_2
(r_2+\beta_2)} \triangle_{\nu_1 \nu_2}^2 \right]^{p/2} \dd x_1
 \dd x_2 \right\}^{\tau/p}
\\
&\ll& \left\{
 \int\limits_{0}^{2\pi}
 \int\limits_{0}^{2\pi}
\left[ \sum\limits_{\xi_2=n_2}^{\infty} 2^{-\xi_2 \beta_2 \tau}
\left< \sum\limits_{\nu_1=n_1+1}^{\infty}
\sum\limits_{\nu_2=n_2+1}^{\xi_2} 2^{2\nu_1 r_1} 2^{2 \nu_2
(r_2+\beta_2)} \triangle_{\nu_1 \nu_2}^2 \right>^{\tau/2}
\right]^{p/\tau} \dd x_1  \dd x_2
\right\}^{\tau/p}.
\end{eqnarray*}
Further,
\begin{eqnarray*}
K_{22} &\ll&  \left\{
 \int\limits_{0}^{2\pi}
 \int\limits_{0}^{2\pi}
\left[ \sum\limits_{\nu_2=n_2+1}^{\infty} 2^{2 \nu_2 (r_2+\beta_2)}
\left< \sum\limits_{\xi_2=\nu_2}^{\infty} 2^{-\xi_2 \beta_2 \tau}
\right>^{2/\tau} \sum\limits_{\nu_1=n_1+1}^{\infty} 2^{2\nu_1
r_1} \triangle_{\nu_1 \nu_2}^2 \right]^{p/2} \dd x_1
\dd x_2 \right\}^{\tau/p}
\\
&\ll&
 \left\|
\sum\limits_{\nu_1=2^{n_1}+1}^{ \infty }
\sum\limits_{\nu_2=2^{n_2}+1}^{ \infty } \nu_1^{r_1} \nu_2^{r_2}
A_{\nu_1 \nu_2} (x_1, x_2) \right\|_{{p}} ^\tau
\\
& \ll & \omega^\tau_{\beta_1, \beta_2}
(f^{(r_1,r_2)},\delta_1,\delta_2)_{{p}}.
\end{eqnarray*}
Similarly,
\begin{eqnarray*}
K_{3} &\ll& 2^{-n_1 \beta_1 \tau}
 \left\|
\sum\limits_{\nu_1=1}^{2^{n_1} } \sum\limits_{\nu_2=2^{n_2}+1}^{
\infty } \nu_1^{r_1+\beta_1} \nu_2^{r_2} A_{\nu_1 \nu_2} (x_1, x_2)
\right\|_{{p}} ^\tau +
 \left\|
\sum\limits_{\nu_1=2^{n_1}+1}^{ \infty }
\sum\limits_{\nu_2=2^{n_2}+1}^{ \infty } \nu_1^{r_1} \nu_2^{r_2}
A_{\nu_1 \nu_2} (x_1, x_2) \right\|_{{p}} ^\tau
\\
&\ll& \omega^\tau_{\beta_1, \beta_2}
(f^{(r_1,r_2)},\delta_1,\delta_2)_{{p}}.
\end{eqnarray*}
Finally, we obtain the estimates for $K_1$ and $K_4$:
\begin{eqnarray*}
K_{1}& \ll& \sum\limits_{\xi_1=n_1}^{\infty}
\sum\limits_{\xi_2=n_2}^{\infty} 2^{\xi_1 r_1\tau} 2^{\xi_2 r_2\tau}
\left\{
 \int\limits_{0}^{2\pi}
 \int\limits_{0}^{2\pi}
\left[ \sum\limits_{\nu_1=\xi_1+1}^{\infty}
\sum\limits_{\nu_2=\xi_2+1}^{\infty} \triangle_{\nu_1 \nu_2}^2
\right]^{p/2} \dd x_1 \dd x_2
\right\}^{\tau/p}
\\
&\ll&
 \left\|
\sum\limits_{\nu_1=2^{n_1}+1}^{ \infty }
\sum\limits_{\nu_2=2^{n_2}+1}^{ \infty } \nu_1^{r_1} \nu_2^{r_2}
A_{\nu_1 \nu_2} (x_1, x_2) \right\|_{{p}} ^\tau
\ll  \omega^\tau_{\beta_1, \beta_2}
(f^{(r_1,r_2)},\delta_1,\delta_2)_{{p}}
\end{eqnarray*}
and
\begin{eqnarray*}
K_{4} &\ll&
 \left\|
\sum\limits_{\nu_1=2^{n_1}+1}^{ \infty }
\sum\limits_{\nu_2=2^{n_2}+1}^{ \infty } \nu_1^{r_1} \nu_2^{r_2}
A_{\nu_1 \nu_2} (x_1, x_2) \right\|_{{p}} ^\tau + 2^{-n_2
\beta_2 \tau}
 \left\|
\sum\limits_{\nu_1=2^{n_1}+1}^{ \infty }
\sum\limits_{\nu_2=1}^{2^{n_2} } \nu_1^{r_1} \nu_2^{r_2+\beta_2}
A_{\nu_1 \nu_2} (x_1, x_2) \right\|_{{p}} ^\tau
\\
&+& 2^{-n_1 \beta_1 \tau}
 \left\|
\sum\limits_{\nu_1=1}^{2^{n_1} } \sum\limits_{\nu_2=2^{n_2}+1}^{
\infty } \nu_1^{r_1+\beta_1} \nu_2^{r_2} A_{\nu_1 \nu_2} (x_1, x_2)
\right\|_{{p}} ^\tau\\ &+& 2^{-n_1 \beta_1 \tau -n_2 \beta_2 \tau}
 \left\|
\sum\limits_{\nu_1=1}^{2^{n_1} } \sum\limits_{\nu_2=1}^{2^{n_2} }
\nu_1^{r_1+\beta_1} \nu_2^{r_2+\beta_2} A_{\nu_1 \nu_2} (x_1, x_2)
\right\|_{{p}}^\tau
\ll \omega^\tau_{\beta_1, \beta_2}
(f^{(r_1,r_2)},\delta_1,\delta_2)_{{p}}.
\end{eqnarray*}
\hfill $\Box$
\\[0.5cm]
{\bf Proof of Theorems \ref{T11} and \ref{T12}.}
These proofs are similar to the proof of Theorem \ref{T10}. We use the following two statements which follow from Lemmas \ref{l2.6}, \ref{l2.7} and Theorems \ref{th7.2}, \ref{th7.3}.

 Let $f \in M_p, 1 < p < \infty, \beta_i, r_i > 0, i=1,2$, $n_1, n_2 \in \N$.
 Then
\begin{eqnarray*}
\omega_{\beta_1,\beta_2} \Big(f^{(r_1,r_2)}, \frac{1}{n_1}, \frac{1}{n_2} \Big)_{L_p(\T^2)}
&\asymp&
\frac{1}{n_1^{\beta_1}n_2^{\beta_2}} \Bigg\{ \sum\limits_{\nu_1=1}^{n_1} \sum\limits_{\nu_2=1}^{n_2} a_{\nu_1,\nu_2}^p
\nu_1^{(\beta_1+r_1+1)p -2}\nu_2^{(\beta_2+r_2+1)p -2} \Bigg\}^{{1}/{p}}
\\
&+&
\frac{1}{n_1^{\beta_1}} \Bigg\{ \sum\limits_{\nu_1=1}^{n_1} \sum\limits_{\nu_2=n_2+1}^{\infty} a_{\nu_1,\nu_2}^p
\nu_1^{(\beta_1+r_1
+1)p -2}\nu_2^{(r_2+1)p -2} \Bigg\}^{{1}/{p}}
\\&+&
\frac{1}{n_2^{\beta_2}} \Bigg\{ \sum\limits_{\nu_1=n_1+1}^{\infty} \sum\limits_{\nu_2=1}^{n_2} a_{\nu_1,\nu_2}^p \nu_1^{(r_1+1)p -2}\nu_2^{(\beta_2+r_2+1)p -2} \Bigg\}^{{1}/{p}}
\\
 &+&\Bigg\{ \sum\limits_{\nu_1=n_1+1}^{\infty} \sum\limits_{\nu_2=n_2+1}^{\infty} a_{\nu_1,\nu_2}^p \nu_1^{(r_1+1)p -2}\nu_2^{(r_2+1)p -2} \Bigg\}^{{1}/{p}}.
\end{eqnarray*}

Let $f \in \Lambda_p, 1 < p < \infty, \beta_i, r_i > 0,  n_i\in \N,   i=1,2$.
 Then
\begin{eqnarray*}
 \omega_{\beta_1,\beta_2} \Big(f^{(r_1,r_2)}, \frac{1}{2^{n_1}}, \frac{1}{2^{n_2}} \Big)_p &\asymp&
\frac{1}{2^{n_1\beta_1+n_2 \beta_2}} \Bigg\{
\sum\limits_{\mu_1=0}^{n_1}\sum\limits_{\mu_2=0}^{n_2}
\lambda_{\mu_1,\mu_2}^2 2^{2(\mu_1(\beta_1+r_1)+\mu_2 (\beta_2+r_2))}
\Bigg\}^{{1}/{2}}
\\
&+&
\frac{1}{2^{n_1\beta_1}} \Bigg\{
\sum\limits_{\mu_1=0}^{n_1}\sum\limits_{\mu_2=n_2+1}^{\infty}
\lambda_{\mu_1,\mu_2}^2 2^{2(\mu_1(\beta_1+r_1) + \mu_2r_2)}\Bigg\}^{{1}/{2}}
\\
&+&\frac{1}{2^{n_2\beta_2}} \Bigg\{
\sum\limits_{\mu_1=n_1+1}^{\infty}\sum\limits_{\mu_2=0}^{n_2}
\lambda_{\mu_1,\mu_2}^2
2^{2(\mu_1r_1
+ \mu_2(\beta_2+r_2))} \Bigg\}^{{1}/{2}}
\\
&+&
\Bigg\{ \sum\limits_{\mu_1=n_1+1}^{\infty}\sum\limits_{\mu_2=n_2+1}^{\infty} \lambda_{\mu_1,\mu_2}^2
2^{2\mu_1r_1+ 2\mu_2r_2} \Bigg\}^{{1}/{2}}.
\end{eqnarray*}

\hfill $\Box$

Note that Theorems \ref{T10}-\ref{T12} deal with the most important case when
 $\omega_{\bf{\beta}}(f^{(\bf{r})}, {\bf{\delta}})_p$ is estimated
in terms of   $\omega_{\bf{r}+\bf{\beta}}(f, {\bf{t}})_p$.
However, to completely solve the general problem on the interrelation between
 $$\omega_{\bf{\beta}}(f^{(\bf{l})}, {\bf{\delta}})_p\qquad\mbox{
and}\qquad \omega_{\bf{r}+\bf{\beta}}(f, {\bf{t}})_p,$$
 we have to consider two more cases:

\begin{itemize}
  \item[(i).] \qquad $\bf{l}=\bf{r}-\bf{\alpha}, \quad\bf{0}<\bf{\alpha}<\bf{r}$;
  \item[(ii).] \qquad $\bf{l}=\bf{r}+\bf{\alpha},\quad \bf{0}<\bf{\alpha}<\bf{\beta}$.
\end{itemize}
The results covering these cases are provided below:

\begin{theorem}
{\textnormal{(i).}}
Let $f \in L_p^0(\T^2),$ $1<p<\infty$, $\theta:=\min(2,p)$, and $\tau:=\max(2,p)$.
Let $\beta_1, \beta_2, >0$, $0<\alpha_1< r_1,$ $0<\alpha_2< r_2$,
 and $\delta_1,\delta_2\in (0,\frac12).$
\\
{\textnormal {A).}} If
$$
\int\limits^{1}_{0}
\int\limits^{1}_{0}
t_1^{-(r_1-\alpha_1)\theta-1}
t_2^{-(r_2-\alpha_2)\theta-1}
\omega^\theta_{r_1+\beta_1, r_2+\beta_2}(f,t_1, t_2)_{L_p(\T^2)}
\,\dd t_1\,  \dd t_2\,<\infty,
$$
then  $f$ has a mixed derivative
$f^{(r_1-\alpha_1,r_2-\alpha_2)} \in {L_p^0(\T^2)}$ in the sense of Weyl and
$$\omega_{\beta_1, \beta_2} (f^{(r_1-\alpha_1,r_2-\alpha_2)},\delta_1, \delta_2)_{L_p(\T^2)}
\ll \, J(\theta),
$$
where
$$
J(s):=
\left\{
\int\limits_{0}^{1}\int\limits_{0}^{1}
t_1^{-(r_1-\alpha_1)s}
\min\Big(1,\frac{\delta_1}{t_1}\Big)^{\beta_1s}
 t_2^{-(r_2-\alpha_2)s}
\min\Big(1,\frac{\delta_2}{t_2}\Big)^{\beta_2s}
 \omega^s_{r_1+\beta_1, r_2+\beta_2}(f,t_1, t_2)_{p} \,\frac{\dd t_1}{t_1}\,
\frac{\dd t_2}{t_2}
\right\}^{{1}/{s}}.
$$
{\textnormal {B).}} If  $f\in {L_p(\T^2)}$ has a mixed derivative
$f^{(r_1-\alpha_1,r_2-\alpha_2)} \in {L_p^0(\T^2)}$  in the sense of Weyl,
then
$$\omega_{\beta_1, \beta_2} (f^{(r_1-\alpha_1,r_2-\alpha_2)},\delta_1, \delta_2)_{L_p(\T^2)}
\gg \, J(\tau).
$$
{\textnormal{(ii).}}
Let $f \in L_p^0(\T^2),$ $1<p<\infty$, $\theta:=\min(2,p)$, and $\tau:=\max(2,p)$.
Let $r_1, r_2>0$, $0<\alpha_1< \beta_1,$ $0<\alpha_2< \beta_2$, and $\delta_1,\delta_2\in (0,\frac12).$
\\
{\textnormal {A).}} If
\begin{equation}\label{wwwww}
 \int\limits^{1}_{0}\int\limits^{1}_{0}
t_1^{-(r_1+\alpha_1)\theta-1} t_2^{-(r_2+\alpha_2)\theta-1}
\omega^\theta
_{r_1+\beta_1, r_2+\beta_2}(f,t_1, t_2)_{L_p(\T^2)} \,\dd t_1\,  \dd t_2\,<\infty,
\end{equation}
then  $f$ has a mixed derivative
$f^{(r_1+\alpha_1,r_2+\alpha_2)} \in {L_p^0(\T^2)}$ in the sense of Weyl and
$$
D(f^{(r_1+\alpha_1,r_2+\alpha_2)}, \tau)\ll
 \, E(f, \theta),
$$
where
$$D(f^{(r_1+\alpha_1,r_2+\alpha_2)}, s)
:=
\left\{
\int\limits_{\delta_1}^{1} \int\limits_{\delta_2}^{1}
\Big(\frac{\delta_1}{t_1}\Big)^{(\beta_1-\alpha_1){s}}
\Big(\frac{\delta_2}{t_2}\Big)^{(\beta_2-\alpha_2){s}}
\omega^s_{\beta_1, \beta_2} (f^{(r_1+\alpha_1,r_2+\alpha_2)},t_1,t_2)_{L_p(\T^2)}\,
\frac{\dd t_1}{t_1}\,\frac{\dd t_2}{t_2} \right\}^{{1}/{s}}
$$
and
$$ \, E(f, s)
:=
 \left\{
\int\limits_{0}^{\delta_1}\int\limits_{0}^{\delta_2}
t_1^{-(r_1+\alpha_1){s}-1}
t_2^{-(r_2+\alpha_2){s}-1}
 \omega^{\theta}_{r_1+\beta_1, r_2+\beta_2}(f,t_1,t_2)_{L_p(\T^2)}\, \dd t_1\,\dd t_2
\right\}^{{1}/{s}}.
$$
{\textnormal {B).}} If  $f\in {L_p(\T^2)}$ has a mixed derivative
$f^{(r_1+\alpha_1,r_2+\alpha_2)} \in {L_p^0(\T^2)}$  in the sense of Weyl,
then
$$
E(f, \tau)
\ll
 \, D(f^{(r_1+\alpha_1,r_2+\alpha_2)}, \theta).
$$
\end{theorem}
The proof can be found in \cite{deriv-st}; see also \cite{deriv-stt}.
As in Theorems \ref{T11} and \ref{T12}, the sharpness of this result can be shown by considering function classes $M_p$ and $\Lambda_p$ defined in Section 2.4.

Moreover, in part (ii),
$D(f^{(r_1+\alpha_1,r_2+\alpha_2)},\theta)$ could not be replaced by
$\omega_{\beta_1,\beta_2}(f^{(r_1+\alpha_1,r_2+\alpha_2)},
\delta_1, \delta_2)_{L_p(\T^2)}$.
Indeed, we consider
$$f_1(x,y) :=  \sum\limits_{n_1=1}^{\infty}
\sum\limits_{n_2=1}^{\infty}
\frac{n_1^{-(r_1+\beta_1+1-\frac{1}{p})}}{\ln^{A_1+\frac{1}{p}}(2n_1)} \,
\frac{n_2^{-(r_2+\beta_2+1-\frac{1}{p})}}{\ln^{A_2+\frac{1}{p}}(2n_2)} \,
\cos n_1x\,
\cos n_2y,$$
where $p\in(1,\infty),$  $r_i,\beta_i, A_i>0\,(i=1,2)$. Then, by Lemma \ref{l2.6},
we have $f\in {L_p^0(\T^2)}$ and, by Theorem  \ref{th7.2},
$$\omega_{r_1+\beta_1, r_2+\beta_2}\left(f_1,  \delta_1, \delta_2 \right)_{L_p(\T^2)}
\quad\asymp \quad
 \delta_1^{r_1+\beta_1} \delta_2^{r_2+\beta_2}
.$$
Hence condition (\ref{wwwww}) holds and
$f_1^{(r_1+\alpha_1,r_2+\alpha_2)}\in {L_p^0(\T^2)}$. Moreover,
one can easily check that $$D(f_1^{(r_1+\alpha_1,r_2+\alpha_2)},p) \asymp E(f_1,p)\asymp
\delta_1^{\beta_1-\alpha_1}\delta_2^{\beta_2-\alpha_2}.$$
On the other hand,
$$\omega_{\beta_1,\beta_2}(f_1^{(r_1+\alpha_1,r_2+\alpha_2)}, \delta_1, \delta_2 )_{L_p(\T^2)}
\quad\asymp \quad
\delta_1^{\beta_1-\alpha_1}\delta_2^{\beta_2-\alpha_2}
\big|\ln \delta_1\big|^{-A_1}
\big|\ln \delta_2\big|^{-A_2}
$$
and thus
$D(f_1^{(r_1+\alpha_1,r_2+\alpha_2)},p)$ and $\omega_{\beta_1,\beta_2}(f_1^{(r_1+\alpha_1,r_2+\alpha_2)})_{L_p(\T^2)}$
indeed have different orders of magnitude.

\vskip 1.5cm

\sect{The mixed moduli of smoothness of $L_p$-functions and their angular approximation}

As we already mentioned, Jackson and Bernstein-Stechkin type results are given by the following inequalities:
\begin{eqnarray*}
Y_{n_1,n_2}(f)_{L_p(\T^2)}&\ll& \omega_{k_1,k_2}(f,\frac{\pi}{n_1+1},\frac{\pi}{n_2+1})_{L_p(\T^2)}
\\
&\ll& \frac{1}{(n_1+1)^{k_1}(n_2+1)^{k_2}}\sum\limits^{n_1+1}_{\nu_1=1}
\sum\limits^{n_2+1}_{\nu_2=1}\nu_1^{k_1-1}\nu_2^{k_2-1}Y_{\nu_1-1,\nu_2-1}(f)_{L_p(\T^2)}.
\end{eqnarray*}
 The next theorem provides sharper estimates (sharp Jackson and sharp inverse inequality).

\begin{theorem}\label{th8.2}
Let $f\in L_p^0(\T^2), 1<p<\infty$,  $\sigma:=\max(2,p), \theta:=\min(2,p)$,   $\alpha_i>0$,  $n_i\in \mathbb{N},$ $i=1,2.$  Then
\begin{multline*}
\frac{1}{n_1^{\alpha_1}}\frac{1}{n_2^{\alpha_2}}  \left\{
\sum\limits^{n_1+1}_{\nu_1=1}\sum\limits^{n_2+1}_{\nu_2=1}\nu_1^{\alpha_1\sigma-1}\nu_2^{\alpha_2\sigma-1}Y^{\sigma}_{\nu_1-1,\nu_2-1}(f)_{L_p(\T^2)}
\right\} ^{1/{\sigma}} \ll
\omega_{\alpha_1,\alpha_2}\left(f,\frac1{n_1},\frac1{n_2} \right)_{L_p(\T^2)}\\
\ll  \frac{1}{n_1^{\alpha_1}}\frac{1}{n_2^{\alpha_2}}  \left\{
\sum\limits^{n_1+1}_{\nu_1=1}\sum\limits^{n_2+1}_{\nu_2=1}\nu_1^{\alpha_1\theta-1}\nu_2^{\alpha_2\theta-1}Y^{\theta}_{\nu_1-1,\nu_2-1}(f)_{L_p(\T^2)}
\right\} ^{1/{\theta}}.
\end{multline*}
\end{theorem}
The sharpness of this result follows from considering  special function classes, see
Section 2.4.
\begin{theorem}\label{th8.3}
Let $f\in M_p, 1<p<\infty, \alpha_i>0, n_i\in \mathbb{N}, i=1,2.$ Then
$$
\omega_{\alpha_1,\alpha_2}\left(f,\frac1{n_1},\frac1{n_2}
\right)_{L_p(\T^2)} \asymp
\frac{1}{n_1^{\alpha_1}}\frac{1}{n_2^{\alpha_2}}  \left\{
\sum\limits^{n_1+1}_{\nu_1=1}\sum\limits^{n_2+1}_{\nu_2=1}
\nu_1^{\alpha_1 p-1}\nu_2^{\alpha_2
p-1}Y^{p}_{\nu_1-1,\nu_2-1}(f)_{L_p(\T^2)}   \right\} ^{1/{p}} .
$$
\end{theorem}

\begin{theorem}\label{th8.4}
Let $f\in \Lambda_p, 1<p<\infty, \alpha_i>0, n_i\in \mathbb{N}, i=1,2.$  Then
$$
\omega_{\alpha_1,\alpha_2}\left(f,\frac1{n_1},\frac1{n_2}
\right)_{L_p(\T^2)} \asymp
\frac{1}{n_1^{\alpha_1}}\frac{1}{n_2^{\alpha_2}}  \left\{
\sum\limits^{n_1+1}_{\nu_1=1}\sum\limits^{n_2+1}_{\nu_2=1}\nu_1^{2\alpha_1
-1} \nu_2^{2\alpha_2 -1}Y^{2}_{\nu_1-1,\nu_2-1}(f)_{L_p(\T^2)}
\right\}^{1/{2}} .
$$

\end{theorem}
{\bf Proof of Theorem \ref{th8.2}.}
Define
$$
I^{\sigma}:=\frac{1}{n^{\alpha_1 \sigma}n^{\alpha_2 \sigma}}\sum\limits^{n_1+1}_{\nu_1=1}\sum\limits^{n_2+1}_{\nu_2=1}\nu_1^{\alpha_1\sigma-1}
\nu_2^{\alpha_2\sigma-1}Y^{\sigma}_{\nu_1-1,\nu_2-1}(f)_{L_p(\T^2)}.
$$
For given  $n_i \in \mathbb{N}$ we find integers  $m_i \geq 0$ such that $2^{m_i}\le n_i< 2^{m_i+1}, i=1,2.$
Then using monotonicity properties of  $Y_{\nu_1,\nu_2}(f)_{L_p(\T^2)}$ we get
$$
I^{\sigma}\ll
2^{-(\alpha_1m_1+\alpha_2m_2)\sigma}\sum\limits^{m_1+1}_{\mu_1=0}\sum\limits^{m_2+1}_{\mu_2=0}2^{(\mu_1\alpha_1+\mu_2\alpha_2)\sigma}Y^{\sigma}_
{\floor{2^{\mu_1-1}},\floor{2^{\mu_2-1}}}(f)_{L_p(\T^2)}.
$$
Using Lemmas \ref{l2.2} and \ref{l2.3}, we have
\begin{eqnarray*}
I^{\sigma}&\ll &
2^{-(\alpha_1m_1+\alpha_2m_2)\sigma}
 \sum\limits^{m_1+1}_{\mu_1=0}\sum\limits^{m_2+1}_{\mu_2=0}2^{(\mu_1\alpha_1+\mu_2\alpha_2)\sigma}
 \left\{\int\limits^{2\pi}_{0}\int\limits^{2\pi}_{0} \left(\sum\limits^{\infty}_{\nu_1=\mu_1}
 \sum\limits^{\infty}_{\nu_2=\mu_2} \Delta^2_{\nu_1\nu_2}\right)^{{p}/{2}}\dd x\dd y  \right\}^{{\sigma}/{p}}
 \\
 &\ll&
  \left\{\int\limits^{2\pi}_{0}\int\limits^{2\pi}_{0} \left(\sum\limits^{\infty}_{\nu_1=m_1+1}
 \sum\limits^{\infty}_{\nu_2=m_2+1} \Delta^2_{\nu_1\nu_2}\right)^{{p}/{2}}\dd x\dd y  \right\}^{{\sigma}/{p}}
 \\
 &+&
 2^{-\alpha_1m_1\sigma}\sum\limits^{m_1}_{\mu_1=0}2^{\alpha_1\mu_1\sigma}
 \left\{\int\limits^{2\pi}_{0}\int\limits^{2\pi}_{0} \left(\sum\limits^{m_1}_{\nu_1=\mu_1}
 \sum\limits^{\infty}_{\nu_2=m_2+1} \Delta^2_{\nu_1\nu_2}\right)^{{p}/{2}}\dd x\dd y  \right\}^{{\sigma}/{p}}
\\
  &+&
 2^{-\alpha_2m_2\sigma}\sum\limits^{m_2}_{\mu_2=0}2^{\alpha_2\mu_2\sigma}
 \left\{\int\limits^{2\pi}_{0}\int\limits^{2\pi}_{0} \left(\sum\limits^{\infty}_{\nu_1=m_1+1}
 \sum\limits^{m_2}_{\nu_2=\mu_2} \Delta^2_{\nu_1\nu_2}\right)^{{p}/{2}}\dd x\dd y  \right\}^{{\sigma}/{p}}
 \\
 &+&
 2^{(-\alpha_1m_1-\alpha_2m_2)\sigma}
 \sum\limits^{m_1+1}_{\mu_1=1}\sum\limits^{m_2+1}_{\mu_2=1}2^{(\mu_1\alpha_1+\mu_2\alpha_2)\sigma}
 \left\{\int\limits^{2\pi}_{0}\int\limits^{2\pi}_{0} \left(\sum\limits^{m_1}_{\nu_1=\mu_1}
 \sum\limits^{m_2}_{\nu_2=\mu_2} \Delta^2_{\nu_1\nu_2}\right)^{{p}/{2}}\dd x\dd y  \right\}^{{\sigma}/{p}}
 \\&=:&I_1+I_2+I_3+I_4.
\end{eqnarray*}
By Lemma \ref{l2.3}, we have
$$
I_1\ll \|f-s_{2^{m_1,\infty}}- s_{\infty,2^{m_2}}+s_{2^{m_1},2^{m_2}} \|^{\sigma}_{L_p(\T^2)}.
$$
Now we estimate
$$
 I_2=2^{-\alpha_1m_1\sigma}\sum\limits^{m_1}_{\mu_1=0}2^{\alpha_1\mu_1\sigma}
 \left\{\int\limits^{2\pi}_{0}\int\limits^{2\pi}_{0} \left(\sum\limits^{m_1}_{\nu_1=\mu_1}
 \sum\limits^{\infty}_{\nu_2=m_2+1} \Delta^2_{\nu_1\nu_2}\right)^{{p}/{2}}\dd x\dd y  \right\}^{{\sigma}/{p}}.
$$
Minkowski's inequality  $\left({\sigma}/{p}\geq 1 \right)$ implies
$$
 I_2\ll  2^{-\alpha_1m_1\sigma}
 \left\{\int\limits^{2\pi}_{0}\int\limits^{2\pi}_{0} \left[\sum\limits^{m_1}_{\mu_1=0}2^{\alpha_1\mu_1\sigma} \left(\sum\limits^{m_1}_{\nu_1=\mu_1}
 \sum\limits^{\infty}_{\nu_2=m_2+1} \Delta^2_{\nu_1\nu_2}\right)^{{\sigma}/{2}}\right]^{{p}/{\sigma}}\dd x\dd y  \right\}^{{\sigma}/{p}}.
$$
Taking into account Lemma \ref{l2.8}, we get
$$
 I_2\ll  2^{-\alpha_1m_1\sigma}
 \left\{\int\limits^{2\pi}_{0}\int\limits^{2\pi}_{0} \left[\sum\limits^{m_1}_{\nu_1=0}2^{\alpha_1\nu_1\sigma} \left(
 \sum\limits^{\infty}_{\nu_2=m_2+1} \Delta^2_{\nu_1\nu_2}\right)^{{\sigma}/{2}}\right]^{{p}/{\sigma}}\dd x\dd y  \right\}^{{\sigma}/{p}}.
$$
Since ${\sigma}/{2}\ge 1$, Lemma \ref{l2.9} yields
$$
 I_2\ll  2^{-\alpha_1m_1\sigma}
 \left\{\int\limits^{2\pi}_{0}\int\limits^{2\pi}_{0} \left[\sum\limits^{m_1}_{\nu_1=0} \sum\limits^{\infty}_{\nu_2=m_2+1}
 2^{\alpha_1\nu_12}\Delta^2_{\nu_1\nu_2}\right]^{{p}/{2}}\dd x\dd y  \right\}^{{\sigma}/{p}}.
$$
Then Lemmas  \ref{l2.2} and \ref{l2.3} imply
$$
 I_2\ll  2^{-\alpha_1m_1\sigma}
\|s^{(\alpha_1,0)}_{2^{m_1},\infty}(f-s_{\infty,2^{m_2}}(f))\|^\sigma_{L_p(\T^2)}.
$$
Similarly,
$$
 I_3\ll  2^{-\alpha_2m_2\sigma}
\|s^{(0,\alpha_2)}_{\infty,2^{m_2}}(f-s_{2^{m_1},\infty}(f))\|^\sigma_{L_p(\T^2)}
$$
and
$$
 I_4\ll  2^{-(\alpha_1m_1+\alpha_2m_2)\sigma}
\|s^{(\alpha_1,\alpha_2)}_{2^{m_1},2^{m_2}}(f)\|^\sigma_{L_p(\T^2)}.
$$
Therefore,
\begin{eqnarray*}
I^\sigma
&\ll &
\|f-s_{2^{m_1,\infty}}- s_{\infty,2^{m_2}}+s_{2^{m_1},2^{m_2}} \|^{\sigma}_{L_p(\T^2)}
\\
&+&
2^{-\alpha_1m_1}\|s^{(\alpha_1,0)}_{2^{m_1},\infty}(f-s_{\infty,2^{m_2}}(f))\|^\sigma_{L_p(\T^2)}+
2^{-\alpha_2m_2}\|s^{(0,\alpha_2)}_{\infty,2^{m_2}}(f-s_{2^{m_1},\infty}(f))\|^\sigma_{L_p(\T^2)}
\\
&+&
 2^{-(\alpha_1m_1+\alpha_2m_2)\sigma}
\|s^{(\alpha_1,\alpha_2)}_{2^{m_1},2^{m_2}}(f)\|^\sigma_{L_p(\T^2)}.
\end{eqnarray*}
Applying Theorem \ref{th4.1}, we get
 $I\ll
\omega_{\alpha_1,\alpha_2}\left(f,\frac{1}{2^{m_1}},\frac{1}{2^{m_2}}
\right)_{L_p(\T^2)}$. Then using properties of the mixed modulus of smoothness,
$$
I\ll  \omega_{\alpha_1,\alpha_2}\left(f,\frac{1}{2^{m_1+1}},\frac{1}{2^{m_2+1}} \right)_{L_p(\T^2)}\ll
\omega_{\alpha_1,\alpha_2}\left(f,\frac{1}{n_1},\frac{1}{n_2} \right)_{L_p(\T^2)}.
$$
Thus, the estimate from below in Theorem \ref{th8.2} is shown. Now we prove the estimate from above.
First,
$$
I_5:=\omega_{\alpha_1,\alpha_2}\left(f,\frac{1}{n_1},\frac{1}{n_2} \right)_{L_p(\T^2)}\ll
 \omega_{\alpha_1,\alpha_2}\left(f,\frac{1}{2^{m_1}},\frac{1}{2^{m_2}} \right)_{L_p(\T^2)},
$$
where $2^{m_i}\le n_i <2^{m_i+1}, i=1,2$.
Theorem  \ref{th4.1} gives
\begin{eqnarray*}
I^p_5&\ll&   2^{-(\alpha_1m_1+\alpha_2m_2)p}
\|s^{(\alpha_1,\alpha_2)}_{2^{m_1},2^{m_2}}(f)\|^p_{L_p(\T^2)}+
2^{-\alpha_1 m_1 p}
\|s^{(\alpha_1,0)}_{2^{m_1},\infty}(f-s_{\infty,2^{m_2}}(f))\|^p_{L_p(\T^2)}
\\
&+&
2^{-\alpha_2m_2p}\|s^{(0,\alpha_2)}_{\infty,2^{m_2}}(f-s_{2^{m_1},\infty}(f))\|^p_{L_p(\T^2)}
+
\|f-s_{2^{m_1,\infty}}- s_{\infty,2^{m_2}}+s_{2^{m_1},2^{m_2}}
\|^{p}_{L_p(\T^2)}
\\
&=:&I_6+I_7+I_8+I_9.
\end{eqnarray*}
Lemmas \ref{l2.2} and \ref{l2.3} imply
$$
I_6\ll 2^{-(m_1\alpha_1+m_2\alpha_2)p}
\int\limits^{2\pi}_{0}\int\limits^{2\pi}_{0} \left(
\sum\limits^{m_1}_{\mu_1=0}\sum\limits^{m_2}_{\mu_2=0}2^{2(\mu_1\alpha_1+\mu_2\alpha_2)}\Delta^2_{\mu_1\mu_2}
\right)^{{p}/{2}}\dd x\dd y.
$$
If $2 \le p < \infty$, we use Minkowski's inequality and Lemma \ref{l2.1}
\begin{eqnarray*}
I_6&\ll & 2^{-(m_1\alpha_1+m_2\alpha_2)p} \left\{
\sum\limits^{m_1}_{\mu_1=0}\sum\limits^{m_2}_{\mu_2=0}2^{2(\mu_1\alpha_1+\mu_2\alpha_2)}
\left( \int\limits^{2\pi}_{0}\int\limits^{2\pi}_{0}
|\Delta_{\mu_1\mu_2}|^p\dd x\dd y \right)^{{2}/{p}}
\right\}^{{p}/{2}}
\\
&=& 2^{(m_1\alpha_1+m_2\alpha_2)p}\left\{
\sum\limits^{m_1}_{\mu_1=0}\sum\limits^{m_2}_{\mu_2=0}2^{2(\mu_1\alpha_1+\mu_2\alpha_2)}
\|s_{2^{\mu_1},2^{\mu_2}}(f)-s_{2^{\mu_1},\floor{2^{\mu_2-1}}}(f)-
\right.
\\
&-&\left.s_{\floor{2^{\mu_1-1}},2^{\mu_2}}(f)+
s_{\floor{2^{\mu_1-1}},\floor{2^{\mu_2-1}}}(f)\|^2_{L_p(\T^2)} \right\}^{{p}/{2}}
\\& \ll &
2^{-(m_1\alpha_1+m_2\alpha_2)p} \left\{
\sum\limits^{m_1+1}_{\mu_1=0}\sum\limits^{m_2+1}_{\mu_2=0}2^{2(\mu_1\alpha_1+\mu_2\alpha_2)}
Y^{2}_{\floor{2^{\mu_1-1}},\floor{2^{\mu_2-1}}}(f)_{L_p(\T^2)}
\right\}^{{p}/{2}}.
\end{eqnarray*}
If  $1<p<2$, we use Lemma \ref{l2.9} and  Lemma \ref{l2.1}
\begin{eqnarray*}
I_6&\ll & 2^{-(m_1\alpha_1+m_2\alpha_2)p}
\sum\limits^{m_1}_{\mu_1=0}\sum\limits^{m_2}_{\mu_2=0}2^{p(\mu_1\alpha_1+\mu_2\alpha_2)}
 \int\limits^{2\pi}_{0}\int\limits^{2\pi}_{0} |\Delta_{\mu_1,\mu_2}|^p\dd x\dd y
\\
&=& 2^{(m_1\alpha_1+m_2\alpha_2)p}
\sum\limits^{m_1}_{\mu_1=0}\sum\limits^{m_2}_{\mu_2=0}2^{p(\mu_1\alpha_1+\mu_2\alpha_2)}
\|s_{2^{\mu_1},2^{\mu_2}}(f)-s_{2^{\mu_1},\floor{2^{\mu_2-1}}}(f)-
\\&-&
 s_{\floor{2^{\mu_1-1}},2^{\mu_2}}(f)+s_{\floor{2^{\mu_1-1}},\floor{2^{\mu_2-1}}}(f)\|^p_{L_p(\T^2)}
\\
&\ll & 2^{-(m_1\alpha_1+m_2\alpha_2)p}
\sum\limits^{m_1+1}_{\mu_1=0}\sum\limits^{m_2+1}_{\mu_2=0}2^{p(\mu_1\alpha_1+\mu_2\alpha_2)}
Y^{p}_{\floor{2^{\mu_1-1}},\floor{2^{\mu_2-1}}}(f)_{L_p(\T^2)} .
\end{eqnarray*}
Thus, we show
$$
I_6\ll 2^{-(m_1\alpha_1+m_2\alpha_2)p} \left\{
\sum\limits^{m_1+1}_{\mu_1=0}\sum\limits^{m_2+1}_{\mu_2=0}2^{\theta(\mu_1\alpha_1+\mu_2\alpha_2)}
Y^{\theta}_{\floor{2^{\mu_1-1}},\floor{2^{\mu_2-1}}}(f)_{L_p(\T^2)}
\right\}^{{p}/{\theta}}.
$$
Now let us estimate  $I_7$. Lemmas  \ref{l2.2} and \ref{l2.3} imply
$$
I_7\ll 2^{-m_1\alpha_1p}
\int\limits^{2\pi}_{0}\int\limits^{2\pi}_{0} \left(
\sum\limits^{m_1}_{\mu_1=0}\sum\limits^{\infty}_{\mu_2=m_2+1}
2^{2\mu_1\alpha_1} \Delta_{\mu_1,\mu_2}^2
 \right)^{{p}/{2}}\dd x\dd y.
$$
Again, for  $2\le p <\infty$, Minkowski's inequality, Lemmas  \ref{l2.3} and \ref{l2.1} give
\begin{eqnarray*}
&I_7\ll & 2^{m_1\alpha_1p} \left\{
\sum\limits^{m_1}_{\mu_1=0}2^{2\mu_1 \alpha_1} \left[
\int\limits^{2\pi}_{0} \int\limits^{2\pi}_{0} \left(
\sum\limits^{\infty}_{\mu_2=m_2+1}\Delta^2_{\mu_1\mu_2}
\right)^{{p}/{2}}\dd x\dd y \right]^{{2}/{p}}
\right\}^{{p}/{2}}
\\&\ll &
2^{m_1\alpha_1p} \left\{ \sum\limits^{m_1}_{\mu_1=0}2^{2\mu_1
\alpha_1} \left[ \int\limits^{2\pi}_{0} \int\limits^{2\pi}_{0}
\left(
\sum\limits^{\infty}_{\nu_1=\mu_1}\sum\limits^{\infty}_{\mu_2=m_2+1}\Delta^2_{\nu_1\mu_2}
\right)^{{p}/{2}}\dd x\dd y \right]^{{2}/{p}}
\right\}^{{p}/{2}}
\\ &\ll &
2^{m_1\alpha_1p} \left\{ \sum\limits^{m_1+1}_{\mu_1=0}2^{2\mu_1 \alpha_1} Y^2_{\floor{2^{\mu_1-1}},2^{m_2}}(f)_{L_p(\T^2)} \right\}^{{p}/{2}}
\\ &\ll &
2^{-(m_1\alpha_1+m_2\alpha_2)p} \left\{
\sum\limits^{m_1+1}_{\mu_1=0}\sum\limits^{m_2+1}_{\mu_2=0}2^{2(\mu_1
\alpha_1+\mu_2\alpha_2)}
Y^2_{\floor{2^{\mu_1-1}},\floor{2^{\mu_2-1}}}(f)_{L_p(\T^2)}
\right\}^{{p}/{2}}.
\end{eqnarray*}
If  $1<p\le 2$, Lemmas  \ref{l2.9} and  \ref{l2.1} imply

\begin{eqnarray*}
 I_7 &\ll&  2^{-m_1\alpha_1p}
\int\limits^{2\pi}_{0}\int\limits^{2\pi}_{0}
\sum\limits^{m_1}_{\mu_1=0}\sum\limits^{\infty}_{\mu_2=m_2+1}
2^{\mu_1\alpha_1 p} |\Delta_{\mu_1\mu_2}|^p
 \dd x\dd y  \\ &\ll&  2^{-m_1\alpha_1p}
\sum\limits^{m_1}_{\mu_1=0}2^{p\mu_1 \alpha_1}
\int\limits^{2\pi}_{0} \int\limits^{2\pi}_{0}
\sum\limits^{\infty}_{\nu_1=\mu_1}\sum\limits^{\infty}_{\mu_2=m_2+1}
|\Delta_{\nu_1\mu_2}|^p \dd x\dd y
\\ &\ll&
2^{-m_1\alpha_1p}  \sum\limits^{m_1+1}_{\mu_1=0}2^{p\mu_1
\alpha_1} Y^p_{\floor{2^{\mu_1-1}},2^{m_2}}(f)_{L_p(\T^2)}
\\ &\ll&
2^{-(m_1\alpha_1+m_2\alpha_2)p}  \sum\limits^{m_1+1}_{\mu_1=0}\sum\limits^{m_2+1}_{\mu_2=0}2^{p(\mu_1 \alpha_1+\mu_2\alpha_2)} Y^p_{\floor{2^{\mu_1-1}},\floor{2^{\mu_2-1}}}(f)_{L_p(\T^2)} .
\end{eqnarray*}
Thus,
$$
I_7\ll 2^{-(m_1\alpha_1+m_2\alpha_2)p} \left\{ \sum\limits^{m_1+1}_{\mu_1=0}\sum\limits^{m_2+1}_{\mu_2=0}2^{\theta(\mu_1 \alpha_1+\mu_2\alpha_2)}
Y^\theta_{\floor{2^{\mu_1-1}},\floor{2^{\mu_2-1}}}(f)_{L_p(\T^2)} \right\}^{{p}/{\theta}}.
$$
Similarly,
$$
I_8\ll 2^{-(m_1\alpha_1+m_2\alpha_2)p} \left\{ \sum\limits^{m_1+1}_{\mu_1=0}\sum\limits^{m_2+1}_{\mu_2=0}2^{\theta(\mu_1 \alpha_1+\mu_2\alpha_2)}
Y^\theta_{\floor{2^{\mu_1-1}},\floor{2^{\mu_2-1}}}(f)_{L_p(\T^2)} \right\}^{{p}/{\theta}}.
$$
By Lemma  \ref{l2.1}, we get

$$
I_9\ll  Y^p_{2^{m_1},2^{m_2}}(f)_{L_p(\T^2)}.
$$
Combining these estimates, we have
$$
I_5^p\ll 2^{-(m_1\alpha_1+m_2\alpha_2)p} \left\{
\sum\limits^{m_1+1}_{\mu_1=0}\sum\limits^{m_2+1}_{\mu_2=0}2^{\theta(\mu_1
\alpha_1+\mu_2\alpha_2)}
Y^\theta_{\floor{2^{\mu_1-1}},\floor{2^{\mu_2-1}}}(f)_{L_p(\T^2)}
\right\}^{{p}/{\theta}},
$$
and using properties of  $Y_{\nu_1,\nu_2}(f)_{L_p(\T^2)}$,
$$
I_5\ll \frac{1}{n_1^{\alpha_1}}\frac{1}{n_2^{\alpha_2}}\left\{
\sum\limits^{n_1+1}_{\nu_1=1}\sum\limits^{n_2+1}_{\nu_2=1}\nu_1^{\alpha_1\theta-1}
\nu_2^{\alpha_2\theta-1} Y^\theta_{\nu_1-1,\nu_2-1}(f)_{L_p(\T^2)}
\right\}^{{1}/{\theta}}.
$$
Thus, the proof of Theorem  \ref{th8.2} is complete.
\hfill $\Box$
\\[0.5cm]
{\bf Proof of Theorem \ref{th8.3}.} Define
$$
I^p:=\frac{1}{n_1^{\alpha_1p}}\frac{1}{n_2^{\alpha_2p}}
\sum\limits^{n_1+1}_{\nu_1=1}\sum\limits^{n_2+1}_{\nu_2=1}\nu_1^{\alpha_1p-1}
\nu_2^{\alpha_2p-1} Y^p_{\nu_1-1,\nu_2-1}(f)_{L_p(\T^2)}.
$$
By  Theorem \ref{ll2.1}, we have
$$
I^p\ll \frac{1}{n_1^{\alpha_1p}}\frac{1}{n_2^{\alpha_2p}}
\sum\limits^{n_1+1}_{\nu_1=1}\sum\limits^{n_2+1}_{\nu_2=1}\nu_1^{\alpha_1p-1}
\nu_2^{\alpha_2p-1} \omega^{p}_{\floor{\alpha_1}+1,\floor{\alpha_2}+1}\left(f,
\frac{1}{\nu_1+1}, \frac{1}{\nu_2+1} \right)_{L_p(\T^2)}.
$$
Theorem \ref{th7.2} implies

\begin{eqnarray*}
I^p
&\ll&
\frac{1}{n_1^{\alpha_1p}}\frac{1}{n_2^{\alpha_2p}}
\sum\limits^{n_1+1}_{\nu_1=1}\sum\limits^{n_2+1}_{\nu_2=1}
\frac{\nu_1^{\alpha_1p-1}
\nu_2^{\alpha_2p-1}
}
{(\nu_1+1)^{(\floor{\alpha_1}+1)p}(\nu_2+1)^{(\floor{\alpha_2}+1)p}}
\left[
\sum\limits^{\nu_1}_{\mu_1=1}\sum\limits^{\nu_2}_{\mu_2=1}
a^{p}_{\mu_1,\mu_2}
\mu_1^{(\floor{\alpha_1}+2)p-2}\mu_2^{(\floor{\alpha_2}+2)p-2}\right]
\\
&+&
\frac{1}{n_1^{\alpha_1p}}\frac{1}{n_2^{\alpha_2p}}  \sum\limits^{n_1+1}_{\nu_1=1}\sum\limits^{n_2+1}_{\nu_2=1}
 \frac{\nu_1^{\alpha_1p-1}
\nu_2^{\alpha_2p-1}}{(\nu_1+1)^{(\floor{\alpha_1}+1)p}}
\left[ \sum\limits^{\nu_1}_{\mu_1=1}\sum\limits^{\infty}_{\mu_2=\nu_2+1}a^{p}_{\mu_1,\mu_2} \mu_1^{(\floor{\alpha_1}+2)p-2}\mu_2^{p-2}\right]
\\
&+&
\frac{1}{n_1^{\alpha_1p}}\frac{1}{n_2^{\alpha_2p}}  \sum\limits^{n_1+1}_{\nu_1=1}\sum\limits^{n_2+1}_{\nu_2=1}
\frac{\nu_1^{\alpha_1p-1}
\nu_2^{\alpha_2p-1} }{(\nu_2+1)^{(\floor{\alpha_2}+1)p}}
\left[ \sum\limits^{\infty}_{\mu_1=\nu_1+1}\sum\limits^{\nu_2}_{\mu_2=1}a^{p}_{\mu_1,\mu_2} \mu_2^{(\floor{\alpha_1}+2)p-2}\mu_1^{p-2}\right]
\\
&+&\frac{1}{n_1^{\alpha_1p}}\frac{1}{n_2^{\alpha_2p}}  \sum\limits^{n_1+1}_{\nu_1=1}\sum\limits^{n_2+1}_{\nu_2=1}\nu_1^{\alpha_1p-1}
\nu_2^{\alpha_2p-1}
\left[ \sum\limits^{\infty}_{\mu_1=\nu_1+1}\sum\limits^{\infty}_{\mu_2=\nu_2+1}a^{p}_{\mu_1,\mu_2} (\mu_1\mu_2)^{p-2}\right]=: A_1+A_2+A_3+A_4.
\end{eqnarray*}
Changing the order of summation, we get
$$
A_1\ll  \frac{1}{n_1^{\alpha_1p}}\frac{1}{n_2^{\alpha_2p}}
\sum\limits^{n_1+1}_{\mu_1=1}\sum\limits^{n_2+1}_{\mu_2=1}
a^{p}_{\mu_1,\mu_2}\mu_1^{(\alpha_1+1)p-2}\mu_2^{(\alpha_2+1)p-2}=:
B_1.
$$
We proceed by estimating $A_2$.
\begin{eqnarray*}
A_2
&\ll&
\frac{1}{n_1^{\alpha_1p}}\frac{1}{n_2^{\alpha_2p}}
\sum\limits^{n_1+1}_{\nu_1=1}\sum\limits^{n_2+1}_{\nu_2=1}
\nu_1^{(\alpha_1-\floor{\alpha_1}-1)p-1} \nu_2^{\alpha_2p-1}
 \sum\limits^{\nu_1}_{\mu_1=1}\sum\limits^{n_2}_{\mu_2=\nu_2}
a^{p}_{\mu_1,\mu_2} \mu_1^{(\floor{\alpha_1}+2)p-2}\mu_2^{p-2}
\\
&+&
\frac{1}{n_1^{\alpha_1p}}\frac{1}{n_2^{\alpha_2p}}  \sum\limits^{n_1+1}_{\nu_1=1}\sum\limits^{n_2+1}_{\nu_2=1}
\nu_1^{(\alpha_1-\floor{\alpha_1}-1)p-1} \nu_2^{\alpha_2p-1}
 \sum\limits^{\nu_1}_{\mu_1=1}\sum\limits^{\infty}_{\mu_2=n_2+1}
a^{p}_{\mu_1,\mu_2} \mu_1^{(\floor{\alpha_1}+2)p-2}\mu_2^{p-2}
=:
A_{11}+A_{12}.
\end{eqnarray*}
Changing the order of summation, since
$\floor{\alpha_1}+1>\alpha_1$ and $\alpha_2>0$, we get
$$
A_{11}\ll  \frac{1}{n_1^{\alpha_1p}}\frac{1}{n_2^{\alpha_2p}}
\sum\limits^{n_1+1}_{\mu_1=1}\sum\limits^{n_2+1}_{\mu_2=1}
a^{p}_{\mu_1,\mu_2}\mu_1^{(\alpha_1+1)p-2}\mu_2^{(\alpha_2+1)p-2}= B_1.
$$
It is clear that
$$
A_{12}\ll
\frac{1}{n_1^{\alpha_1p}}  \sum\limits^{n_1+1}_{\nu_1=1}
\nu_1^{(\alpha_1-\floor{\alpha_1}-1)p-1}
 \sum\limits^{\nu_1}_{\mu_1=1}\sum\limits^{\infty}_{\mu_2=n_2+1}a^{p}_{\mu_1,\mu_2} \mu_1^{(\floor{\alpha_1}+2)p-2}\mu_2^{p-2}.
$$
Once again, changing the order of summation, we have
$$
A_{12}\ll
\frac{1}{n_1^{\alpha_1p}}
 \sum\limits^{n_1+1}_{\mu_1=1}\sum\limits^{\infty}_{\mu_2=n_2+1}a^{p}_{\mu_1,\mu_2} \mu_1^{(\alpha_1+1)p-2}\mu_2^{p-2}=: B_2.
$$
Thus, $A_2\ll B_1+B_2$ .

Similarly, we get $A_3\ll B_1+B_3$  and $A_4\ll B_1+B_2+B_3+B_4,$
where
$$
B_3:=\frac{1}{n_2^{\alpha_2p}}
 \sum\limits^{\infty}_{\mu_1=n_1+1}\sum\limits^{n_2+1}_{\mu_2=1}a^{p}_{\mu_1,\mu_2}\mu_1^{p-2} \mu_2^{(\alpha_1+1)p-2}
$$
and
$$
B_4:=
 \sum\limits^{\infty}_{\mu_1=n_1+1}\sum\limits^{\infty}_{\mu_2=n_2+1}a^{p}_{\mu_1,\mu_2}(\mu_1\mu_2)^{p-2}.
$$
Combining the estimates for $A_1$, $A_2$ , $A_3$, and $A_4$,
\begin{eqnarray*}
I^p
&\ll&
\frac{1}{n_1^{\alpha_1 p}n_2^{\alpha_2 p}} \sum\limits^{n_1}_{\nu_1=1}\sum\limits^{n_2}_{\nu_2=1}
a_{\nu_1,\nu_2}^p \nu_1^{(\alpha_1+1)p-2}\nu_2^{(\alpha_2+1)p-2}
+
\frac1{n_1^{\alpha_1 p}} \sum\limits^{n_1}_{\nu_1=1}\sum\limits^{\infty}_{\nu_2=n_2+1}
a_{\nu_1,\nu_2}^p \nu_1^{(\alpha_1+1)p-2}\nu_2^{p-2}
\\
&+&\frac{1}{n_2^{\alpha_2 p}} \sum\limits^{\infty}_{\nu_1=n_1+1}\sum\limits^{n_2}_{\nu_2=1}
a_{\nu_1,\nu_2}^p \nu_1^{p-2}\nu_2^{(\alpha_2+1)p-2}
+
\sum\limits^{\infty}_{\nu_1=n_1+1}\sum\limits^{\infty}_{\nu_2=n_2+1}
a_{\nu_1,\nu_2}^p (\nu_1\nu_2)^{p-2}.
\end{eqnarray*}
Theorem  \ref{th7.2} yields
$ I\ll  \omega_{\alpha_1,\alpha_2}\left(f, \frac1{n_1} , \frac1{n_2} \right)_{L_p(\T^2)}$,
i.e., the estimate from below is obtained.

To get the estimate from above, by Theorem  \ref{th7.2}, we have
\begin{eqnarray*}
 a_{n_1,n_2}^p(n_1n_2)^{p-1}&\ll&
\frac1{n_1^{(\alpha_1]+1)p}n_2^{(\floor{\alpha_2}+1])p}}
\sum\limits^{n_1}_{\nu_1=\floor{{n_1}/{2}}+1}
\sum\limits^{n_2}_{\nu_2=\floor{{n_2}/{2}}+1} a_{\nu_1,\nu_2}^p
\nu_1^{(\floor{\alpha_1}+2)p-2} \nu_2^{(\floor{\alpha_2}+2)p-2}
\\
& \ll&
\omega^p_{\floor{\alpha_1}+1,\floor{\alpha_2}+1}\left(f, \frac1{n_1} , \frac1{n_2} \right)_{L_p(\T^2)}.
\end{eqnarray*}
Moreover,
\begin{eqnarray*}
n_1^{p-1} \sum\limits^{\infty}_{\nu_2=n_2+1} a_{n_1,\nu_2}^p \nu_2^{p-2}
&\ll&  \frac1{n_1^{(\floor{\alpha_1}+1)p}}
\sum\limits^{n_1}_{\nu_1=\floor{\frac{n_1}{2}}+1} \sum\limits^{\infty}_{\nu_2=n_2+1}
a_{\nu_1,\nu_2}^p \nu_1^{(\floor{\alpha_1}+2)p-2} \nu_2^{p-2}
\\
&
\ll& \omega^p_{\floor{\alpha_1}+1,\floor{\alpha_2}+1}\left(f, \frac1{n_1},\frac1{n_2} \right)_{L_p(\T^2)},
\\n_2^{p-1} \sum\limits^{\infty}_{\nu_1=n_1+1} a_{\nu_1,n_2}^p \nu_1^{p-2}
&\ll&  \omega^p_{\floor{\alpha_1}+1,\floor{\alpha_2}+1}\left(f, \frac1{n_1} , \frac1{n_2} \right)_{L_p(\T^2)},
\end{eqnarray*}
and
\begin{eqnarray*}
&&
\sum\limits^{\infty}_{\nu_1=n_1+1} \sum\limits^{\infty}_{\nu_2=n_2+1} a_{\nu_1,\nu_2}^p
\nu_1^{p-2}\nu_2^{p-2} \ll  \omega^p_{\floor{\alpha_1}+1,\floor{\alpha_2}+1}\left(f, \frac1{n_1} , \frac1{n_2} \right)_{L_p(\T^2)}.
\end{eqnarray*}
Again, using  Theorem \ref{th7.2}, we have
\begin{eqnarray*}
\omega^p_{\alpha_1,\alpha_2}\left(f, \frac1{n_1} ,
\frac1{n_2} \right)_{L_p(\T^2)}
&\ll&
\frac{1}{n_1^{\alpha_1p}n_2^{\alpha_2p}}
\sum\limits^{n_1}_{\nu_1=1}\sum\limits^{n_2}_{\nu_2=1}
\nu_1^{\alpha_1p-1}\nu_2^{\alpha_2p-1}
(\nu_1^{p-1}\nu_2^{p-1}a_{\nu_1,\nu_2}^p)
\\
+
\frac{1}{n_1^{\alpha_1p}}
\sum\limits^{n_1}_{\nu_1=1}
\nu_1^{\alpha_1p-1}\Big(\nu_1^{p-1}
\sum\limits^{\infty}_{\nu_2=n_2+1}
a_{\nu_1,\nu_2}^p   \nu_2^{p-2}\Big)
&+&
\frac{1}{n_2^{\alpha_2p}}
\sum\limits^{n_2}_{\nu_2=1}
\nu_2^{\alpha_2p-1}
\Big(\nu_2^{p-1}
\sum\limits^{\infty}_{\nu_1=n_1+1}
a_{\nu_1,\nu_2}^p \nu_1^{p-2}\Big)
\\&+&
\sum\limits^{\infty}_{\nu_1=n_1+1} \sum\limits^{\infty}_{\nu_2=n_2+1}
a_{\nu_1,\nu_2}^p (\nu_1\nu_2)^{p-2}.
\end{eqnarray*}
Then applying the above mentioned inequalities, we have
\begin{eqnarray*}
\omega^p_{\alpha_1,\alpha_2}\left(f, \frac1{n_1} ,
\frac1{n_2} \right)_{L_p(\T^2)}
 &\ll& \frac{1}{n_1^{\alpha_1p}n_2^{\alpha_2p}}
\sum\limits^{n_1}_{\nu_1=1}\sum\limits^{n_2}_{\nu_2=1}
\nu_1^{\alpha_1p-1}\nu_2^{\alpha_2p-1}
\omega^p_{\floor{\alpha_1}+1,\floor{\alpha_2}+1}\left(f, \frac1{\nu_1} ,
\frac1{\nu_2} \right)_{L_p(\T^2)}
\\&+&
\frac{1}{n_1^{\alpha_1p}} \sum\limits^{n_1}_{\nu_1=1}
 \nu_1^{\alpha_1p-1} \omega^p_{\floor{\alpha_1}+1,\floor{\alpha_2}+1}\left(f, \frac1{\nu_1} , \frac1{n_2} \right)_{L_p(\T^2)}
\\&+&\frac{1}{n_2^{\alpha_2p}} \sum\limits^{n_2}_{\nu_2=1}
\nu_2^{\alpha_2p-1}\omega^p_{\floor{\alpha_1}+1,\floor{\alpha_2}+1}\left(f,
\frac1{n_1} , \frac1{\nu_2} \right)_{L_p(\T^2)}
\\&+&
\omega^p_{\floor{\alpha_1}+1,\floor{\alpha_2}+1}\left(f, \frac1{n_1} , \frac1{n_2} \right)_{L_p(\T^2)}.
\end{eqnarray*}
Properties of the mixed moduli of smoothness imply
$$
\omega^p_{\alpha_1,\alpha_2}\left(f, \frac1{n_1} ,
\frac1{n_2} \right)_{L_p(\T^2)}
 \ll \frac{1}{n_1^{\alpha_1p}n_2^{\alpha_2p}}
\sum\limits^{n_1}_{\nu_1=1}\sum\limits^{n_2}_{\nu_2=1}
\nu_1^{\alpha_1p-1}\nu_2^{\alpha_2p-1}
\omega^p_{\floor{\alpha_1}+1,\floor{\alpha_2}+1}\left(f, \frac1{\nu_1} ,
\frac1{\nu_2} \right)_{L_p(\T^2)}.
$$
Then, by Theorem \ref{ll2.1},
\begin{align*}
&
\omega^p_{\alpha_1,\alpha_2}\left(f, \frac1{n_1} ,
\frac1{n_2} \right)_{L_p(\T^2)}
 \\&\ll \frac{1}{n_1^{\alpha_1p}n_2^{\alpha_2p}}
\sum\limits^{n_1}_{\nu_1=1}\sum\limits^{n_2}_{\nu_2=1}
\frac{\nu_1^{\alpha_1p-1}\nu_2^{\alpha_2p-1}}{\nu_1^{(\floor{\alpha_1}+1)p}\nu_2^{(\floor{\alpha_2}+1)p}} \left(
\sum\limits^{\nu_1}_{\mu_1=1} \sum\limits^{\nu_2}_{\mu_2=1}
\mu_1^{\floor{\alpha_1}}\mu_2^{\floor{\alpha_2}}
Y_{\mu_1-1,\mu_2-1}(f)_{L_p(\T^2)}\right)^p.
\end{align*}
By Hardy's inequality (Lemma  \ref{l2.8}), we get
$$
\omega^p_{\alpha_1,\alpha_2}\left(f, \frac1{n_1} ,
\frac1{n_2} \right)_{L_p(\T^2)}
 \ll \frac{1}{n_1^{\alpha_1p}n_2^{\alpha_2p}}
\sum\limits^{n_1}_{\mu_1=1}\sum\limits^{n_2}_{\mu_2=1}
\mu_1^{\alpha_1p-1}\mu_2^{\alpha_2p-1}
Y^p_{\mu_1-1,\mu_2-1}(f)_{L_p(\T^2)}.
$$
Thus, the proof of Theorem  \ref{th8.3} is now complete.
\hfill $\Box$
\\[0.5cm]
{\bf Proof of Theorem  \ref{th8.4}.}
Here, it is enough to show that
$
\omega_{\alpha_1,\alpha_2}\left( f, \frac{1}{2^{m_1}},
\frac{1}{2^{m_2}} \right)_{L_p(\T^2)}\asymp
I$, where $$I:= 2^{-(m_1\alpha_1+m_2\alpha_2)} \left\{ \sum\limits^{m_1}_{\mu_1=1}\sum\limits^{m_2}_{\mu_2=1}
2^{2(\mu_1\alpha_1+\mu_2\alpha_2)} Y^{2}_{2^{\mu_1}-1,2^{\mu_2}-1}(f)_{L_p(\T^2)} \right\}^{1/2}.
$$
Lemmas  \ref{l2.1} and \ref{l2.7} imply

\begin{eqnarray*}
I &\asymp& 2^{-(m_1\alpha_1+m_2\alpha_2)}
\\
&&
\left\{
\sum\limits^{m_1}_{\mu_1=1}\sum\limits^{m_2}_{\mu_2=1}
2^{2(\mu_1\alpha_1+\mu_2\alpha_2)} \|f- s_{2^{\mu_1}-1, \infty}(f)-s_{\infty, 2^{\mu_2}-1}(f)+
s_{2^{\mu_1}-1,2^{\mu_2}-1}(f)\|^2_{L_p(\T^2)} \right\}^{1/2} \\
&\asymp&
2^{-(m_1\alpha_1+m_2\alpha_2)} \left\{ \sum\limits^{m_1}_{\mu_1=1}\sum\limits^{m_2}_{\mu_2=1}
2^{2(\mu_1\alpha_1+\mu_2\alpha_2)} \sum\limits^{\infty}_{\nu_1=\mu_1} \sum\limits^{\infty}_{\nu_2=\mu_2}
\lambda^2_{\mu_1,\mu_2} \right\}^{1/2}
\\
&\asymp&
2^{-(m_1\alpha_1+m_2\alpha_2)} \left\{ \sum\limits^{m_1}_{\mu_1=1}\sum\limits^{m_2}_{\mu_2=1}  \lambda^2_{\mu_1,\mu_2}
2^{2(\mu_1\alpha_1+\mu_2\alpha_2)}  \right\}^{1/2}
+ 2^{-m_1\alpha_1}\left\{ \sum\limits^{m_1}_{\mu_1=1}\sum\limits^{\infty}_{\mu_2=m_2+1}  \lambda^2_{\mu_1,\mu_2}
2^{2\mu_1\alpha_1}  \right\}^{1/2}
\\&+&
 2^{-m_2\alpha_2} \left\{ \sum\limits^{\infty}_{\mu_1=m_1+1}\sum\limits^{m_2}_{\mu_2=1}  \lambda^2_{\mu_1,\mu_2}
2^{2\mu_2\alpha_2}  \right\}^{1/2}
+\left\{ \sum\limits^{\infty}_{\mu_1=m_1+1} \sum\limits^{\infty}_{\mu_2=m_2+1} \lambda^2_{\mu_1,\mu_2}  \right\}^{1/2}.
\end{eqnarray*}
Using now Theorem  \ref{th7.3}, we finally get  $ I \asymp
\omega_{\alpha_1,\alpha_2}\left( f, \frac{1}{2^{m_1}},
\frac{1}{2^{m_2}} \right)_{L_p(\T^2)}$.
\hfill $\Box$

Theorem  \ref{th8.2} was proved in the case of  $\alpha_i \in \N$  in the paper \cite{25}.
For the proof of Theorem \ref{th8.3}, see \cite{23}.
The one-dimensional sharp Jackson inequality was proved in \cite{timan}.
The history of this topic and new results can be found in \cite{jat-dai}.

\vskip 1.5cm

\sect{Interrelation between the mixed moduli of smoothness of different orders in $L_p$}

We have stated above (see Theorem \ref{th6.1}) that the following properties of mixed moduli are well known:
\begin{eqnarray*}
  \qquad \omega_{\beta_1,\beta_2}(f, \delta_1,\delta_2)_{L_p(\T^2)} &\ll&  \omega_{\alpha_1, \alpha_2}(f,
\delta_1,\delta_2)_{L_p(\T^2)};
\\
\qquad  \frac{\omega_{\alpha_1, \alpha_2}(f,\delta_1,\delta_2)_{L_p(\T^2)}}{ \delta_1^{\alpha_1}\delta_2^{\alpha_2}} &\ll&
\frac{\omega_{\beta_1,\beta_2}(f, \delta_1,\delta_2)_{L_p(\T^2)} }{ \delta_1^{\beta_1}\delta_2^{\beta_2}};
\\
{\omega_{\alpha_1,\alpha_2}(f,\delta_1,\delta_2)_{L_p(\T^2)}} &\ll&
\delta_1^{\alpha_1} \delta_2^{\alpha_2}\int_{\delta_1}^1\int_{\delta_2}^1
\frac{\omega_{\beta_1,\beta_2}(f,t_1,t_2)_{L_p(\T^2)}}{t_1^{\alpha_1}t_2^{\alpha_2}}
\frac{\dd t_1}{t_1}\frac{\dd t_2}{t_2},
\end{eqnarray*}
where $f \in L_p^0(\T^2)$, $1 < p < \infty,$
and  $0 < \alpha_i < \beta_i$, $\delta_i \in (0, \frac12) , i=1,2$.

The following theorem sharpens all these estimates.
\begin{theorem}\label{th10.2}
Let $f \in L_p^0(\T^2)$, $1 < p < \infty,$
$\tau := \max(2,p)$, $\theta := \min(2,p),$ $0 < \alpha_i < \beta_i,$ $\delta_i \in (0, \frac12)$, $i=1,2$. Then
$$ \delta_1^{\alpha_1}\delta_2^{\alpha_2} \left\{ \int\limits_{\delta_1}^1 \int\limits_{\delta_2}^1 \Big[ t_1^{- \alpha_1} t_2^{- \alpha_2}
\omega_{\beta_1,\beta_2}(f,t_1,t_2)_{L_p(\T^2)} \Big]^{\tau} \frac{\dd t_1}{t_1} \frac{\dd t_2}{t_2}
\right\}^{{1}/{\tau}} \ll \omega_{\alpha_1, \alpha_2}(f,\delta_1,\delta_2)_{L_p(\T^2)}
$$
$$
\ll
\delta_1^{\alpha_1}\delta_2^{\alpha_2}  \left\{ \int\limits_{\delta_1}^1 \int\limits_{\delta_2}^1 \Big[ t_1^{- \alpha_1} t_2^{- \alpha_2}
\omega_{\beta_1,\beta_2}(f,t_1,t_2)_{L_p(\T^2)} \Big]^{\theta} \frac{\dd t_1}{t_1} \frac{\dd t_2}{t_2}
\right\}^{{1}/{\theta}}.
$$
\end{theorem}
These estimates are the best possible in the sense of order. This can be shown using the following two equivalence results for  the
function classes $M_p$ and $\Lambda_p$ defined in Section 2.4.
\begin{theorem}\label{th10.3}
Let $f \in M_p, 1 < p < \infty,$ $0 < \alpha_i < \beta_i,$ $\delta_i \in (0, \frac12), i=1,2$. Then

$$  \omega_{\alpha_1,\alpha_2}(f,\delta_1,\delta_2)_{L_p(\T^2)} \asymp
 \delta_1^{\alpha_1}\delta_2^{\alpha_2} \left\{ \int\limits_{\delta_1}^1 \int\limits_{\delta_2}^1 \Big[ t_1^{- \alpha_1} t_2^{- \alpha_2}
\omega_{\beta_1,\beta_2}(f,t_1,t_2)_{L_p(\T^2)} \Big]^{p} \frac{\dd t_1}{t_1} \frac{\dd t_2}{t_2}
\right\}^{{1}/{p}}.
$$
\end{theorem}

\begin{theorem}\label{th10.4}
Let  $f \in \Lambda_p, 1 < p < \infty,$ $0 < \alpha_i < \beta_i,$ $\delta_i \in (0, \frac12), i=1,2$. Then

$$  \omega_{\alpha_1,\alpha_2}(f,\delta_1,\delta_2)_{L_p(\T^2)} \asymp
 \delta_1^{\alpha_1}\delta_2^{\alpha_2} \left\{ \int\limits_{\delta_1}^1 \int\limits_{\delta_2}^1 \Big[ t_1^{- \alpha_1} t_2^{- \alpha_2}
\omega_{\beta_1,\beta_2}(f,t_1,t_2)_{L_p(\T^2)} \Big]^{2} \frac{\dd t_1}{t_1} \frac{\dd t_2}{t_2}
\right\}^{{1}/{2}}.
$$
\end{theorem}
{\bf Proof of Theorem \ref{th10.2}.}
We consider
$$
I:=\delta_1^{\alpha_1}\delta_2^{\alpha_2} \left( \int\limits_{\delta_1}^1 \int\limits_{\delta_2}^1 \Big[ t_1^{- \alpha_1} t_2^{- \alpha_2}
\omega_{\beta_1,\beta_2}(f,t_1,t_2)_{L_p(\T^2)} \Big]^{\tau} \frac{\dd t_1}{t_1} \frac{\dd t_2}{t_2}\right)^{1/\tau}.
$$
For given $\delta_i \in (0,\frac12)$, we take integers $n_i$ such that $\frac{1}{n_i+1} \le \delta_i < \frac1{n_i}, i=1,2$.
Then
$$
I^{\tau}\ll  n_1^{-\alpha_1 \tau}n_2^{-\alpha_2 \tau} \sum\limits^{n_1}_{\mu_1=1} \sum\limits^{n_2}_{\mu_2=1}
\mu_1^{\alpha_1\tau-1} \mu_2^{\alpha_2\tau-1} \omega_{\beta_1,\beta_2}^{\tau}
\left(f,\frac{1}{\mu_1},\frac{1}{\mu_1} \right)_{L_p(\T^2)}.
$$
Further, Theorem \ref{th8.2} implies
$$
I^{\tau}\ll  n_1^{-\alpha_1 \tau}n_2^{-\alpha_2 \tau} \sum\limits^{n_1+1}_{\mu_1=1} \sum\limits^{n_2+1}_{\mu_2=1}
\mu_1^{\alpha_1\tau-\beta_1\tau-1} \mu_2^{\alpha_2\tau-\beta_2\tau-1}
\left\{\sum\limits^{\mu_1+1}_{\nu_1=1}\sum\limits^{\mu_2+1}_{\nu_2=1} \nu_1^{\beta_1\theta-1} \nu_2^{\beta_2\theta-1}
Y^{\theta}_{\nu_1-1,\nu_2-1}(f)_{L_p(\T^2)}\right\}^{{\tau}/{\theta}}.
$$
Since ${\tau}/{\theta} \geq 1$, using Lemma \ref{l2.8}, we get
$$
I^{\tau}\ll  n_1^{-\alpha_1 \tau}n_2^{-\alpha_2 \tau} \sum\limits^{n_1+1}_{\mu_1=1} \sum\limits^{n_2+1}_{\mu_2=1}
\mu_1^{\alpha_1\tau-1} \mu_2^{\alpha_2\tau-1}
Y^{\tau}_{\mu_1-1,\mu_2-1}(f)_{L_p(\T^2)}.
$$
By Theorem \ref{th8.2}, we have
$$
I^{\tau}\ll  \omega_{\alpha_1,\alpha_2}^{\tau} \left(f,\frac{1}{n_1}, \frac{1}{n_2} \right)_{L_p(\T^2)},
$$
which implies
$$
I\ll  \omega_{\alpha_1,\alpha_2} \left(f,\frac{1}{n_1+1}, \frac{1}{n_2+1} \right)_{L_p(\T^2)}\le  \omega_{\alpha_1,\alpha_2}(f,\delta_1,\delta_2)_p.
$$
The estimate ``$\ll  \omega_{\alpha_1,\alpha_2}(f,\delta_1,\delta_2)_p$'' is proved.

Let us verify the part ``$ \omega_{\alpha_1,\alpha_2}(f,\delta_1,\delta_2)_p\ll$''.
For given $\delta_i\in (0,\frac12)$, we take integers  $n_i$ such that $\frac{1}{n_i+1}\le \delta_i <\frac1{n_i}, i=1,2$.
 Then, by Theorem \ref{th8.2}, we get
\begin{eqnarray*}
\omega_{\alpha_1,\alpha_2}(f,\delta_1,\delta_2)_{L_p(\T^2)}
&\ll&
\omega_{\alpha_1,\alpha_2}\left(f,\frac{1}{n_1},\frac{1}{n_2} \right)_{L_p(\T^2)}
\\&\ll&  n_1^{-\alpha_1}n_2^{-\alpha_2}
\left\{ \sum\limits^{n_1+1}_{\nu_1=1} \sum\limits^{n_2+1}_{\nu_2=1} \nu_1^{\alpha_1\theta-1}
\nu_2^{\alpha_2\theta-1} Y^{\theta}_{\nu_1-1,\nu_2-1}(f)_{L_p(\T^2)}\right\}^{{1}/{\theta}}.
\end{eqnarray*}
Then Theorem \ref{ll2.1} gives us
\begin{eqnarray*}
\omega_{\alpha_1,\alpha_2}(f,\delta_1,\delta_2)_{L_p(\T^2)} &\ll&
n_1^{-\alpha_1}n_2^{-\alpha_2} \left\{ \sum\limits^{n_1+1}_{\nu_1=1}
\sum\limits^{n_2+1}_{\nu_2=1} \nu_1^{\alpha_1\theta-1}
\nu_2^{\alpha_2\theta-1}
\omega^{\theta}_{\beta_1,\beta_2}\left(f,\frac{1}{\nu_1},\frac{1}{\nu_2}
\right)_{L_p(\T^2)}\right\}^{{1}/{\theta}}
\\
&\ll&  \delta_1^{\alpha_1}\delta_2^{\alpha_2} \left\{
\int\limits^{1}_{\delta_1} \int\limits^{1}_{\delta_2} \Big[
t_1^{-\alpha_1}t_2^{-\alpha_2}
\omega_{\beta_1,\beta_2}\left(f,t_1,t_2 \right)_{L_p(\T^2)}
\Big]^{\theta} \frac{\dd t_1}{t_1} \frac{\dd t_2}{t_2}
\right\}^{{1}/{\theta}},
\end{eqnarray*}
which completes the proof of Theorem \ref{th10.2}.
\hfill $\Box$
\\[0.5cm]
{\bf Proof of Theorem  \ref{th10.3}.}
Denoting
$$
A:=\delta_1^{\alpha_1}\delta_2^{\alpha_2} \left\{
\int\limits^{1}_{\delta_1} \int\limits^{1}_{\delta_2} \Big[
t_1^{-\alpha_1}t_2^{-\alpha_2}
\omega_{\beta_1,\beta_2}\left(f,t_1,t_2 \right)_{L_p(\T^2)} \Big]^{p}
\frac{\dd t_1}{t_1} \frac{\dd t_2}{t_2} \right\}^{{1}/{p}},
$$
and choosing integers $n_i$ such that
 $\frac{1}{n_i+1}\le \delta_i< \frac{1}{n_i}, i=1,2$,
we have
$$
A^p\asymp n_1^{-\alpha_1p}n_2^{-\alpha_2p}
\sum\limits^{n_1}_{\mu_1=1} \sum\limits^{n_2}_{\mu_2=1}
\mu_1^{\alpha_1p-1} \mu_2^{\alpha_2p-1}
\omega^{p}_{\beta_1,\beta_2}\big(f,\frac{1}{\mu_1},\frac{1}{\mu_2}
\big)_{L_p(\T^2)}.
$$
Then Theorem  \ref{th8.3} yields
$$
A^p\asymp n_1^{-\alpha_1p}n_2^{-\alpha_2p}  \sum\limits^{n_1}_{\mu_1=1} \sum\limits^{n_2}_{\mu_2=1} \mu_1^{\alpha_1p-\beta_1p-1}
\mu_2^{\alpha_2p-\beta_2p-1} \sum\limits^{\mu_1+1}_{\nu_1=1}\sum\limits^{\mu_2+1}_{\nu_2=1}
\nu_1^{\beta_1p-1}\nu_2^{\beta_2p-1} Y^p_{\nu_1-1,\nu_2-1}(f)_{L_p(\T^2)}
$$
$$
\asymp  n_1^{-\alpha_1p}n_2^{-\alpha_2p}
\sum\limits^{n_1+1}_{\nu_1=1}\sum\limits^{n_2+1}_{\nu_2=1}
\nu_1^{\alpha_1p-1}\nu_2^{\alpha_2p-1}
Y^p_{\nu_1-1,\nu_2-1}(f)_{L_p(\T^2)}.
$$
Finally, by Theorem \ref{th8.3}, we have
$A^p \asymp
\omega_{\alpha_1,\alpha_2} \big(f, \frac1{n_1},\frac{1}{n_2}
\big)_{L_p(\T^2)}$, which completes the proof.
\hfill $\Box$
\\[0.5cm]
{\bf Proof of Theorem \ref{th10.4}.}
As above, let
$$
B:=\delta_1^{\alpha_1}\delta_2^{\alpha_2} \left\{
\int\limits^{1}_{\delta_1} \int\limits^{1}_{\delta_2} \Big[
t_1^{-\alpha_1}t_2^{-\alpha_2}
\omega_{\beta_1,\beta_2}\left(f,t_1,t_2 \right)_{L_p(\T^2)} \Big]^{2}
\frac{\dd t_1}{t_1} \frac{\dd t_2}{t_2} \right\}^{{1}/{2}}.
$$
and
  $\frac{1}{2^{n_i+1}}\le \delta_i<
\frac{1}{2^{n_i}}, i=1,2$. Then
$$
B^2\asymp 2^{-n_1\alpha_1}2^{-n_2\alpha_2}  \sum\limits^{n_1}_{\mu_1=1} \sum\limits^{n_2}_{\mu_2=1} 2^{2\mu_1\alpha_1+2\mu_2\alpha_2}
 \omega^{2}_{\beta_1,\beta_2}\left(f,\frac{1}{2^{\mu_1}},\frac{1}{2^{\mu_2}} \right)_{L_p(\T^2)}.
$$
Using now Theorem \ref{th8.4}, we have
$$
B^2\asymp 2^{-n_1\alpha_1}2^{-n_2\alpha_2}  \sum\limits^{n_1}_{\mu_1=1} \sum\limits^{n_2}_{\mu_2=1} 2^{2\mu_1\alpha_1+2\mu_2\alpha_2-2\mu_1\beta_1-2\mu_2\beta_2}
\sum\limits^{\mu_1+1}_{\nu_1=1}\sum\limits^{\mu_2+1}_{\nu_2=1} 2^{2\nu_1\beta_1+2\nu_2\beta_2}
Y^{2}_{2^{\nu_1}-1,2^{\nu_2}-1}(f)_{L_p(\T^2)}
$$
$$
\asymp 2^{-2(n_1\alpha_1+n_2\alpha_2)}
\sum\limits^{n_1+1}_{\nu_1=1} \sum\limits^{n_2+1}_{\nu_2=1}
2^{2(\nu_1\alpha_1+\nu_2\alpha_2)}
Y^{2}_{2^{\nu_1}-1,2^{\nu_2}-1}(f)_{L_p(\T^2)}.
$$
By Theorem  \ref{th8.4}, we get $B \asymp
\omega_{\alpha_1,\alpha_2}\left(f,\frac{1}{2^{n_1}},
\frac{1}{2^{n_2}} \right)_{L_p(\T^2)}\asymp \omega_{\alpha_1,\alpha_2}\Big(
f,\delta_1,\delta_2 \Big)_{L_p(\T^2)}$.
\hfill $\Box$

The one-dimensional version of Theorem \ref{th10.3} for functions with general monotone coefficients was stated in \cite{16}.

\vskip 1.5cm

\sect{Interrelation between the mixed moduli of smoothness
 in various ($L_p$,$L_q$) metrics}

For functions on $\T$, the classical Ul'yanov inequality (\cite{14}, see also \cite{dit-ul}) states that
$$
 \omega_{\alpha }(f, \delta)_{L_q(\T)} \ll \left( \int\limits_0^{\delta}
\Big[ t^{- \theta} \omega_{\alpha
}(f, t)_{L_p(\T)} \Big]^q
\frac{\dd t}{t}  \right)^{{1}/{q}},
$$
where
$f \in L_p(\T), 1 < p < q < \infty,$
$\theta := \frac{1}{p} - \frac{1}{q}$, and  $ \alpha\in \mathbb{N}.$
Very recently, this inequality was generalized using the fractional modulus of smoothness ($\alpha>0$) as follows
(see \cite{sim-2, sim-jat, trebels}):
$$
 \omega_{\alpha }(f, \delta)_{L_q(\T)} \ll \left( \int\limits_0^{\delta}
\Big[ t^{- \theta} \omega_{\alpha +\theta
}(f, t)_{L_p(\T)} \Big]^q
\frac{\dd t}{t}  \right)^{{1}/{q}}.
$$
Below, we present an analogue of the sharp Ul'yanov inequality for the mixed moduli of smoothness (see \cite{26}).
\begin{theorem}\label{th11.1}
Let $f \in L_p^0(\T^2), 1 < p < q < \infty,$
$\theta := \frac{1}{p} - \frac{1}{q},  \alpha_1>0,\alpha_2>0.$
Then
\begin{equation}\label{f10.1}
 \omega_{\alpha_1,\alpha_2}(f, \delta_1,\delta_2)_{L_q(\T^2)} \ll \left( \int\limits_0^{\delta_1} \int\limits_0^{\delta_2}
\Big[ (t_1t_2)^{- \theta} \omega_{\alpha_1 +\theta,\alpha_2+ \theta}(f, t_1,t_2)_{L_p(\T^2)} \Big]^q
\frac{\dd t_1}{t_1} \frac{\dd t_2}{t_2} \right)^{{1}/{q}}.
\end{equation}

\end{theorem}
{\bf Proof.}
For given  $\delta_i\in (0,1)$, there exist integers $n_i$  such that
  $\frac{1}{2^{n_i+1}}\le \delta_i< \frac1{2^{n_i}}, i=1,2$.
Then, by Theorem \ref{th4.1}, we get
\begin{eqnarray*}
I&:=&\omega_{\alpha_1,\alpha_2}(f,\delta_1,\delta_2)_{L_q(\T^2)}\ll
\omega_{\alpha_1,\alpha_2} \left(f,\frac1{2^{n_1}},\frac1{2^{n_2}}
\right)_{L_q(\T^2)}
\\
&\ll &
\|f-s_{2^{n_1},\infty}-s_{\infty,2^{n_2}}+s_{2^{n_1},2^{n_2}}\|_{L_q(\T^2)}
\\&+&
2^{-\alpha_1n_1}\|s^{(\alpha_1,0)}_{2^{n_1},\infty}(f-s_{\infty,2^{n_2}})\|_{L_q(\T^2)}+2^{-\alpha_2n_2}
\|s_{\infty,2^{n_2}}^{(0,\alpha_2)}(f-s_{2^{n_1},\infty})\|_{L_q(\T^2)}
\\ &+&
2^{-\alpha_1n_1-\alpha_2n_2}\|s^{(\alpha_1,\alpha_2)}_{2^{n_1},2^{n_2}}(f)\|_{L_q(\T^2)}=:
I_1+I_2+I_3+I_4.
\end{eqnarray*}
To estimate  $I_1$, we use Lemmas  \ref{l2.1} and \ref{l2.10}
$$
I_1\ll  \left\{ \sum\limits^{\infty}_{\nu_1=n_1}\sum\limits^{\infty}_{\nu_2=n_2}
2^{(\nu_1+\nu_2)\theta q}Y^{q}_{2^{\nu_1}-1,2^{\nu_2}-1}(f)_{L_p(\T^2)} \right\}^{{1}/{q}}.
$$
Further, by
Theorem  \ref{ll2.1}, we get
$$
I_1\ll  \left\{ \sum\limits^{\infty}_{\nu_1=n_1}\sum\limits^{\infty}_{\nu_2=n_2}
2^{(\nu_1+\nu_2)\theta q}\omega^{q}_{\alpha_1+\theta,\alpha_2+\theta}
\left(f,\frac{1}{2^{\nu_1}},\frac{1}{2^{\nu_2}} \right)_{L_p(\T^2)} \right\}^{{1}/{q}}
=: J.
$$
Let us now estimate $I_2$. Define
$\varphi(x,y):=s^{(\alpha_1,0)}_{2^{n_1},\infty}(f)$. Then $I_2=
2^{-\alpha_1n_1}\|\varphi-s_{\infty,2^{n_2}}(\varphi)\|_{L_q(\T^2)}$.
Lemma  \ref{l3.5} implies, for a.e. $x$,
$$
\int\limits^{2\pi}_{0} |\varphi-s_{\infty,2^{n_2}}(\varphi)|^q \dd y\ll
 \sum\limits^{\infty}_{\nu_2=n_2} 2^{\nu_2\theta q}
\left( \int\limits^{2\pi}_{0} |\varphi-s_{\infty,2^{\nu_2}}(\varphi)|^p \dd y \right)^{{q}/{p}}.
$$
Therefore,
$$
\int\limits^{2\pi}_{0}\int\limits^{2\pi}_{0} |\varphi-s_{\infty,2^{n_2}}(\varphi)|^q \dd y\dd x\ll  \sum\limits^{\infty}_{\nu_2=n_2} 2^{\nu_2\theta q}
\int\limits^{2\pi}_{0}
\left( \int\limits^{2\pi}_{0} |\varphi-s_{\infty,2^{\nu_2}}(\varphi)|^p \dd y \right)^{{q}/{p}} \dd x.
$$
Using Minkowski's inequality, we have
$$
\int\limits^{2\pi}_{0} \left( \int\limits^{2\pi}_{0}
|\varphi-s_{\infty,2^{\nu_2}}(\varphi)|^p \dd y \right)^{{q}/{p}}
\dd x\ll  \left( \int\limits^{2\pi}_{0} \left( \int\limits^{2\pi}_{0}
|\varphi-s_{\infty,2^{\nu_2}}(\varphi)|^q \dd x\right)^{{p}/{q}}
\dd y \right)^{{q}/{p}}.
$$
Since $\varphi-s_{\infty,2^{\nu_2}}(\varphi)=s^{(\alpha_1,0)}_{2^{n_1},
\infty}(f-s_{\infty,2^{\nu_2}}(f))$, then
Lemma \ref{l3.6} implies, for a.e.  $y$,
$$
\left( \int\limits^{2\pi}_{0} |s^{(\alpha_1,0)}_{2^{n_1},
\infty}(f-s_{\infty,2^{\nu_2}}(f))|^q \dd x\right)^{{1}/{q}} \ll
\left( \int\limits^{2\pi}_{0} |s^{(\alpha_1+\theta,0)}_{2^{n_1},
\infty}(f-s_{\infty,2^{\nu_2}}(f))|^p \dd x\right)^{{1}/{p}}.
$$
This gives
$$
\int\limits^{2\pi}_{0} \left( \int\limits^{2\pi}_{0}
|s^{(\alpha_1,0)}_{2^{n_1},\infty}(f-s_{\infty,2^{\nu_2}}(f))|^q \dd x\right)^{{p}/{q}} \dd y  \ll
\int\limits^{2\pi}_{0} \int\limits^{2\pi}_{0}
|s^{(\alpha_1+\theta,0)}_{2^{n_1},\infty}(f-s_{\infty,2^{\nu_2}}(f))|^p \dd x\dd y.
$$
Thus,
\begin{eqnarray*}
I_2
&\ll&  2^{-\alpha_1n_1} \left\{ \sum\limits^{\infty}_{\nu_2=n_2}
2^{\nu_2\theta q} \|
s^{(\alpha_1+\theta,0)}_{2^{n_1},\infty}(f-s_{\infty,2^{\nu_2}}(f))\|^q_{L_p(\T^2)}
\right\}^{1/{q}}
\\
&=&
2^{n_1\theta}
 \left\{ \sum\limits^{\infty}_{\nu_2=n_2} 2^{\nu_2\theta q}
 \Big[ 2^{-(\alpha_1+\theta)n_1} \| s^{(\alpha_1+\theta,0)}_{2^{n_1},
 \infty}(f-s_{\infty,2^{\nu_2}}(f))\|_{L_p(\T^2)}
\Big]^q  \right\}^{1/{q}}.
\end{eqnarray*}
Now using Theorem \ref{th4.1}, we get

$$
I_2\ll  \left\{ 2^{n_1\theta q} \sum\limits^{\infty}_{\nu_2=n_2} 2^{\nu_2\theta q}
\omega^{q}_{\alpha_1+\theta,\alpha_2+\theta}\left(f,\frac{1}{2^{n_1}},\frac{1}{2^{\nu_2}} \right)_{L_p(\T^2)}
 \right\}^{1/q}\ll  J.
$$
Similarly, one can show that  $I_3\ll J$  and  $I_4\ll J$. Combining the estimates for $I_1,I_2,I_3,$ and $I_4$, we get
$$
I\ll \left\{\sum\limits^{\infty}_{\nu_1=n_1} \sum\limits^{\infty}_{\nu_2=n_2} 2^{(\nu_1+\nu_2)\theta q}
\omega^{q}_{\alpha_1+\theta,\alpha_2+\theta}\left(f,\frac{1}{2^{\nu_1}},\frac{1}{2^{\nu_2}} \right)_{L_p(\T^2)}
 \right\}^{1/q}.
$$
Taking into account monotonicity properties of the mixed moduli of smoothness, we finally have
$$
I\ll  \left\{ \int\limits^{2^{\frac1{n_1+1}}}_{0}
\int\limits^{2^{\frac1{n_2+1}}}_{0} \Big[(t_1t_2)^{-\theta}
\omega_{\alpha_1+\theta,\alpha_2+\theta}\left(f,t_1,t_2
\right)_{L_p(\T^2)}  \Big]^q \frac{\dd t_1}{t_1}\frac{\dd t_2}{t_2}
\right\}^{1/q}\ll
$$
$$
\ll  \left\{ \int\limits^{\delta_1}_{0} \int\limits^{\delta_2}_{0}
\Big[(t_1t_2)^{-\theta}
\omega_{\alpha_1+\theta,\alpha_2+\theta}\left(f,t_1,t_2
\right)_{L_p(\T^2)}  \Big]^q \frac{\dd t_1}{t_1}\frac{\dd t_2}{t_2}
\right\}^{1/q}.
$$
\hfill $\Box$

\vskip 0.5cm

\subsection{Sharpness}
Let us show that it is impossible to obtain the reverse part  $\gg$ of inequality (\ref{f10.1}), i.e., in general, (\ref{f10.1}) is not an equivalence.
%Let us construct the function  $f_0(x,y)\in
Let us construct a function  $f_0(x,y)\in
L^0_p(\T^2)$ such that  the terms on the left- and right-hand side of (\ref{f10.1})
 have different orders as functions of $\delta_1$ and $\delta_2$.
    Consider
$$
f_0(x,y):=\sum\limits^{\infty}_{\nu_1=0} \sum\limits^{\infty}_{\nu_2=0} a_{\nu_1,\nu_2}
\cos{2^{\nu_1}x} \cos{2^{\nu_2}y},
\; \textnormal{where} \;
a_{\nu_1,\nu_2}:=\frac{(\nu_1+1)^{\beta_1}(\nu_2+1)^{\beta_2}}{2^{\alpha_1\nu_1+\alpha_2\nu_2}},
$$
$1<p<q< \infty$,
$\theta:=\frac{1}{p}-\frac{1}{q}$,
$\beta_i >-\frac{1}{2}, \alpha_i>\theta$, $(i=1,2)$.
By Theorem  \ref{th7.3}, we get
$$
\omega_{\alpha_1,\alpha_2}\left(f_0,\frac{1}{2^{n_1}},\frac{1}{2^{n_2}}
\right)_{L_q(\T^2)} \asymp
\frac{(n_1+1)^{\beta_1+\frac12}(n_2+1)^{\beta_2+\frac12}}{2^{\alpha_1\nu_1+\alpha_2\nu_2}}
$$
and
$$
\omega_{\alpha_1+\theta,\alpha_2+\theta}\left(f_0,\frac{1}{2^{n_1}},\frac{1}{2^{n_2}}
\right)_{L_p(\T^2)} \asymp
\frac{(n_1+1)^{\beta_1}(n_2+1)^{\beta_2}}{2^{\alpha_1\nu_1+\alpha_2\nu_2}}.
$$
Using these estimates, one can easily see that
$$
\omega_{\alpha_1,\alpha_2}(f_0,\delta_1,\delta_2
)_{L_q(\T^2)}\asymp \delta_1^{\alpha_1} \delta_2^{\alpha_2}
\left(\ln{\frac{2}{\delta_1}} \right)^{\beta_1+\frac12}
\left(\ln{\frac{2}{\delta_2}} \right)^{\beta_2+\frac12}
$$
and
\begin{equation}\label{f10.1111}
  \left\{ \int\limits^{\delta_1}_{0} \int\limits^{\delta_2}_{0}
\Big[(t_1t_2)^{-\theta}
\omega_{\alpha_1+\theta,\alpha_2+\theta}\left(f_0,t_1,t_2
\right)_{L_p(\T^2)}  \Big]^q \frac{\dd t_1}{t_1}\frac{\dd t_2}{t_2}
\right\}^{1/q}
\asymp
\delta_1^{\alpha_1-\theta}
\delta_2^{\alpha_2-\theta} \left(\ln{\frac{2}{\delta_1}}
\right)^{\beta_1} \left(\ln{\frac{2}{\delta_2}} \right)^{\beta_2}.
\end{equation}
Thus, the inequality
$$ \omega_{\alpha_1,\alpha_2}(f, \delta_1,\delta_2)_{L_q(\T^2)} \gtrsim \left( \int\limits_0^{\delta_1} \int\limits_0^{\delta_2}
\Big[ (t_1t_2)^{- \theta} \omega_{\alpha_1 +\theta,\alpha_2+ \theta}(f, t_1,t_2)_{L_p(\T^2)} \Big]^q
\frac{\dd t_1}{t_1} \frac{\dd t_2}{t_2} \right)^{{1}/{q}}
$$
does not hold for $f_0$.

Finally, we would like to mention that inequality (\ref{f10.1}) improves the classical Ul'yanov in\-equality for the mixed moduli of smoothness given by (\cite{9})
\begin{equation}\label{f10.2}
 \omega_{\alpha_1,\alpha_2}(f, \delta_1,\delta_2)_{L_q(\T^2)} \ll \left( \int\limits_0^{\delta_1} \int\limits_0^{\delta_2}
\Big[ (t_1t_2)^{- \theta} \omega_{\alpha_1,\alpha_2}(f, t_1,t_2)_{L_p(\T^2)} \Big]^q
\frac{\dd t_1}{t_1} \frac{\dd t_2}{t_2} \right)^{{1}/{q}}.
\end{equation}
Indeed, for $f_0$ we have
$$\omega_{\alpha_1,\alpha_2}(f_0,t_1,t_2)_{L_p(\T^2)}\asymp t_1^{\alpha_1} t_2^{\alpha_2}
\left(\ln{\frac{2}{t_1}} \right)^{\beta_1+\frac12}
\left(\ln{\frac{2}{t_2}} \right)^{\beta_2+\frac12}
$$
and therefore for this function
$$
  \left\{ \int\limits^{\delta_1}_{0} \int\limits^{\delta_2}_{0}
\Big[(t_1t_2)^{-\theta}
\omega_{\alpha_1,\alpha_2}(f_0,t_1,t_2)_{L_p(\T^2)}  \Big]^q \frac{\dd t_1}{t_1}\frac{\dd t_2}{t_2}
\right\} ^{1/q} \asymp \delta_1^{\alpha_1-\theta}
\delta_2^{\alpha_2-\theta} \left(\ln{\frac{2}{\delta_1}}
\right)^{\beta_1+\frac12} \left(\ln{\frac{2}{\delta_2}} \right)^{\beta_2+\frac12}.
$$
Thus, for the function $f_0$
the integrals in the right-hand sides of (\ref{f10.1}) and (\ref{f10.2})
have different orders of magnitude as functions of
$\delta_1$ and $\delta_2$ (cf. (\ref{f10.1111}) and (\ref{f10.2})).

%%%%%%%%%%%%%%%%%%%%%%%%%%%%%%%%%%%%%%%%%%%%%%%%%%%%%%%%%%%%%%%%%%%%%%
%%%%%%%%%%%%%%%%%%%%%% YOUR MATHEMATICS  %%%%%%%%%%%%%%%%%%%%%%%%%%%%%
%%%%%%%%%%%%%%%%%%%%%%%%%%%%%%%%%%%%%%%%%%%%%%%%%%%%%%%%%%%%%%%%%%%%%%

{

\bigskip
\hskip1.4 em\vbox{\noindent  M.~Potapov\\
Department of Mechanics and Mathematics\\
Moscow State University
\\Moscow  119991
 Russia\\
 {\tt mkpotapov@mail.ru }}

\bigskip
% first author's name, affiliation, and address, including email and web page
\hskip1.4 em\vbox{\noindent  B.~Simonov\\
Department of Applied Mathematics
\\Volgograd State Technical University,
\\Volgograd 400005 Russia\\
 {\tt simonov-b2002@yandex.ru }}

\bigskip
\hskip1.4 em\vbox{\noindent  S.~Tikhonov\\
ICREA and Centre de Recerca Matem\`{a}tica
\\Apartat 50, Bellaterra, Barce\-lona 08193 Spain\\
 {\tt stikhonov@crm.cat }
 \\{\tt http://www.icrea.cat/Web/ScientificStaff/Sergey-Tikhonov-479}}

}

%%%%%%%%%%%%%%%%%%%%%%%%%%%%%%%%%%%%%%%%%%%%%%%%%%%%%%%%%%%%%%%%%%%%%%
\endddoc
%%%%%%%%%%%%%%%%%%%%%%%%%%%%%%%%%%%%%%%%%%%%%%%%%%%%%%%%%%%%%%%%%%%%%%
\vskip 1.5cm



\begin{thebibliography}{100}

\bibitem{Am}
T. I. Amanov, Spaces of Differentiable Functions With Dominating Mixed Derivatives.
Nauka Kaz. SSR, Alma-Ata, 1976.

\bibitem{ant} A. P. Antonov,
Smoothness of sums of trigonometric series
with monotone coefficients, Russian Mathematics (Iz. VUZ) {\bf 51 (4)} (2007), 18--26.


\bibitem{ba}
K. I. Babenko, Approximation by trigonometric polynomials in a certain class of
periodic functions of several variables. Dokl. Akad. Nauk SSSR  {\bf 132}  (1960), 982--
985; English transl. in Soviet Math. Dokl.  {\bf 1}  (1960), 672--675.


\bibitem{2}
 N. S. Bakhvalov,  {Embedding theorems for
    classes of functions with several bounded derivatives},
    Vestn. Mosk. Univ., Ser. I, \textbf{18 (4)} (1963), 7--16.

\bibitem{bazar}
D. B. Bazarkhanov, Characterizations of the Nikol'skij-Besov and
Lizorkin-Triebel Function Spaces of Mixed Smoothness, Proc. Steklov
Inst.  {\bf  243} (2003), 46--58; translation from Function spaces,
approximations, and differential equations, Collected papers.
Dedicated to the 70th birthday of Oleg Vladimirovich Besov,
corresponding member of RAS, Tr. Mat. Inst. Steklova,  {\bf 243}
(2003), 53--65.


\def\JAT{J. Approx.\ Theory}
\def\CA{Constr.\ Approx.}

\bibitem{ben-sha}
C. Bennett, R. Sharpley, Interpolation of Operators, Academic Press, Boston, 1988.

\bibitem{berisha}
M. Berisha,
On the coefficients of double lacunary trigonometric series,
Serdica,  {\bf 14(1)} (1988), 68--74.

\bibitem{be}
O. V. Besov, On some conditions for derivatives of periodic functions to
belong to $L_p$, Nauchn. Dokl. Vyssh. Shkoly Fiz.-Mat. Nauki, {\bf 1} (1959), 13--17
(in Russian).

\bibitem{besov}
O. V. Besov, V. P. Il'in, S. M. Nikol'skii, Integral Representations of Functions and Imbedding Theorems,
J. Wiley and Sons, New York, 1978, 1979; translated from Russian: Nauka, Moscow, 1975.

\bibitem{brud}
A. Brudnyi, Y. Brudnyi, Methods of Geometric Analysis in Extension and Trace Problems, Vol. 1, Springer, 2012.

\bibitem{butzer}
P. L. Butzer, H. Dyckhoff, E. G\"{o}rlich, R. L. Stens, Best trigonometric approximation,
fractional order derivatives and Lipschitz classes, Canad. J. Math.,  {\bf 29(4)} (1977), 781--793.

\bibitem{chen}
W. Chen, Z. Ditzian, Mixed and directional derivatives, Proc. Amer. Math. Soc.  {\bf 108} (1990), 177--185.

\bibitem{13}
 C. Cottin, Mixed K-functionals: a measure of smoothness for
blending-type approximation, Math. Z.,  {\bf 204} (1990), 69--83.

\bibitem{jat-dai}
F. Dai, Z. Ditzian, S. Tikhonov,
{Sharp Jackson inequality}, {\JAT},   {\bf 151(1)} (2008), 86--112.

\bibitem{dai}
F. Dai, Y. Xu,
Analysis on h-harmonics and Dunkl Transforms, Birkh\"{a}user Verlag,
to appear.

\bibitem{dah}
W. Dahmen, R. DeVore, K. Scherer,
Multi-dimensional spline approximation, SIAM J. Numer. Anal., {\bf 17(3)} (1980), 380--402.

\bibitem{dav}
O. V. Davydov, Sequences of rectangular Fourier sums of continuous functions with given majorants of the mixed moduli of smoothness,
Sb. Math.,
{\bf 187(7)} (1996), 981--1004.

\bibitem{de1}
R.~A. DeVore, S.V. Konyagin, V.N. Temlyakov, Hyperbolic wavelet approximation, \CA,
{\bf 14} (1998), 1--26.

\bibitem{de-lo}
R.~A. DeVore, G. G. Lorentz, Constructive Approximation, Springer, 1993.

\bibitem{de2}
R.~A. DeVore, P.P. Petrushev, V.N. Temlyakov, Multivariate trigonometric polynomial approximations
with frequencies from the hyperbolic cross, Math. Notes,
{\bf 56}  (1994), 900--918.

\bibitem{dit}
Z. Ditzian,
Moduli of continuity in $\R^n$ and $D\subset \R^n$,
Trans. Amer. Math. Soc.,
{\bf 282(2)}  (1984), 611--623.

\bibitem{march}
Z. Ditzian, On the Marchaud-type inequality, Proc. Amer. Math. Soc.,
{\bf 103}  (1988), 198--202.

\bibitem{ima}
Z. Ditzian,
The modulus of smoothness and discrete data in a square domain,
IMA Journal of Numerical Analysis, {\bf 8}  (1988), 311--319.


\bibitem{ditzian}
Z. Ditzian, V. H. Hristov, K. G. Ivanov, Moduli of smoothness and K-functionals in
$L_p$, $0 < p < 1$, \CA,
{\bf 11(1)}  (1995), 67--83.

\bibitem{dit-ul}
Z. Ditzian, S. Tikhonov,
{Ul'yanov and Nikol'skii-type inequalities}, {\JAT,
{\bf 133(1)}  (2005), 100---133.}

\bibitem{dit-studia}
Z. Ditzian, S. Tikhonov,
{Moduli of smoothness of functions and their derivatives},
{Studia Math.}, {\bf 180(2)} (2007), 143--160.

\bibitem{dung}
D. Dung, T. Ullrich,
Whitney type inequalities for local anisotropic polynomial approximation,
\JAT, {\bf 163(11)}  (2011), 1590--1605.

\bibitem{dyach1996}
M. I. Dyachenko, Uniform convergence of double Fourier series for classes of functions with anisotropic smoothness, Mat. Zametki,
{\bf 59(6)} (1996), 937--943.

\bibitem{4}
 M. I. Dyachenko, Some problems in the theory of multiple
trigonometric series, Russian Math. Surveys
{\bf 47(5)} (1992), 103--171; translated from Uspekhi Mat. Nauk
{\bf 47(5)} (1992), 97--162.

\bibitem{5}
 M.~Dyachenko, S.~Tikhonov. {A Hardy-Littlewood theorem for
multiple series}, {J. Math. Anal. Appl.}, {\bf 339} (2008), 503--510.

\bibitem{23}
M. G. Esmaganbetov, Existence conditions of the mixed Weyl derivatives in
 $L_p([0,2\pi]^2) (1<p<\infty)$ and their structural properties. Deposited at the VINITI, manuscript N 1675-82, 17.02.1982.



\bibitem{ho2}
G. Garrig\'{o}s, R. Hochmuth, A. Tabacco,
Wavelet characterizations for anisotropic Besov spaces with $0<p<1$,
Proc. Edinb. Math. Soc., II. Ser.
{\bf 47(3)} (2004), 573--595.


\bibitem{gorb}
D. Gorbachev, S. Tikhonov, {Moduli of smoothness and growth properties of Fourier transforms: two-sided estimates},
{\JAT},
{\bf 164(9)}
 (2012), 1283--1312.



\bibitem{ha}
M. Hansen, Nonlinear approximation and function spaces of dominating mixed smoothness, Ph.D. thesis, FSU, Jena, Germany, 2010.


\bibitem{han}
M. Hansen, W. Sickel,
Best $m$-Term Approximation and Sobolev–Besov Spaces of dominating mixed smoothness --- the Case of Compact Embeddings,
\CA,
{\bf 36(1)} (2012), 1--51.

\bibitem{ho1}
R. Hochmuth,
Wavelet characterizations for anisotropic Besov spaces,
Appl. Comput. Harmon. Anal.,
{\bf 12(2)} (2002), 179--208.

\bibitem{johnen}
H. Johnen, K. Scherer, On the equivalence of the K-functional and moduli of
continuity and some applications, Constructive theory of functions of several variables
(Proc. Conf., Math. Res. Inst., Oberwolfach 1976), Lecture Notes in Math.,
{\bf 571},
Springer-Verlag, Berlin–Heidelberg 1977,  119--140.


\bibitem{teml0}
B. S. Kashin, V. N. Temlyakov, On a norm and approximate characteristics of classes of multivariable functions, Theory of functions, CMFD, 25, PFUR, M., (2007), 58--79; English version:
Journal of Mathematical Sciences,
{\bf 155(1)} (2008), 57--80.


\bibitem{kolyada}
V. Kolyada, F. P\'{e}rez L\'{a}zaro,
Inequalities for partial moduli of continuity and partial derivatives,
\CA,
{\bf 34(1)} (2011),
23--59.

\bibitem
{krbec}
M. Krbec, H.-J. Schmeisser,
Imbeddings of Brezis–Wainger type. The case of missing derivatives,
Proceedings of the Royal Society of Edinburgh, Section: A Mathematics,
{\bf 131(3)} (2001),  667--700.

\bibitem{8}
L. Leindler, Generalization of inequalities of Hardy and Littlewood, Acta Sci. Math., {\bf 31} (1970),
279--285.


\bibitem
{li1}
P. I. Lizorkin,  Properties of functions in the spaces $\Lambda ^r_{p\,\theta }$, Trudy Mat. Inst. Steklov
{\bf 131}  (1974), 158--181; English transl. in Proc. Steklov Inst. Math. {\bf 131} (1974), 165--188.

\bibitem
{li2}
P. I. Lizorkin,  S. M. Nikol'skii, Classification of
differentiable functions on the basis of spaces with dominating
mixed smoothness, Trudy Mat. Inst. Steklov, {\bf 77} (1965),
143--167.

\bibitem
{li3}
P. I. Lizorkin,  S. M. Nikol'skii, Function spaces of mixed smoothness from the decomposition point of view,
Trudy Mat. Inst. Steklov, {\bf 187} (1989), 143--161; English transl. in Proc. Steklov Inst. Math., {\bf 187(3)} (1990), 163--184.

\bibitem
{ma}
J. Marcinkiewicz, Sur quelques int\'{e}grales du type de Dini, Ann. Soc. Polon.
Math., {\bf 17} (1938), 42--50.

\bibitem{15}
 F. Moricz, On double cosine, sine, and Walsh series with
monotone coefficients,
 Proc. Amer. Math. Soc., {\bf 109}  (1990), 417--425.


\bibitem{veres}
F. Moricz, A. Veres,
On the absolute convergence of multiple Fourier series
Acta Math. Hungar.,
{\bf 117(3)}  (2007), 275--292.

\bibitem{mus}
J. Musielak, On the absolute convergence of multiple Fourier series, Ann. Polon. Math., {\bf 5}  (1958), 107--120.

\bibitem{1}
S. M. Nikol'skii, {Functions with dominating mixed derivatives satisfying multiple H\"{o}lder
conditions},  Sib. Mat. Zh., \textbf{4(6)} (1963),  1342--1364; English transl. in
Amer. Math. Soc. Transl., Ser. 2, {\bf 102} (1973), 27--51].

\bibitem{1111}
S. M. Nikol'skii,
Stable boundary-value problems of a differentiable function of several variables,
Mat. Sb.  {\bf 61(103)}, N 2 (1963), 224--252.

\bibitem{7}
S. M. Nikol'skii, { Approximation of Functions of Several Variables and Embedding Theorems}.
Springer-Verlag, 1975.

\bibitem{17}
E. D. Nursultanov, On the coefficients of multiple Fourier
series from $L\sb p$-spaces, Izv. Ross. Akad. Nauk Ser. Mat., {\bf 64 (1)}
(2000), 95--122; translation in Izv. Math., {\bf 64 (1)} (2000), 93--120.

\bibitem{6}
M. K. Potapov, On ``angular" approximation, Proc. Conf. Constructive Function Theory
(Budapest, 1969), Akad. Kiado, Budapest, 1971, 371--379.

\bibitem{9}
M. K. Potapov, Imbedding of classes of functions with a dominating mixed modulus of smoothness,
 Trudy Mat. Inst. Steklov., {\bf 131} (1974), 199--210.

\bibitem{25}
M. K. Potapov,
The Hardy-Littlewood and Marcinkiewicz-Littlewood-Paley theorems, angular approximation,
and embedding of some classes of functions, Mathematica (Cluj), {\bf 14(37)}, 2, (1972), 339--362.

\bibitem{21}
M. K. Potapov, B.V. Simonov, B. Lakovich,
On estimates for the mixed modulus of continuity of a function with a transformed Fourier series,
Publ. Inst. Math., Nouv. S\'{e}r., {\bf  58(72)} (1995), 167--192.


\bibitem{22}
M. K. Potapov, B. V. Simonov,
Generalized classes of the Besov-Nikol'skij and Weyl-Nikol'skij functions: Their interrelations,
Proc. Steklov Inst. Math. {\bf 214} (1996), 243--259; translation from Tr. Mat. Inst. Steklova {\bf 214} (1997), 250--266.




%\bibitem{24}
%M.K. Potapov; B. Lakovich; B.V. Simonov,
%On the interrelation of the Besov-Nikol'skii and Weyl-Nikol'skii function classes in the mixed norm,
% Mathematica Montisnigri, 2000, vol. XII,  63-85.


\bibitem{deriv-st}
M. K. Potapov, B. V. Simonov, S. Yu. Tikhonov,
On Besov, Besov-Nikol'skii classes and on the estimates of mixed smoothness moduli of fractional
derivatives, Proc. Steklov Inst. Math., {\bf 243} (2003), 234--246; translated from
Tr. Mat. Inst. Steklova, {\bf 243} (2003), 244--256.


\bibitem{deriv-stt}
M. K. Potapov, B. V. Simonov, S. Yu. Tikhonov,
Transformation of Fourier series using power and weakly oscillating sequences, Mat.
Zametki, {\bf 77(1)}  (2005), 99--116; translation in Math. Notes, {\bf 77(1)} (2005), 90--107.

\bibitem{sim-2}
M. K. Potapov, B. V. Simonov, S. Yu. Tikhonov, {Relations between moduli of smoothness in different
metrics}, Vestn. Mosk. Univ., Ser. 1: Mat., Mekh., No. 3 (2009), 17--25; translation in Moscow Univ. Math. Bull. {\bf 64} (2009), 105--112.

\bibitem{26}
M. K. Potapov, B.V. Simonov, S. Yu. Tikhonov, {
Relations between the mixed moduli of smoothness and
embedding theorems for Nikol'skii classes},
Proceedings of the Steklov Institute of Mathematics, {\bf 269} (2010),
197--207; translation from Russian: Trudy Matem. Inst. V. A. Steklova, {\bf 269} (2010), 204--214.


\bibitem{isaac}
M. K. Potapov, B.V. Simonov, S. Yu. Tikhonov, {
Constructive characteristics of mixed moduli of smoothness of positive orders},
to be published in Proceedings of VII ISAAC Congress, Moscow, 2011.

\bibitem{pust}
N. N. Pustovoitov, The orthoprojection widths of some classes of periodic functions of two variables with
a given majorant of the mixed moduli of continuity, Izvestiya: Mathematics, {\bf 64(1)} (2000), 121--141.

\bibitem{ram}
A. K. Ramazanov, A. R. K. Ramazanov,
Mixed moduli of continuity and their applications in polynomial approximations with interpolation,
Analysis Mathematica,
{\bf 35 (3)} (2009), 213--232.

\bibitem{12}
K. V. Runovskii, {Several questions of approximation theory}, Disser. Cand. Nauk, Moscow, MGU,
1989.

\bibitem{runov}
K. V. Runovskii,
A direct theorem on approximation ``by angle'' in the spaces $L_p, 0 < p < 1$,
Math. Notes, {\bf 52 (5)} (1992), 1140--1142.

\bibitem
{samko}
S. G. Samko, A. A. Kilbas, O. I. Marichev, Fractional Integrals and Derivatives:
Theory and Applications. Gordon and Breach Science Publishers, Yverdon,
1993.

\bibitem{sch}
H.-J. Schmeisser,  On spaces of functions and distributions with mixed smoothness properties of Besov–Triebel–
Lizorkin type. I, II, Math. Nachr.
{\bf 98} (1980), 233--250;
{\bf 106} (1982), 187--200.


\bibitem{sch1}
H.-J. Schmeisser,
Recent developments in the theory of function spaces with dominating mixed smoothness, Nonlinear Analysis,
Function Spaces and Applications, Proceedings of the Spring School held in Prague, May 30-June 6, 2006. {\bf 8}.
Czech Academy of Sciences, Mathematical Institute, Praha (2007), 145--204.


\bibitem{sch3}
H.-J. Schmeisser, W. Sickel,
Spaces of functions of mixed smoothness and
approximation from hyperbolic crosses,
\JAT, {\bf 128} (2004). 115--150.

\bibitem{book1}
H.-J. Schmeisser, H. Triebel. Topics in Fourier Analysis and Function Spaces, Wiley,
Chichester, 1987.


\bibitem{dah1}
L. Schumaker,
Spline Functions: Basic Theory, Cambridge Univ. Press, 2007.

\bibitem{19}
B. Simonov, I. Simonova, On the Fourier coefficients of some functions. Deposited at the VINITI, manuscript N 1675-82, 17.04.1982.

\bibitem{sbornik} B. V. Simonov,  S. Yu. Tikhonov, { Embedding
theorems in constructive approximation}, Sb. Math., {\bf 199(9)} (2008),  1367--1407; translation from Mat. Sb.,
{\bf 199(9)} (2008),  107--148.

\bibitem{sim-jat}
B. Simonov, S. Tikhonov, {Sharp Ul'yanov-type inequalities using fractional smoothness}, {\JAT},
{\bf 162(9)} (2010), 1654--1684.

\bibitem{stro}
J. O. Str\"{o}mberg,
Computation with wavelets in higher dimensions. Proceedings
of the International Congress of Mathematicians, Vol. III (Berlin, 1998). Doc. Math.
1998, Extra Vol. III, 523--532.

\bibitem{sun}
Y. Sun,  H. Wang,
Representation and approximation of multivariate periodic functions with bounded mixed moduli of smoothness,
Proc. Steklov Inst. Math.,
{\bf 219} (1997),
 350--371; translation from Tr. Mat. Inst. Steklova {\bf 219} (1997), 356--377.


\bibitem{10}
{R. Taberski,} { Differences, moduli and derivatives of
fractional orders},  Comment. Math. Prace Mat.
{\bf 19(2)} (1976/77), 389--400.


\bibitem{teml1}
{V. N. Temlyakov,}
Approximation of functions with a bounded mixed difference by trigonometric polynomials, and the widths of some classes of functions,
Mathematics of the USSR-Izvestiya,  {\bf 20(1)} (1983), 173--187.

\bibitem{teml2}
V. N. Temlyakov, Approximations of functions with bounded mixed derivative, Trudy Mat. Inst. Steklov., {\bf 178} (1986), 3--113.

\bibitem{teml3}
V. N. Temlyakov,
Approximation of Periodic Functions,  Nova Science Pub. Inc. 1994.

\bibitem{real1}
 S. Tikhonov, {On modulus of smoothness of fractional order},
Real. Analysis Exchange, {\bf 30(2)} (2004/2005), 507--518.


\bibitem{16}
 S. Tikhonov, {Trigonometric series with general monotone
coefficients}, J. Math. Anal. Appl., {\bf 326(1)} (2007),
721--735.

\bibitem{real3}
S. Tikhonov, Best approximation and moduli of smoothness: computation and
equivalence theorems, {\JAT}, {\bf 153} (2008), 19--39.


\bibitem{jfaa}
S. Tikhonov, {Weak type inequalities for moduli of smoothness: the case of limit value parameters},
J. Fourier Anal. Appl., {\bf 16 (4)} (2010), 590--608.


\bibitem{timan-book}
A. F. Timan, Theory of Approximation of Functions of a Real Variable,
Pergamon Press: Oxford, 1963.


\bibitem
{timan}
M. F. Timan, On Jackson's theorem in $L_p$ spaces, Ukrain. Mat. Zh. {\bf 18 (1)} (1966), 134--137 (in Russian).


\bibitem
{timan2}
M. F. Timan,
Difference properties of functions of several variables,
Math. USSR - Izvestiya, {\bf 3(3)} (1969), 633--642; translation from
Izv. Akad. Nauk SSSR, Ser. Mat. {\bf 33(3)} (1969),
667--676.

\bibitem
{tom}
M. Tomi\'{c},
The converse theorem of approximation by angle in various metrics for non-periodic functions,
Facta Univ., Ser. Math. Inf., {\bf 16} (2001), 45--60.

\bibitem
{tom1}
M. Tomi\'{c},
On representation of derivatives of functions in $L_p$,
 Matematicki Vesnik,
{\bf 62(3)} (2010), 235--250.


\bibitem{trebels}
W. Trebels, {Inequalities for moduli of smoothness versus embeddings of function spaces},
Arch. Math., {\bf 94} (2010), 155--164.

\bibitem{11}
W. Trebels, On the approximation behavior of the Riesz means
in $L^p(\R^n)$, in: Approximation theory (Proc. Internat. Colloq., Inst. Angew. Math. Univ. Bonn, Bonn, 1976),
R. Schaback and K. Scherer (eds.), Lect. Notes Math. 556, 428--438,  Springer, Berlin, 1976.


\bibitem{book2}
H. Triebel, Bases in Function Spaces, Sampling, Discrepancy, Numerical Integration, EMS Tracts in Mathematics, 2010.

\bibitem{tri1}
H. Triebel,  A diagonal embedding theorem for function spaces with dominating mixed
smoothness properties, Banach Center Publications {\bf 22}, Warsaw, (1989), 475--486.

\bibitem{trigub}
R. M. Trigub, E. S. Belinsky,
Fourier Analysis and Approximation of Functions,
 Springer, 2010.

\bibitem{ullrich}
T. Ullrich, Function spaces with dominating mixed smoothness, characterization by differences,
Jenaer Schriften zur Mathematik und Informatik, Math/Inf/05/06 (2006).

\bibitem{14}
P. L. Ul'yanov, The imbedding of certain function classes $H^{\omega}_{p}$,
Izv. Akad. Nauk SSSR Ser. Mat., {\bf 32(3)} (1968),  649--686.

\bibitem{vyb}
J. Vyb\'{i}ral, Function spaces with dominating mixed smoothness, Dissertationes Math., {\bf 436} (2006).

\bibitem{vyb1}
J. Vyb\'{i}ral, W. Sickel,
Traces of functions with a dominating mixed derivative in $\R^3$,
Czech. Math. J., {\bf 57(4)} (2007), 1239--1273.


\bibitem{wang}
H. Wang,
Widths between the anisotropic spaces and the spaces of functions with mixed smoothness, \JAT, {\bf 164(3)} (2012), 406--430.

\bibitem{wilmes}
G. Wilmes, On Riesz-type inequalities and K-functionals related to Riesz potentials in $\R^n$ , Numer. Funct. Anal.
Optimization, {\bf 1} (1979), 57--77.

\bibitem{w1}
G. Wilmes, Some inequalities for Riesz potentials of trigonometric polynomials of several variables, in: Harmonic
Analysis in Euclidean Spaces, (Proc. Sympos. Pure Math., Williamstown, 1978), Part I, Amer. Math. Soc.,
Providence, 1979, 175--182.


\bibitem
{timan1}
I. E. Zhak, M. F. Timan, On summation of double series, Mat. Sb. (N.S.), {\bf 35(77)}, 1 (1954), 21--56.

\bibitem{zhizhi-s}
L. V. Zhizhiashvili,
Some problems in the theory of simple and multiple trigonometric and orthogonal series, Russian Mathematical Surveys,
{\bf 28 (2)(170)} 1973, 65--119.

\bibitem{zhizhi}
L. V. Zhizhiashvili,
Some Problems of Multidimensional Harmonic Analysis, Tbilisi Univ. Press, Tbilisi, 1996.

\bibitem{zhizhi-b}
L. V. Zhizhiashvili,
Trigonometric Fourier series and Their Conjugates, Kluwer Acad. Publ., 1996.

\bibitem{18}
  A. Zygmund, { Trigonometric Series}.
    Cambridge Univ. Press, Cambridge, 1959.

\end{thebibliography}
